\documentclass{article}
\usepackage[utf8]{inputenc}
\usepackage[T1]{fontenc} 
\usepackage{amsmath}
\usepackage{amsfonts}
\usepackage{amssymb}
\usepackage{amsthm}
\usepackage{color}
\usepackage{esint}
\usepackage{enumerate}
\usepackage{dsfont}
\usepackage{stmaryrd}
\usepackage{hyperref}
\usepackage{graphicx}
\usepackage{graphics}
\usepackage{epstopdf}
\usepackage{epsfig}
\usepackage{setspace}
\usepackage{comment} 
\usepackage{tikz}
\usetikzlibrary{positioning,arrows}
\usetikzlibrary{calc,arrows}
\usetikzlibrary{er,positioning,bayesnet}

\usepackage{pgfplots}
\usepackage{pgfplotstable}
\pgfplotsset{compat=1.11} 

\newtheorem{theorem}{Theorem}
\newtheorem{lemma}[theorem]{Lemma}
\newtheorem{remark}[theorem]{Remark}
\newtheorem{proposition}[theorem]{Proposition}
\newtheorem{definition}[theorem]{Definition}

\newcommand{\RR}{\mathbb{R}}
\newcommand{\R}{\mathbb{R}}

\newcommand{\Z}{\mathbb{Z}}

\newcommand{\N}{\mathbb{N}}

\newcommand{\eps}{\varepsilon}

\newcommand{\dps}{\displaystyle}
\newcommand{\dis}{\displaystyle}

\renewcommand{\div}{\operatorname{div}}

\begin{document}

\date{\today}
\author{Frédéric Legoll$^1$\footnote{corresponding author}\ , Pierre-Lo\"ik Rothé$^1$, Claude Le Bris$^1$ and Ulrich Hetmaniuk$^2$
\\
{\footnotesize $^1$ \'Ecole des Ponts \& Inria, 6 et 8 avenue Blaise Pascal, 77455 Marne-La-Vall\'ee Cedex 2, France}
\\  
{\footnotesize $^2$ Shift-Invert LLC, Colorado Springs, CO 80919, USA}
\\
{\footnotesize \tt \{frederic.legoll,claude.le-bris\}@enpc.fr, ulrich@shiftinvert.com}\\
}
\title{An MsFEM approach enriched using Legendre polynomials}

\maketitle


\begin{abstract}
We consider a variant of the conventional MsFEM approach with enrichments based on Legendre polynomials, both in the bulk of mesh elements and on their interfaces. A convergence analysis of the approach is presented. Residue-type {\em a posteriori} error estimates are also established. Numerical experiments show a significant reduction in the error at a limited additional off-line cost. In particular, the approach developed here is less prone to resonance errors in the regime where the coarse mesh size $H$ is of the order of the small scale $\eps$ of the oscillations.   
\end{abstract}

\section{Introduction}

We consider the problem
\begin{equation}
\label{eq:lap-coeff}
-\div (A\nabla u) = f \ \ \text{in $\Omega$}, \qquad u = 0 \ \ \text{on $\partial \Omega$},
\end{equation}
where $\Omega$ is a bounded polygonal domain in~$\mathbb{R}^2$, $f$ is a given right-hand side and the symmetric matrix-valued elliptic coefficient $A(x)$ presents heterogeneities at small scales (henceforth denoted by $\eps$) compared with the characteristic size of $\Omega$. Classical approximation techniques such as finite elements are known to poorly perform in such cases, unless the mesh size is taken (possibly prohibitively) small. Multiple alternative dedicated approaches have therefore been introduced. Among those, the multiscale finite element method (henceforth abbreviated as MsFEM), introduced in~\cite{Hou_book_msfem_2009,Hou:1997:A-maa}, uses a Galerkin approach of~\eqref{eq:lap-coeff} on a pre-computed basis. The basis functions are obtained by solving \emph{local} problems mimicking~\eqref{eq:lap-coeff} at the scale of mesh elements, with carefully chosen right-hand sides and boundary conditions. The vanilla version of the approach, called \emph{linear MsFEM}, uses as basis functions the solutions to these local problems, posed on each mesh element, with null right-hand sides and with the coarse P1 elements as Dirichlet boundary conditions (see~\eqref{eq:lin_msfem} below for the precise definition of these basis functions). Various improvements of that version are possible. In particular, the so-called \emph{oversampling} variant, which solves local problems on larger domains and restricts their solutions to the considered element, is very effective. The down side is that the approach is not conformal and the size of the oversampling area must be carefully calibrated, which can be a delicate practical issue. 

\medskip

Our purpose here is to introduce and study an MsFEM method improved differently. It essentially elaborates upon the \emph{Approximate Component Mode Synthesis} (ACMS) method introduced in~\cite{hetmaniuk2010special} and fully analyzed in~\cite{hetmaniuk2014error}. In that approach, the linear MsFEM basis is enriched with local eigenvectors related to the scalar product associated with the variational formulation of~\eqref{eq:lap-coeff}. The approach of~\cite{hetmaniuk2014error,hetmaniuk2010special}, which is related to domain decomposition methods for elliptic partial differential equations~\cite{gander_loneland,gander2015analysis,gervasio_spectral_1997,quva:99,Toselli2005} and to modern component mode synthesis methods~\cite{Bennighof:2004:An-aa,bour:92}, is very effective. The resolution of eigenproblems for each element of the coarse discretization can however be computationally challenging, even considering that the work is performed off-line. This is the reason why the approach we present here complements the linear MsFEM basis with enrichments that are not eigenvectors, but solutions of edge and bulk problems using polynomials either as boundary condition or right-hand side (in Section~\ref{sec:pres_dis_method}, see e.g.~\eqref{eq:def_phi_B} for the definition of some of the basis functions of our MsFEM approach and~\eqref{eq:def_VBHM}--\eqref{eq:approx_bubble} for the global Galerkin approximation). Similarly to the other MsFEM variants, all basis functions for such enrichments can be computed in parallel. One cannot indeed too much emphasize that, if the dogma of multiscale approaches is to drastically reduce the on-line cost at the expense of an increase of the off-line cost, it might be the case for a large class of complex enough problems that the approach is doomed because of a prohibitively computationally expensive off-line stage. Another advantage of the approach presented here is that the classical Legendre interpolation results apply, allowing one to get rigorous \emph{a priori} and \emph{a posteriori} error estimates for the approach more easily. A similar approach has been introduced independently in~\cite{gao2018high}, for the specific case of quadrangles, Legendre polynomials and Gauss-Lobatto quadratures. The approach of~\cite{gao2018high} shows promising results in time-domain acoustic-wave modeling: it yields approximations that compare well with reference solutions computed with the spectral finite element method. Our aim is to push further the approach by expanding it to triangular meshes and to provide a detailed convergence analysis (see Proposition~\ref{lem:Error_tot}), along with some theoretical tools for adaptivity based on suitable {\em a posteriori} error estimates (see Proposition~\ref{prop:error_indicator}). We emphasize that our method is both \emph{local} and \emph{conformal}: the support of the enrichment function is either the two elements associated with the edge when an edge element is considered (see~\eqref{eq:def_phi_gamma}), or the element itself when a bulk element is considered (see~\eqref{eq:def_phi_B}). Also, as said above, it is \emph{fully parallel} in the off-line stage. We mention that another, very interesting and efficient, line of thought is exemplified by the approach called \emph{Localized Orthogonal Decomposition method} (LOD) introduced in~\cite{maalqvist2014localization}. There, the classical finite elements are enriched with solutions to specific PDEs. These functions have global supports, in contrast to MsFEM basis functions. They however turn out to rapidly decay away from the element considered. This property allows one to design an approximation space with functions solution to PDEs with smaller, truncated supports (typically of size of order $O(H \, | \ln H |)$). The associated error estimates in the energy norm are then independent of the scales of the heterogeneities.

\medskip

Our article is organized as follows. The enriched MsFEM variant that we introduce is presented in Section~\ref{sec:pres_dis_method}. We prove (in Section~\ref{sec:main_res}) that, 
with a sufficiently large number of enrichment functions, we can get a convergence rate that does not depend on the oscillations of $ A$. The numerical experiments we present in Section~\ref{sec:Num_exp} moreover show that already a small number of enrichment functions significantly reduces the error. Our analysis applies to both quadrangular and triangular meshes, the latter being more flexible and allowing one to discretize more complex geometries than those accessible to quadrangular meshes. Furthermore, we propose an \emph{a posteriori} estimator that can be used to locally adapt the level of enrichment.

The numerical experiments of Section~\ref{sec:Num_exp} show that the proposed approach outperforms the linear MsFEM especially in the regime where $H \approx \varepsilon$, allowing for results of comparable quality to those obtained using the ACMS method of~\cite{hetmaniuk2014error,hetmaniuk2010special}, is on par with non-conformal approaches such as the variant of MsFEM using oversampling, and that it achieves all this at a reasonable additional computational cost. Our numerical results also seemingly indicate that the \emph{a posteriori} estimator we propose reproduces truly the trend of the error in energy norm. 

\section{Discretization approach}
\label{sec:pres_dis_method}

We define a family of meshes $\left(\mathcal{T}_H\right)_H$ of the two-dimensional domain $\Omega$, i.e. a decomposition of $\Omega$ into a finite number of convex elements (quadrangles or triangles) with straight edges. Note that, throughout the article, we work in two dimensions, both for the analysis and for the numerical tests. Our approach and our analysis can presumably be extended to some three-dimensional meshes (such as cartesian meshes), but we will not proceed in this direction here.

The mesh is assumed {\em conformal} (there is no hanging nodes and each internal edge is shared by exactly two elements of the mesh) and {\em regular} in the following sense:
\begin{equation} \label{eq:def_mesh_regular}
  \begin{array}{c}
    \text{for any element $K$, there exists an affine transformation $F : K_{\rm ref} \mapsto K$},
    \\
    \text{where $K_{\rm ref}$ is the reference element (here the reference square or triangle),}
    \\
    \text{such that} \qquad \| \nabla F \|_{L^\infty} \leq \gamma \, H \qquad \text{and} \qquad \| \nabla F^{-1} \|_{L^\infty} \leq \gamma \, H^{-1},
    \\
    \text{where $\gamma>1$ is a constant independent of both $K$ and $H$.}
  \end{array}
\end{equation}
In practice, this property is ensured using a mesh with quadrangular (or triangular) elements with a minimum angle condition (see e.g.~\cite[Section~4.4]{BS}). We denote by $\Gamma$ the interior skeleton, that is
$$
\Gamma = \left( \cup_{K \in \mathcal{T}_H} \partial K \right) \setminus \partial \Omega.
$$

The variational formulation of~\eqref{eq:lap-coeff} is expressed using, for $u,v \in H^1_0(\Omega)$, the bilinear form~$\dps a(u,v) = \int_\Omega (\nabla v)^T A \nabla u$. The associated energy norm is denoted by $\dps \|v\|_E = \sqrt{a(v,v)}$. Since $A$ is assumed symmetric, the unique solution $u$ to~\eqref{eq:lap-coeff} also satisfies
$$
u = \underset{v\in H_0^1(\Omega)}{\mathrm{argmin}} \ \left(\frac{1}{2} a(v,v) - \langle f,v \rangle_{L^2(\Omega)}\right).
$$
We introduce the set of bubble functions
$$
V_B = \left\{ v \in H^1_0(\Omega), \ \ v|_K \in H^1_0(K) \ \text{for any $K \in \mathcal{T}_H$} \right\},
$$
where the subscript $B$ stands, understandably, for \emph{bubbles}. We also define
\begin{equation} \label{eq:def-VGamma}
V_\Gamma = \left\{ E_\Omega \tau \in H^1_0(\Omega), \ \ \tau \in H^{1/2}_{00}(\Gamma) \right\},
\end{equation}
which is the subspace of energy-minimizing extensions of trace functions on $\Gamma$, where the extension $E_\Omega(\tau)$ solves the minimization problem $\dps \inf_{v\in H_0^1(\Omega)} a(v,v)$ subject to $v|_{\Gamma}=\tau$, that is
\begin{equation}
  \label{eq:def-harmonic}
  \begin{cases}
    -\div \left(A \nabla (E_\Omega \tau) \right) = 0 \quad \text{in $K$, for any $K \in \mathcal{T}_H$},
    \\
    E_\Omega \tau = \tau \quad \mbox{on $\Gamma$},
    \\
    E_\Omega \tau = 0 \quad \mbox{on $\partial \Omega$},
  \end{cases}
\end{equation}
in the weak sense. We recall that, in~\eqref{eq:def-VGamma}, the space $H^{1/2}_{00}(\Gamma)$ is the space of functions $\tau \in H^{1/2}(\Gamma)$ such that their extension $\overline{\tau}$ on $\Gamma \cup \partial \Omega$ defined by
$$
\overline{\tau} = \tau \ \ \text{on $\Gamma$}, \qquad \overline{\tau} = 0 \ \ \text{on $\partial \Omega$},
$$
is a function in $H^{1/2}(\Gamma \cup \partial \Omega)$. In particular, any function $\tau$ in $H^{1/2}_{00}(\Gamma)$ is such that its extension $\overline{\tau}$ belongs to $H^{1/2}(\partial K)$ for any element $K$ (thus the well-posedness of~\eqref{eq:def-harmonic}). For the convenience of the reader, we collect in Appendix~\ref{sec:H_un_demi} some more details on $H^{1/2}$ spaces.

\medskip

Both spaces $V_B$ and $V_\Gamma$ are infinite-dimensional. We readily observe that the decomposition
\begin{equation}
\label{eq:orthogonal_decomposition}
H_0^1(\Omega) = V_B \oplus V_\Gamma
\end{equation}
holds and is orthogonal with respect to the scalar product defined by $a(\cdot,\cdot)$ because of the definition of the energy-minimizing extension. Indeed, it holds that
$$
\forall v_B \in V_B, \ \ \forall v_\Gamma \in V_\Gamma, \quad a(v_B,v_\Gamma)=0,
$$
by using the variational formulation of~\eqref{eq:def-harmonic} with a test function in $H^1_0(K)$.

Following the decomposition~\eqref{eq:orthogonal_decomposition}, the solution $u$ to~\eqref{eq:lap-coeff} can be uniquely expressed as $u=u_B+u_\Gamma$ with the bubble part
$$
u_B = \underset{w\in V_B}{\mathrm{argmin}} \left(\frac{1}{2} a(w,w) - \langle f,w\rangle_{L^2(\Omega)} \right)
$$
and the interface part
$$
u_\Gamma = \underset{w\in V_\Gamma}{\mathrm{argmin}} \left( \frac{1}{2} a(w,w) - \langle f,w \rangle_{L^2(\Omega)} \right).
$$

In our approach, instead of approximating $u$ directly, we approximate $u_B$ and $u_\Gamma$ {\em separately}. This splitting is motivated as follows.

First, the decomposition~\eqref{eq:orthogonal_decomposition} implies a natural splitting of the error.
If we indeed consider a numerical approximation $u_{B,H}$ of $u_B$ in a finite dimensional space $V_{B,H} \subset V_B$ defined by
$$
u_{B,H} = \underset{w\in V_{B,H}}{\mathrm{argmin}} \left(\frac{1}{2} a(w,w) - \langle f,w\rangle_{L^2(\Omega)}\right),
$$
and likewise a numerical approximation $u_{\Gamma,H}$ of $u_\Gamma$ in some $V_{\Gamma,H} \subset V_\Gamma$, it is then natural to define our approximation of $u = u_B + u_\Gamma$ as $u_H = u_{B,H} + u_{\Gamma,H}$. We then have 
$$
u_H = \underset{w\in V_H}{\mathrm{argmin}} \left(\frac{1}{2} a(w,w) - \langle f,w\rangle_{L^2(\Omega)}\right)
$$
for the discretization space $V_H = V_{B,H} \oplus V_{\Gamma,H}$ and the error in energy norm reads as
\begin{equation} \label{eq:decoupling}
  \|u-u_H\|^2_E = \|u_B-u_{B,H}\|_E^2+\|u_\Gamma-u_{\Gamma,H}\|_E^2.
\end{equation}

Second, the analysis of the classical MsFEM suggests that the interface part $u_\Gamma$ is more difficult to approximate than the bubble part $u_B$. Even roughly approximating $u_B$ by $u_{B,H}=0$ already gives an energy error of order $O(H)$ (see~\eqref{eq:error_u_bubble0} below), which is often considered as a sufficiently small error for multiscale problems. Moreover, $u_B$ is the collection of solutions to {\em independent} local problems with homogeneous Dirichlet boundary conditions. Hence, $u_B$ can be computed effectively in parallel by using a FE solver for the Dirichlet problems.
The situation is drastically different for $u_\Gamma$. Approximating $u_\Gamma$ by $u_{\Gamma,H}=0$ yields an error that remains of order $O(1)$ and does {\em not} decay with $H$. Moreover (and this is now an argument specific to the multiscale context), when approximating $u_\Gamma$ by $u_{\Gamma,H} = u_{\rm MsFEM-lin}$, which is the best approximation obtained when considering extensions of continuous and piecewise affine functions on $\Gamma$ (corresponding to the linear MsFEM approximation introduced in~\cite{Hou:1997:A-maa} and recalled below, see~\eqref{eq:lin_msfem}), then the error is of order $O(1)$ when $H$ approaches the small scale $\varepsilon$ (this is what the classical theoretical error bound predicts, and this is also what is observed numerically, see e.g.~\cite[Table~II]{Hou:1997:A-maa}). MsFEM type methods are in essence directed towards finding the correct bulk solutions assuming a certain, unknown shape of the solution along the interfaces. The recent history of the development of this category of methods can be revisited as the quest to determine the ``right'' interface conditions.

\medskip

Our approach designs two independent approximation spaces:
\begin{itemize}
\item on the one hand, a space to approach $u_B$ by solving problems similar to~\eqref{eq:lap-coeff} though localized on the elements and with high order polynomials as right-hand sides (a similar idea is used in the recent work~\cite{omnes}). This space is denoted $V_{B,H,\{M_K\}} \subset V_B$, where $\{M_K\}$ is a set of positive integers associating a polynomial degree $M_K$ to each element $K \in \mathcal{T}_H$.
\item on the other hand, a space that approximates $u_\Gamma$ using an harmonic lifting (namely the $A$-harmonic lifting defined by~\eqref{eq:def-harmonic}) of high order polynomials. This space is denoted $V_{\Gamma,H,\{N_e\}} \subset V_\Gamma$, where $\{N_e\}$ is a set of positive integers associating a polynomial degree $N_e$ to each edge $e \subset \Gamma$.
\end{itemize}
We now detail these two approximation spaces.

\medskip

We first consider the bubble space $V_B$. For any element $K$, we choose a positive integer $M_K$ and consider the space of polynomial functions on $K$ of degree lower than or equal to $M_K$. The degree $M_K$ may depend on the element $K$ that we consider because we have in mind local refinement strategies, based on the {\em a posteriori} error estimates we introduce below, which lead to non-uniform discretization parameters. Throughout the article, we adopt the following convention: by degree, we mean {\em total degree} if $K$ is a triangle, and {\em partial degree} in each variable if $K$ is a quadrangle. We denote by ${\cal N}_{M_K}$ the dimension of this space of polynomials and introduce a basis of this space, which we denote $\{P_i\}_{i=1,\dots,{\cal N}_{M_K}}$. For any $1 \leq i \leq {\cal N}_{M_K}$, we introduce the function $\phi_{K,i}^B \in H^1_0(K)$, which is supported in $K$, and which is the solution to
\begin{equation} \label{eq:def_phi_B}
\phi_{K,i}^B = 0 \ \ \text{on $\partial K$} \quad \text{and} \quad \forall v \in H^1_0(K), \ \ \int_K (\nabla v)^T A \nabla \phi_{K,i}^B = \int_K P_i \, v.
\end{equation}
If $K$ is a quadrangular element, we readily note that, in practice, $P_i$ can be chosen as the polynomial that has value 1 at the $i^{\rm th}$ Gauss-Lobatto point and $0$ at the other Gauss-Lobatto points within $K$. Note that we do not consider the case $M_K=0$.

Then, we define the finite dimensional space
\begin{equation} \label{eq:def_VBHM}
V_{B,H,\{M_K\}} =\mathrm{Span} \Big\{ \phi_{K,i}^B, \ \ 1 \leq i \leq {\cal N}_{M_K}, \ \ K \in \mathcal{T}_H \Big\} \subset V_B
\end{equation}
and the approximation~$u_{B,H,\{M_K\}} \in V_{B,H,\{M_K\}}$ of $u_B \in V_B$ as the solution to
\begin{equation} \label{eq:approx_bubble}
\forall v_{B,H,\{M_K\}} \in V_{B,H,\{M_K\}}, \quad 
\int_\Omega (\nabla v_{B,H,\{M_K\}})^T A \nabla u_{B,H,\{M_K\}} = \int_\Omega f \, v_{B,H,\{M_K\}},
\end{equation}
which can equivalently be defined as
\begin{equation} \label{eq:notre_msfem_1}
u_{B,H,\{M_K\}} = \underset{w\in V_{B,H,\{M_K\}}}{\mathrm{argmin}} \left(\frac{1}{2} a(w,w) - \langle f,w\rangle_{L^2(\Omega)}\right).
\end{equation}
Besides considering the above finite dimensional space~\eqref{eq:def_VBHM}, it is also possible to choose $V_{B,H,\{M_K\}}=\{0\}$, in which case $u_B$ is approximated by $u_{B,H,\{M_K\}} = 0$. This crude approximation may be sufficient in some situations since, as briefly mentioned above and as will be detailed below in~\eqref{eq:error_u_bubble0}, we have in this case $\|u_B-u_{B,H,\{M_K\}}\|_E \leq C H$, which is a small error. 

\begin{remark}
We have mentioned above that, in the case of quadrangles, we can choose polynomials $P_i$ associated with the Gauss-Lobatto points. Indeed, such a choice makes the quadrature formulas (to compute the local integrals needed to assemble the stiffness matrix and the right-hand side of~\eqref{eq:approx_bubble}) particularly simple, since $P_i$ vanishes at all but one integration point. From a theoretical viewpoint, any choice of basis is of course possible.
\end{remark}

\medskip

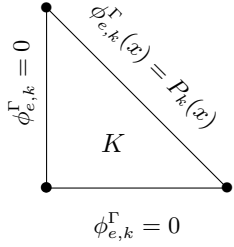
\begin{figure}[htbp]
  \centering
  \begin{tikzpicture}[scale=1.2]
    \draw (0,0) node{$\bullet$} 
    -- (0,2) node{$\bullet$}  node[midway,above=0.2cm,rotate=90,font=\small]{$\phi_{e,k}^\Gamma=0$} 
    -- (2,0) node{$\bullet$} node[midway,above=0.2cm,font=\small ,rotate=-45] {$\phi_{e,k}^\Gamma(x)=P_k(x)$} 
    -- (0,0) node{$\bullet$} node[midway,below=0.2cm,font=\small]{$\phi_{e,k}^\Gamma=0$} ;
    \draw (0.5,0.5) node[anchor=west]{$K$};
  \end{tikzpicture}
  \caption{Local problem defining an edge enrichment $\phi_{e,k}^\Gamma$ for some $k \geq 2$. This enrichment is supported on the two elements sharing the edge $e$. Only one of them is represented. \label{fig:phi_ek}}
\end{figure}

We now turn to the interface space $V_\Gamma$. For any interior edge $e$ of the coarse mesh, we choose a positive integer $N_e$. For any $2 \leq k \leq N_e$, we define the edge enrichment function $\phi_{e,k}^\Gamma$, which is supported on the two elements sharing the edge $e$ (see Figure~\ref{fig:phi_ek}), and which satisfies
\begin{equation}
\label{eq:def_phi_gamma}
\begin{cases}
-\div(A \nabla \phi_{e,k}^\Gamma) = 0 \ \ \text{in $K$}, \\
\phi_{e,k}^\Gamma = P_k \ \ \text{on $e$},\\
\phi_{e,k}^\Gamma = 0  \ \ \text{on $\partial K\setminus e$},
\end{cases}
\end{equation}
where $K$ is any of the two elements containing the edge $e$, and where $P_k$ is a polynomial function of degree $k$ that vanishes at the vertices of the edge $e$. In practice, we work with a so-called {\em boundary-adapted} basis of the polynomial functions of degree lower than or equal to $N_e$. By definition (see e.g. the discussion at the bottom of p.~82 and the central column of Fig.~2.12 on p.~83 of~\cite{canuto2010spectral}), such a basis is composed of two {\em vertex functions} (which are affine on the edge, vanish at one end-point and have value 1 at the other end-point) plus {\em internal functions}, which vanish at both end-points of the edge. The polynomial $P_k$ in~\eqref{eq:def_phi_gamma} is chosen to be an internal basis function, and in practice the internal basis function $\eta_k$ ($k \geq 2$) shown on~\cite[central column of Fig.~2.12 p.~83]{canuto2010spectral} and that we denote below as the $k^{\rm th}$ internal Legendre polynomial.

Formally, the cases $k=0$ and $k=1$ correspond to the linear MsFEM nodal basis functions associated with the two vertices of $e$. Denoting $i_e$ and $j_e$ these two vertices, we set $\phi_{e,0}^\Gamma = \phi_{i_e}^{\rm MsFEM}$ and $\phi_{e,1}^\Gamma = \phi_{j_e}^{\rm MsFEM}$, where $\phi_i^{\rm MsFEM}$ is the solution on any element $K$ to 
\begin{equation} \label{eq:lin_msfem}
-\div (A \nabla \phi_i^{\rm MsFEM}) = 0 \ \ \text{in $K$}, \qquad \phi_i^{\rm MsFEM} = \phi_i \ \ \text{on $\partial K$},
\end{equation}
where $\phi_i$ is the nodal P1 Finite Element basis function associated with the vertex $i$. Note that the support of $\phi_i^{\rm MsFEM}$ is the set of elements having the vertex $i$ as a vertex.

We next define the finite dimensional space
\begin{align}
  V_{\Gamma,H,\{N_e\}}
  &=
  \mathrm{Span} \Big\{ \phi_j^{\rm MsFEM}, \ 1 \leq j \leq N_{\rm vertex}, \quad \phi_{e,k}^\Gamma, \ 2 \leq k \leq N_e, \ e \subset \Gamma \Big\}
  \nonumber
  \\
  &=
  \mathrm{Span} \Big\{ \phi_{e,k}^\Gamma, \ \ 0 \leq k \leq N_e, \ \ e \subset \Gamma \Big\},
\label{def:approx_space_int}
\end{align}
where $N_{\rm vertex}$ is the number of internal vertices of the mesh. We note that $V_{\Gamma,H,\{N_e\}}$ is a subset of $V_\Gamma$. We then define the approximation $u_{\Gamma,H,\{N_e\}} \in V_{\Gamma,H,\{N_e\}}$ of $u_\Gamma \in V_\Gamma$ as the solution to
\begin{equation} \label{eq:approx_interface}
\forall v_{\Gamma,H,\{N_e\}} \in V_{\Gamma,H,\{N_e\}}, \quad  
\int_\Omega (\nabla v_{\Gamma,H,\{N_e\}})^T A \nabla u_{\Gamma,H,\{N_e\}} = \int_\Omega f \, v_{\Gamma,H,\{N_e\}}.
\end{equation}
Again, $u_{\Gamma,H,\{N_e\}}$ can equivalently be defined as
\begin{equation} \label{eq:notre_msfem_2}
u_{\Gamma,H,\{N_e\}} = \underset{w\in V_{\Gamma,H,\{N_e\}}}{\mathrm{argmin}} \left(\frac{1}{2} a(w,w) - \langle f,w\rangle_{L^2(\Omega)}\right).
\end{equation}
Similarly to the degrees $\{ M_K \}$, the degrees $\{ N_e \}$ may vary from one edge to the next, because of possible local refinement strategies.

\medskip

The next two remarks respectively discuss our choice of using polynomials on the edges and our choice to specifically use internal Legendre polynomials.

\begin{remark}
In order to approximate $V_\Gamma$, we have decided to use liftings of polynomials defined on $\Gamma$. Our choice has been motivated by the versatility of polynomials (simplicity of implementation and efficient approximation properties). Other choices could however have been made. The main challenge here is to build an approximation space that accurately captures the oscillations of $u_\Gamma$ on $\Gamma$. For instance, one can think of approaching such oscillating functions by sine functions with increasing frequencies like in Fourier approximation. It turns out that, if we enrich the MsFEM linear basis with liftings of $P_N(x) = \sin(\pi N x /H)$, then we get similar numerical results as with our polynomials. We have chosen to work with polynomials because proving approximation properties for a basis made of sine functions is more delicate than for a polynomial basis, for which we can rely on the extensive theory of polynomial approximation (see e.g. Lemma~\ref{lem:inter_poly} below).
\end{remark}

\begin{remark}
Note that the boundary condition imposed in~\eqref{eq:def_phi_gamma} is continuous on $\partial K$, since we have considered polynomial functions $P_k$ that vanish at the two ends of the edge $e$. If the boundary condition had some jumps on $\partial K$, then the problem~\eqref{eq:def_phi_gamma} would be ill-posed in $H^1(K)$. We have chosen to work with the internal Legendre polynomials (see discussion below~\eqref{eq:def_phi_gamma}) but other choices can be made, as long as the boundary conditions vanish at the two ends of the edge $e$. Our specific choice is motivated by the fact that the internal Legendre polynomials are easy to compute (there is an explicit and simple recursion relation to compute their coefficients, see e.g.~\cite[Eq.~(2.3.31)]{canuto2010spectral}).

Since we see our approach as an enrichment of the MsFEM linear method, the affine nodal functions $\phi_i^{\rm MsFEM}$ must be part of the space spanned by our boundary conditions on $\partial K$. What matters for the analysis is that the space spanned by the boundary conditions on each edge $e$ is the space of polynomial functions of degree lower than or equal to some $N_e$. 
\end{remark}

Denoting by
$$
V_{H,\{M_K\},\{N_e\}} = V_{B,H,\{M_K\}} \oplus V_{\Gamma,H,\{N_e\}},
$$
our approximation $u_{H,\{M_K\},\{N_e\}}$ of $u$ is defined by
\begin{equation} \label{eq:notre_msfem_3}
u_{H,\{M_K\},\{N_e\}}=u_{B,H,\{M_K\}}+u_{\Gamma,H,\{N_e\}}.
\end{equation}
Note that the choices $V_{B,H,\{M_K\}}=\{0\}$ and $N_e=1$ for each edge $e$ leads to an approximation space (and therefore a discrete solution) which is identical to the space used in the classical linear MsFEM approach.

The sets of positive integers $\{M_K, \ K \in \mathcal{T}_H\}$ and $\{N_e, \ e\in \Gamma\}$ define the approximation spaces that are used in the variational problems. As pointed out above, we consider the general case when $M_K$ and $N_e$ may differ from one another to allow for local refinement with respect to the parameters of the method (see Proposition~\ref{prop:error_indicator} and Section~\ref{sec:num_a_posteriori}). For the sake of clarity, in the case when we choose $M_K=M$ for any element $K \in \mathcal{T}_H$ (resp. $N_e=N$ for any edge $e \subset \Gamma$), we replace the notation $\{M_K\}$ by $M$ (resp. $\{N_e\}$ by $N$). 

\medskip

We conclude this section by two general remarks.

\begin{remark}
In practice, one does not have access to the space $V_{H,M,N}$ itself. Indeed, the enrichments $\phi_{e,i}^\Gamma$ or $\phi_{K,i}^B$ are solutions to local problems and must be approximated by a finite element approach using a fine mesh of size $h$ adapted to the characteristic length of variation of the diffusion coefficient $A$.
  
Hence, in practice, for the numerical implementation, we use the space $V_{H,M,N,h} = V_{B,H,M,h} \oplus V_{\Gamma,H,N,h}$ spanned by the functions $\phi_{e,i}^{\Gamma,h}$ and $\phi_{K,i}^{B,h}$, which are the approximation (on the mesh of size $h$) of $\phi_{e,i}^\Gamma$ and $\phi_{K,i}^B$. The study of the convergence of the approach with respect to the parameter $h$ is standard and therefore not considered further in this article. 
\end{remark}

\begin{remark}
The construction of our basis during the off-line stage can be performed in parallel. Indeed, the basis functions for either the bubble or the interface approximation spaces are solutions to independent local problems. 
The stiffness matrix and the right-hand side term for $f=1$ (for the two problems~\eqref{eq:approx_bubble} and~\eqref{eq:approx_interface}) can also be precomputed in parallel during this off-line stage.
\end{remark}

\section{\emph{A priori} and \emph{a posteriori} estimates}
\label{sec:main_res}


We henceforth assume that there exists $0 < \alpha_{\rm min} \leq \alpha_{\rm max}$ such that
\begin{equation} \label{eq:A_ellip}
  \forall \xi \in \R^2, \quad \alpha_{\rm min} \, |\xi|^2 \leq A(x) \xi \cdot \xi \leq \alpha_{\rm max} \, |\xi|^2 \quad \text{a.e. in $\Omega$}.
\end{equation}
Our central \emph{a priori} error estimate reads as follows: 

\begin{proposition} \label{lem:Error_tot}
Assume that $A$ satisfies~\eqref{eq:A_ellip}, that the solution $u$ to~\eqref{eq:lap-coeff} belongs to $H_0^1(\Omega) \cap H^s(\Omega)$ for some $s>3/2$ and that the right hand side $f$ belongs to $H^\ell(\Omega)$ for some integer $\ell \geq 0$. We consider our MsFEM approach~\eqref{eq:notre_msfem_1}--\eqref{eq:notre_msfem_2}--\eqref{eq:notre_msfem_3} in the case when $M_K = M$ for all elements $K$ and $N_e = N$ for all edges $e$, for some $M,N \in \N^\star$. We then have
\begin{multline} \label{eq:mercredi}
  \|u-u_{H,M,N}\|_E \leq \frac{C}{\sqrt{\alpha_{\rm min}}} \ \frac{H^{\min(\ell,M+1)+1}}{M^{\ell+1}} \ \|f\|_{H^\ell(\Omega)}
  \\  
  + C \sqrt{\alpha_{\rm max}} \ \frac{H^{\min(s,N+1)-1}}{N^{s-1}} \ \|u\|_{H^s(\Omega)},
\end{multline}
where the constant $C$ is independent of $H$, $M$, $N$, $A$, $u$ and $f$ (but depends on $\ell$ and $s$). 

In the case when no bubble enrichments are used (that is when $V_{B,H,M}=\{0\}$), we have 
\begin{equation} \label{eq:mercredi_M0}
\|u-u_{H,M,N}\|_E \leq \frac{C}{\sqrt{\alpha_{\rm min}}} \ H \ \|f\|_{L^2(\Omega)} + C \sqrt{\alpha_{\rm max}} \ \frac{H^{\min(s,N+1)-1}}{N^{s-1}} \ \|u\|_{H^s(\Omega)},
\end{equation}
where the constant $C$ is again independent of $H$, $N$, $A$, $u$ and $f$ (but depends on $s$).

The two estimates~\eqref{eq:mercredi} and~\eqref{eq:mercredi_M0} hold in both cases when we use triangular elements throughout the domain, or quadrangular elements.
\end{proposition}

Some remarks are in order.

\medskip

We note that, for $f$ only in $L^2(\Omega)$ (that is $\ell=0$), increasing the polynomial degree $M$ decreases the error at a rate $O(1/M)$. When $f$ is a more regular function, the error decreases with respect to $M$ at a better rate.

\medskip

Extending the above result to the case when the degrees $M_K$ may be different from one element $K$ to the next is not difficult. In contrast, it is unclear to us how to extend it to the case when the degrees $N_e$ used on the edges differ from one edge to the next.

\medskip

We temporarily consider and discuss the classical case of a rescaled periodic matrix coefficient (that is $A(x) = A_{\rm per}(x/\varepsilon)$ for some $\Z^d$-periodic, symmetric coercive matrix $A_{\rm per}$) and a coarse mesh of size $H$ comparable to $\varepsilon$. In that regime, it is observed numerically (see e.g.~\cite[Table~II]{Hou:1997:A-maa}), and this is consistent with the theoretical analysis available, that the classical linear MsFEM approach suffers from an error that does not decrease when $H$ and $\varepsilon$ simultaneously tend to 0 while remaining of the same order of magnitude. In contrast, it is possible in our approach to increase $N$ in order to still have a converging approximation. It is indeed expected that $\| u \|_{H^s(\Omega)} \approx O(\varepsilon^{1-s})$. Choosing $N$ of the order of $1/\varepsilon$ thus guarantees a small error. 

\medskip

Note finally that the efficiency of our approach sensitively depends on the regularity of $u$ \emph{and} on the norm of its derivatives. On the bright side, this implies that the more regular $u$ is, the more efficient our approach is. This unfortunately also means, on the other hand, that the more oscillatory the solution is, the larger the norm of the derivatives of the solution is and thus the larger $N$ has to be taken to obtain a given accuracy.
In this respect, the LOD method~\cite{maalqvist2014localization} is way more robust, since the accuracy only depends on $H$, $f$ and the contrast of $A$ but neither on the regularity nor on the scale of the oscillations. These robustness and accuracy are however obtained at the price of computing ``not so'' local solutions elsewhere than in the given element. 

\begin{remark}
In the periodic case $A(x)= A_{\rm per}(x/\varepsilon)$ mentioned above, we typically have that $\|u\|_{H^s(\Omega)}$ is of the order of $\varepsilon^{1-s}$. In such a case, for given $H$, $M$ and $N$, the right-hand side in the error estimate~\eqref{eq:mercredi} blows up when $\eps \to 0$. This is however not the case of the actual error on the left-hand side. Recall indeed that our approximation space $V_{H,M,N}$ contains the linear MsFEM approximation space, for which the estimate $\dis \| u - u_{\rm MsFEM-lin} \|_E \leq C \left( H + \sqrt{\varepsilon} + \sqrt{\varepsilon/H} \right)$ holds. The error in our approach being smaller than the linear MsFEM error, our approximation does not blow up when $\varepsilon$ goes to $0$ and $H$, $M$ and $N$ are fixed.

This observation questions the sharpness of our error estimate~\eqref{eq:mercredi} in the periodic setting. In the present state of our understanding, we have been unable to derive a sharper estimate, even in this restricted setting.
\end{remark}

\medskip

The proof of Proposition~\ref{lem:Error_tot} is a direct consequence of~\eqref{eq:decoupling} and of the following Lemma~\ref{lem:Error_bubble} and Lemma~\ref{lem:Error_Interface}, which respectively address the bubble approximation and the interface approximation. The proofs of these two lemmas are postponed until Appendices~\ref{sec:proofs_un} and~\ref{sec:proofs_deux}.

\begin{lemma} \label{lem:Error_bubble}
Assume that $A$ satisfies~\eqref{eq:A_ellip} and that $f \in H^\ell(\Omega)$ for some integer $\ell \geq 0$. In the case when $M \geq 1$, the components $u_B$ and $u_{B,H,M}$ satisfy 
\begin{equation}
  \| u_B - u_{B,H,M} \|_E \leq \frac{C_\ell}{\sqrt{\alpha_{\rm min}}} \ \frac{H^{\min(\ell,M+1)+1}}{M^{\ell+1}} \ \|f\|_{H^\ell(\Omega)},
  \label{eq:error_u_bubble}
\end{equation}
for some $C_\ell$ independent of $H$, $M$, $A$ and $f$. If $V_{B,H,M}=\{0\}$, then
\begin{equation}
  \| u_B - u_{B,H,M} \|_E \leq \frac{C}{\sqrt{\alpha_{\rm min}}} \ H \ \|f\|_{L^2(\Omega)},
  \label{eq:error_u_bubble0}
\end{equation}
for some universal constant $C$ (with of course $u_{B,H,M} = 0$).
\end{lemma}

\begin{lemma} \label{lem:Error_Interface}
Assume that $A$ satisfies~\eqref{eq:A_ellip} and that the solution $u$ to~\eqref{eq:lap-coeff} belongs to $H_0^1(\Omega) \cap H^s(\Omega)$ for some $s > 3/2$. Then, the components $u_\Gamma$ and $u_{\Gamma,H,N}$ satisfy 
\begin{equation}
  \| u_\Gamma - u_{\Gamma,H,N} \|_E \leq C_s \sqrt{\alpha_{\rm max}} \ \frac{H^{\min(s,N+1)-1}}{N^{s-1}} \ \|u\|_{H^s(\Omega)},
  \label{eq:error_u_gamma_quad}
\end{equation}
where the constant $C_s$ is independent of $H$, $N$, $A$ and $u$.
\end{lemma}

In sharp contrast with the estimates~\eqref{eq:error_u_bubble} and~\eqref{eq:error_u_bubble0} which do not depend on the oscillations of $A$, the estimate~\eqref{eq:error_u_gamma_quad} depends on the norm of derivatives of $u$, hence, indirectly on the oscillations of $A$. As expected, the interface component $u_\Gamma$ is more delicate to capture than the bubble component $u_B$.

For the ACMS method, the estimate shown in~\cite{hetmaniuk2014error} depends on the $k^{\rm th}$ largest eigenvalue $\lambda_e^k$ for the associated edge eigenproblem. The rate of decrease of $\lambda_e^k$ with respect to $k$ and $H$ is not known, although the numerical experiments empirically suggest that it is $O(k/H)$, which would give an error estimate similar to~\eqref{eq:error_u_gamma_quad}.

The proof of Lemma~\ref{lem:Error_bubble} and Lemma~\ref{lem:Error_Interface} 
essentially follows, and it is not unexpected, the pattern of the proof of the classical C\' ea's Lemma. The best approximation is estimated using the Legendre projection (for Lemma~\ref{lem:Error_bubble}) or the Legendre interpolant on the bulk and the lifting of the interpolant along the edges (for Lemma~\ref{lem:Error_Interface} in the case of quadrangles). 
Some technicalities arise for Lemma~\ref{lem:Error_Interface} in the case of triangular meshes and an alternative proof (which actually also covers the case of quadrangles) using $hp$-Finite Element methods must be used. 
We will return to this in Appendix~\ref{sec:proofs_deux}.

\medskip

We now turn to our \emph{a posteriori} error estimator. In contrast to our \emph{a priori} estimates above, we now consider the general case when the polynomial degrees $N_e$ (resp. $M_K$) associated to each edge $e$ (resp. each element $K$) can be different. For some technical reasons (in particular due to the use of Scott-Zhang interpolation results, see Lemma~\ref{lem:clem_interp}), we assume that the polynomial degrees of the edges are comparable on neighboring edges, in the sense that
\begin{equation} \label{eq:degre_edge_proche}
\forall e, e' \in \Gamma \ \text{s.t.} \ \overline{e} \cap \overline{e'} \neq \emptyset, \quad \frac{N_e}{\sqrt{\gamma}} \leq N_{e'} \leq \sqrt{\gamma} \, N_e,
\end{equation}
where $\gamma$ is the mesh regularity constant of~\eqref{eq:def_mesh_regular}.

\begin{proposition} \label{prop:error_indicator}
Assume that the diffusion coefficient matrix $A$ satisfies~\eqref{eq:A_ellip} and belongs to $(C^1(\overline{\Omega}))^{d \times d}$. We also assume that there exists some integer $\overline{\ell}$ such that, for any element $K$ of the coarse mesh, $f \in H^{\ell_K}(K)$ for some integer $\ell_K \geq 0$ which satisfies $\ell_K \leq \overline{\ell}$. 

Consider the MsFEM approach on the discrete space $V_{H,\{M_K\},\{N_e\}}$, where $M_K > 0$ is the maximal degree of the polynomial functions used as right-hand sides for the bubble basis functions in the element $K$, and $N_e > 0$ is the maximal degree of the polynomial functions used as boundary conditions for the interface basis functions associated to the edge $e$. We assume that the degrees $\{ N_e \}$ satisfy~\eqref{eq:degre_edge_proche}.

For any $\eta > 0$, the discrete solution $u_{H,\{M_K\},\{N_e\}}$ satisfies the \emph{a posteriori} estimate
\begin{align}
  & \| u-u_{H,\{M_K\},\{N_e\}} \|_E
  \nonumber
  \\
  & \leq 
  C_{\eta,A} 
  \left\{ \sum_{K\in\mathcal{T}_H} H_K^2 \, \frac{H_K^{\min(\ell_K,M_K+1)}}{M_K^{\ell_K}} \, \| f + \div \left( A \nabla u_{B,H,\{M_K\}} \right) \|_{L^2(K)} \, \| f \|_{H^{\ell_K}(K)} \right.
  \nonumber
  \\ 
  &\left.
  + \sum_{K \in \mathcal{T}_H} \|f\|_{L^2(K)}^2 \left( \sum_{e\subset \partial K} \frac{H_e \, H_K}{N_e^{1-2\eta} \, p_e} \right) + \sum_{e\subset \Gamma} \frac{H_e}{p_e} \left\| J_e \big( \nu^T A \nabla u_{\Gamma,H,\{N_e\}} \big) \right\|_{L^2(e)}^2 \right\}^{1/2}
  \label{eq:a_posteriori_estimator}
\end{align}
where $J_e(\psi)$ denotes the jump of a given function $\psi$ across the edge $e$, and $\nu$ is a normal vector to the edge. In the above estimate, $H_K$ is the diameter of the element $K$, $H_e$ is the length of the edge $e$ and we have set $p_e = \min \{ N_{\widetilde{e}} \ | \ \widetilde{e} \subset \partial K^1_e \cup \partial K^2_e \}$ where $K^1_e$ and $K^2_e$ are the two elements sharing the edge $e$. The constant $C_{\eta,A}$ depends only on $\eta$, on the regularity parameter $\gamma$ of the mesh (see~\eqref{eq:def_mesh_regular}), on $\overline{\ell}$ and on the diffusion coefficient $A$ through $\alpha_{\rm min}$ and $\| A \|_{C^1(\overline{\Omega})}$. 

Without bubble enrichment, that is when $V_{B,H,\{M_K\}} = \{ 0 \}$, we have the estimate
\begin{align}
  &\| u-u_{H,\{M_K\},\{N_e\}} \|_E
  \leq 
  C_{\eta,A}
  \left\{ \sum_{K\in\mathcal{T}_H} H_K^2 \, \| f \|_{L^2(K)}^2 \right.
  \label{eq:a_posteriori_estimator_no_bubble}
  \\ 
  &\left.
  \qquad + \sum_{K \in \mathcal{T}_H} \|f\|_{L^2(K)}^2 \left( \sum_{e \subset \partial K} \frac{H_e \, H_K}{N_e^{1-2\eta} \, p_e} \right) + \sum_{e\subset \Gamma} \frac{H_e}{p_e} \left\| J_e \big( \nu^T A \nabla u_{\Gamma,H,\{N_e\}} \big) \right\|_{L^2(e)}^2 \right\}^{1/2}.
\nonumber
\end{align}
\end{proposition}

In contrast to Proposition~\ref{lem:Error_tot}, where our assumptions on the regularity of $f$ and $u$ somehow only \emph{implicitly} presuppose some regularity of $A$, we \emph{explicitly} assume in Proposition~\ref{prop:error_indicator} some given regularity of $A$.

\medskip

The right-hand side of~\eqref{eq:a_posteriori_estimator} actually defines an error indicator: the actual error is bounded from above by the product of a computable indicator within the brackets (involving the two components of the numerical solution) times a constant independent of $H_K$, $H_e$, $M_K$ and $N_e$.

The proof of Proposition~\ref{prop:error_indicator}, which is postponed until Appendix~\ref{sec:proofs_trois}, follows the analogous proof performed for the ACMS method in~\cite{hetmaniuk2014error}. However, Scott-Zhang type polynomial interpolation has to be introduced instead of classical polynomial interpolation. 

Some illustrations regarding the behavior (and in particular the effectivity) of the \emph{a posteriori} estimator are presented in Section~\ref{sec:Num_exp}. 

\section{Numerical experiments}
\label{sec:Num_exp}

This section is divided into two parts. We first compare our approach to standard MsFEM approximations (linear MsFEM and oversampling MsFEM approaches), and to the ACMS method of~\cite{hetmaniuk2010special,hetmaniuk2014error}. Second, we investigate the performance of the \emph{a posteriori} estimator proposed in Proposition~\ref{prop:error_indicator}. All our numerical experiments have been performed with FreeFem++~\cite{freefem++}.

\subsection{Comparison with other MsFEM approaches and with the ACMS approach}

\begin{figure}[htbp]
  \hskip -1 truecm
  \includegraphics[width=0.6\textwidth]{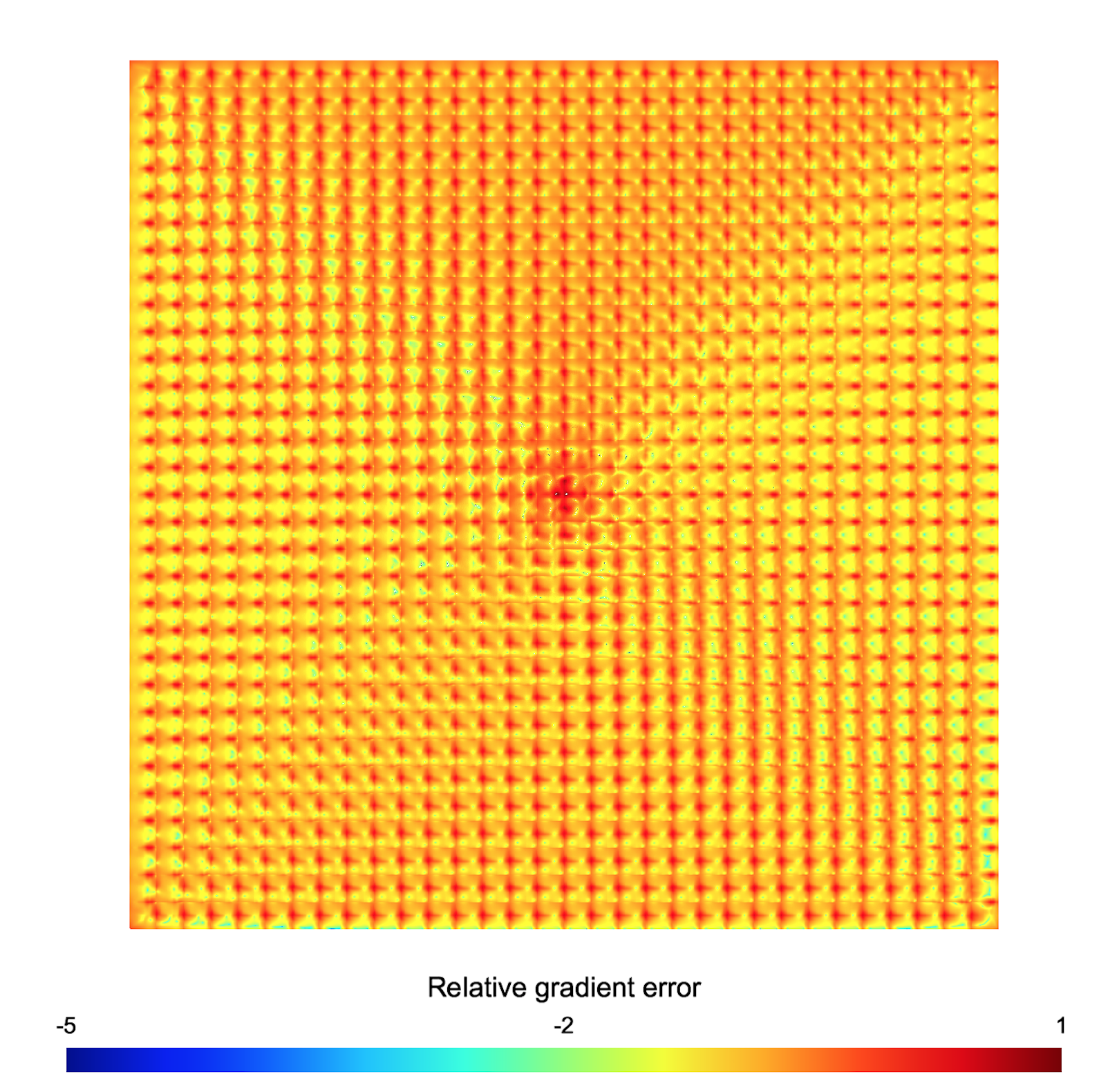}
  \includegraphics[width=0.55\textwidth]{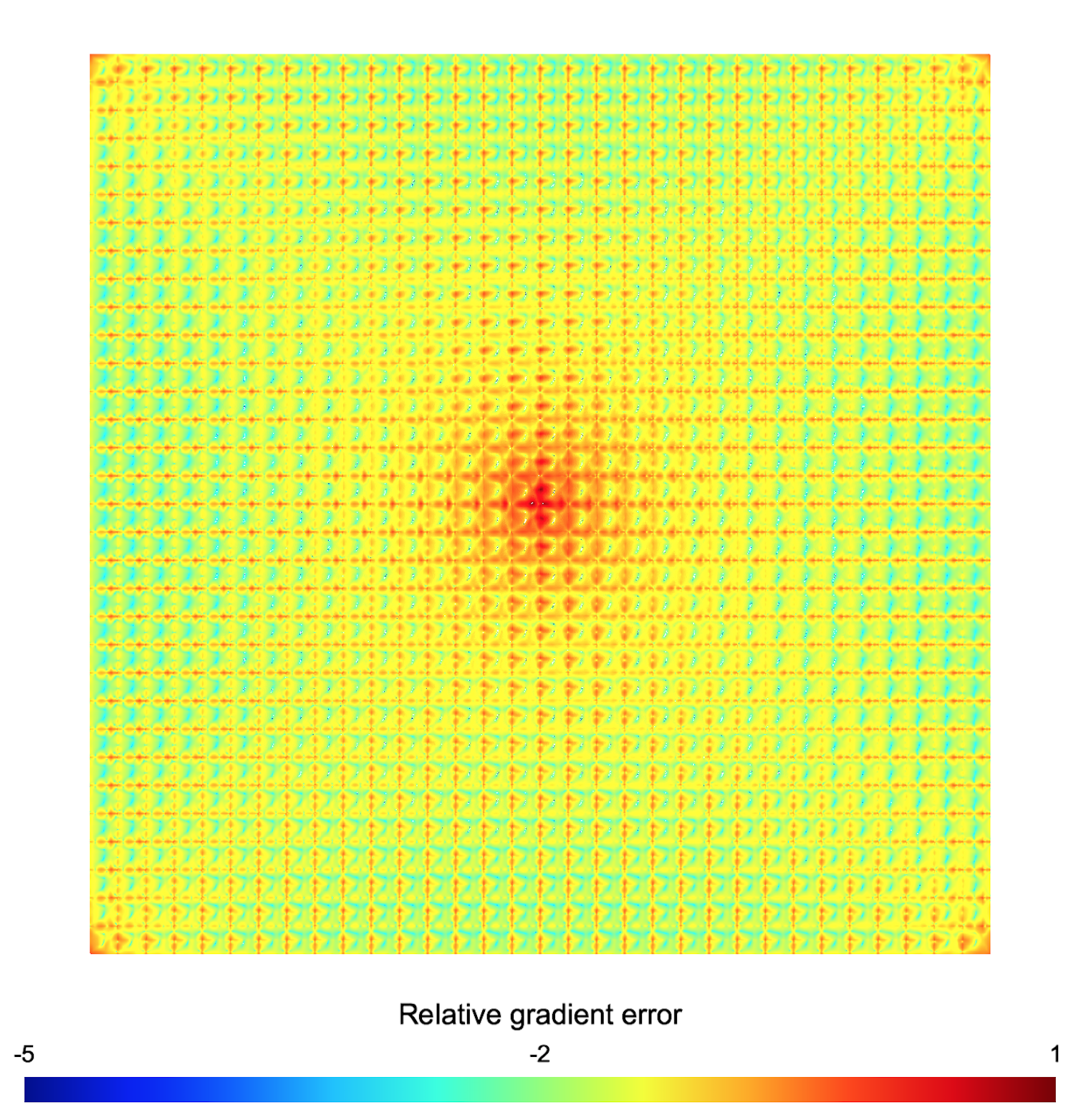}
  \caption{Error on the approximation of the gradient for MsFEM-lin (left) and for our approach with 4 polynomials (right): the accuracy is poor in the red regions, much better in the light green regions, and excellent in the dark blue regions. \label{error_grad_MsFEM}}
\end{figure}

In our numerical experiments, the emphasis is put on the enrichment by edge functions. As already mentioned above, the bubble error (that is, the first term in the right-hand side of~\eqref{eq:decoupling}) when no bubble enrichments are used behaves like classical FE estimates for the Poisson problem: it decreases linearly with respect to $H$, with a prefactor that only depends on the $L^2$ norm of the right-hand side and the coercivity constant of the diffusion coefficient $A$ (see~\eqref{eq:error_u_bubble0}).
In contrast, the interface error (that is, the second term in the right-hand side of~\eqref{eq:decoupling}) depends on the oscillations of $A$ and has a more intricate behavior. Moreover, in the classical MsFEM approaches (linear and oversampling), the basis functions belong to $V_\Gamma$. Such approaches can hence also be enriched by bubble elements. In order to compare their respective effectiveness, it thus appears that it is best not to consider bubble enrichments. For our tests, we therefore only act on $V_\Gamma$ and use no bubble enrichment, that is we keep $V_{B,H,\{M_K\}} = \{ 0 \}$. Such a choice will be illustrated below (see Figure~\ref{error_grad_MsFEM} and the associated discussion).

In what follows, we work with the uniform choice $\{N_e\}=N$ for all edges $e \subset \Gamma$. The linear MsFEM approach corresponds to the choice $N=1$ (and $V_{B,H,\{M_K\}} = \{ 0 \}$).

We solve~\eqref{eq:lap-coeff} for a classical benchmark test introduced in~\cite{Hou:1997:A-maa}, where $A$ is periodic and oscillates at the scale $\varepsilon$. More specifically, we consider 
\begin{equation} \label{eq:def_Aeps}
A_\varepsilon(x) = a\left(\frac{x}{\varepsilon},\frac{y}{\varepsilon}\right)I_2, \qquad a(x,y) = \frac{2+1.8\sin(2\pi x)}{2+1.8\cos(2\pi y)} + \frac{2+\sin(2\pi y)}{2+1.8\sin(2\pi x)},
\end{equation}
where $I_2$ is the identity matrix, and solve
\begin{equation}
\label{pb:num_pb}
-\div (A_\varepsilon \nabla u_\varepsilon)= -1 \quad \text{in $\Omega$}, \qquad u_\varepsilon=0 \quad \text{on $\partial \Omega$},
\end{equation}
on the domain $\Omega = (0,1)^2$. We consider $\varepsilon$ ranging from $1/32$ to $1/128$.

In order to compute errors, we have computed a reference solution to~\eqref{pb:num_pb} using P2 Finite Elements with a mesh of size $h=1/2048$. Note that $h \ll \varepsilon$ for the range of values of $\varepsilon$ that we consider, so that this classical finite element approach can be considered accurate. Similarly, on each element $K$, the interface basis functions $\phi_{e,i}^{\Gamma}$ have no analytical expression and are approximated using P1 Finite Elements on a mesh of a small size (of the order of $h$). 

\medskip

To start with, we wish to illustrate our above somewhat intuitive statement regarding the fact that multiscale approaches typically do a better job at approximating the solution in the bulk than on the interfaces, thus the interest of focusing our study and our efforts on the enrichment (by Legendre polynomials) on the edges. To support this claim, we show on Figure~\ref{error_grad_MsFEM} the typical error obtained using the linear version of MsFEM (left). We specifically show the relative error $\dps \log_{10} \left[ \left| \nabla \left(u_\varepsilon-u_\varepsilon^{H,M,N} \right) \right| / \left| \nabla u_\varepsilon \right| \right]$ as a function of $x \in \Omega$. The largest errors are evidently concentrated on the interfaces (that is on the edges of our quadrangular mesh, which is clearly visible on the figure). Already an enrichment of three Legendre polynomials per edge (i.e. using polynomials of degree up to $N=4$ on the edges) allows one to dramatically reduce the latter error, as shown on the right of Figure~\ref{error_grad_MsFEM}. 
Further enriching the description of the solution along the edges with a larger number of Legendre polynomials, say $N=10$, would typically render the error almost homogeneous throughout the computational domain. All in all, the above set of comments justify our tactical choice to keep $V_{B,H,M} = \{ 0 \}$ and focus on increasing $N$.

\medskip

\begin{figure}[htbp]
  \hskip -1truecm
  \includegraphics[width=0.6\textwidth]{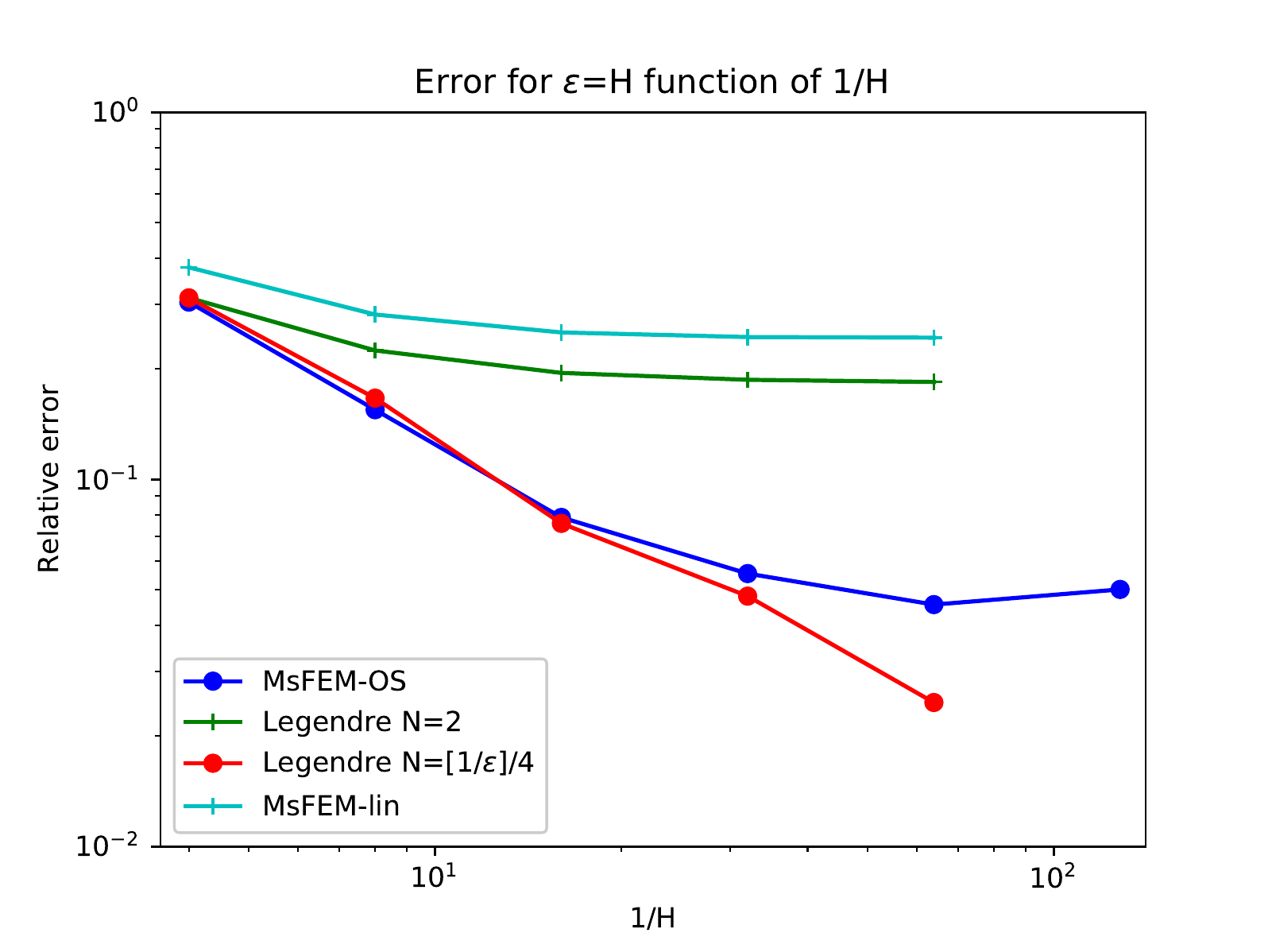}
  \includegraphics[width=0.6\textwidth]{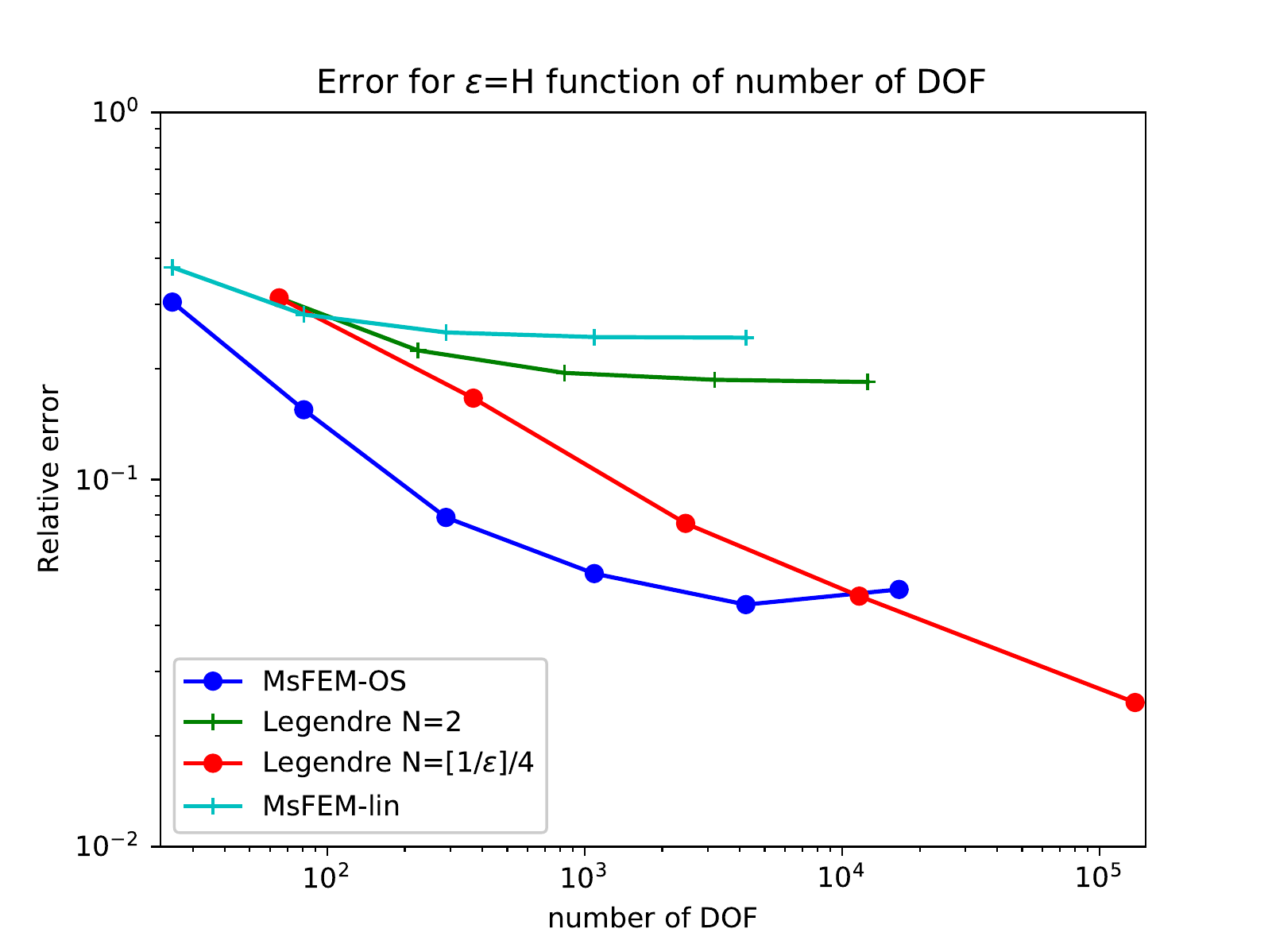}
  \caption{Compared performances in the regime $H\approx \varepsilon$. \label{res_sin}}
\end{figure}

Our next observation is the purpose of Figure~\ref{res_sin} and concerns the relatively poor performance of the linear MsFEM and of the oversampling MsFEM approaches in the regime $H\approx\varepsilon$ (often called the \emph{resonance regime}) where the coarse mesh size matches the typical size of oscillations. Whether one argues in terms of the meshsize $H$ (left of Figure~\ref{res_sin}), or in terms of the number of degrees of freedom (right of that figure), these two MsFEM variants saturate, while the Legendre enriched approach performs increasingly better. In Figure~\ref{res_sin}, and likewise in Figures~\ref{res_leg_msfem__hou},~\ref{res_leg_eigen_eps_32} and~\ref{res_leg_tri_H} below, all errors are relative errors in the energy norm
\begin{equation}
\label{def:err_relative}
\mathcal{E}_{\rm rel} = \sqrt{ \frac{a(u_\varepsilon-u_\varepsilon^{H,M,N},u_\varepsilon-u_\varepsilon^{H,M,N})}{a(u_\varepsilon,u_\varepsilon)} },
\end{equation}
while the approaches we test are respectively denominated as \emph{MsFEM-lin} for the standard version of linear MsFEM, \emph{MsFEM-OS} for its variant using oversampling (where the oversampling domain is 3 times larger in each direction than the original coarse element: we thus consider quadrangles of size $(3H)^2$ rather than $H^2$, and likewise for triangles), \emph{Legendre N = \dots} for the approach presented here using the corresponding degree $N$ of Legendre polynomials on the edges (for instance, \emph{Legendre N=2} corresponds to adding one enrichment per edge wrt the MsFEM-lin approach), and \emph{Eigen \dots} for the ACMS approach.  

The relative error~\eqref{def:err_relative} is easy to compute. Introduce indeed the energy 
$$
\mathcal{E}(v) = \frac{1}{2} \int_\Omega (\nabla v)^T A_\varepsilon \nabla v - \int_\Omega fv,
$$
defined for any $v \in H^1_0(\Omega)$. Denoting $\mathcal{E}^\star = \mathcal{E}(u_\varepsilon)$ and using that the matrix $A_\varepsilon$ is symmetric and the variational formulation of~\eqref{eq:lap-coeff}, we compute that, for any $v \in H^1_0(\Omega)$,
$$
\mathcal{E}(u_\eps + v) - \mathcal{E}^\star = \frac{1}{2} \int_\Omega (\nabla v)^T A_\varepsilon \nabla v = \frac{1}{2} a(v,v),
$$
and thus $\dps \mathcal{E}(v) - \mathcal{E}^\star = \frac{1}{2} a(u_\varepsilon-v,u_\varepsilon-v)$ for any $v \in H^1_0(\Omega)$. We therefore deduce from~\eqref{def:err_relative} that $\mathcal{E}_{\rm rel} = \sqrt{ \left( \mathcal{E}\left(u_\varepsilon^{H,M,N}\right) -\mathcal{E}^\star \right)/(-\mathcal{E}^\star) }$.

\begin{remark} \label{rem:energy_comput}
This definition of the error is in practice very useful because computing~\eqref{def:err_relative} only requires to compare the two scalar quantities $\mathcal{E}^\star$ and $\mathcal{E}(u_\varepsilon^{H,M,N})$, which can be obtained independently. We get $\mathcal{E}^\star$ by computing the energy for our reference solution, while $\mathcal{E}(u_\varepsilon^{H,M,N})$ can be computed in parallel over the coarse elements $K$ once the global problem has been solved at the on-line stage. In particular, it is neither necessary to store the reference and numerical solutions, nor to compute their difference on a fine common Finite Element space, an operation which would be computationally expensive.
\end{remark}

\medskip

We next perform, for $\varepsilon$ fixed (namely at the value $\varepsilon=1/32$), and $H$ decreasing from $1/4$ to $1/64$ (or, correspondingly, the number of degrees of freedom increasing), full comparisons of the accuracy obtained for the various methods considered, including the ACMS method. The results are shown on Figures~\ref{res_leg_msfem__hou} and~\ref{res_leg_eigen_eps_32} respectively.

For any fixed $H$, our approach is more accurate than the MsFEM oversampling method when $N$ is large enough, say here $N \geq 9$ (see left side of Figure~\ref{res_leg_msfem__hou}). For $N=5$, our approach and the MsFEM oversampling method essentially share the same accuracy. The oversampling variant is more accurate for smaller values of $N$. However, for a fixed $H$, our approach uses more degrees of freedom than the oversampling approach. We thus compare the approaches for a given number of degrees of freedom on the right side of Figure~\ref{res_leg_msfem__hou}. When $H$ is not too small (and thus the number of degrees of freedom is not too large), the MsFEM oversampling method provides better results than our approach. However, for smaller values of $H$ (and thus larger numbers of degrees of freedom, say larger than $10^4$), our approach outperforms the MsFEM oversampling method. We also notice that the oversampling approach suffers from the resonance effect mentioned above (the error is essentially the same for any $H$ between $1/128$ and $1/32$), whereas our approach provides an error which is monotonically decreasing with $H$, a fact which is advantageous from a practical viewpoint.

Our tests of Figure~\ref{res_leg_eigen_eps_32} clearly show that our approach is equally accurate as (and in some cases more accurate than) the ACMS method, for each given level of enrichment. We recall that it is however less expensive, since solving boundary value problems is less expensive than solving eigenvalue problems.

\begin{figure}[htbp]
  \hskip -1truecm
  \includegraphics[width=0.6\textwidth]{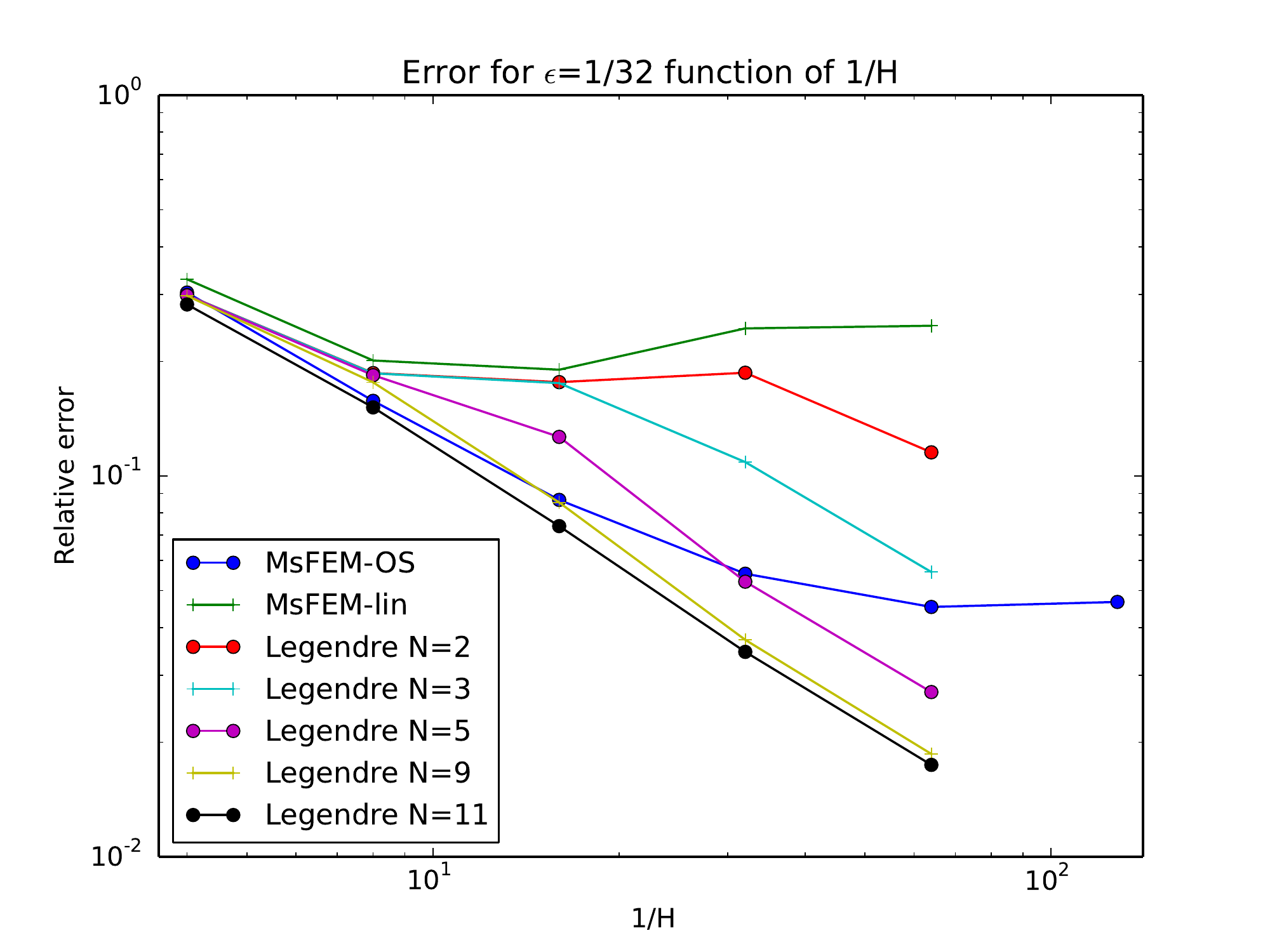}
  \includegraphics[width=0.6\textwidth]{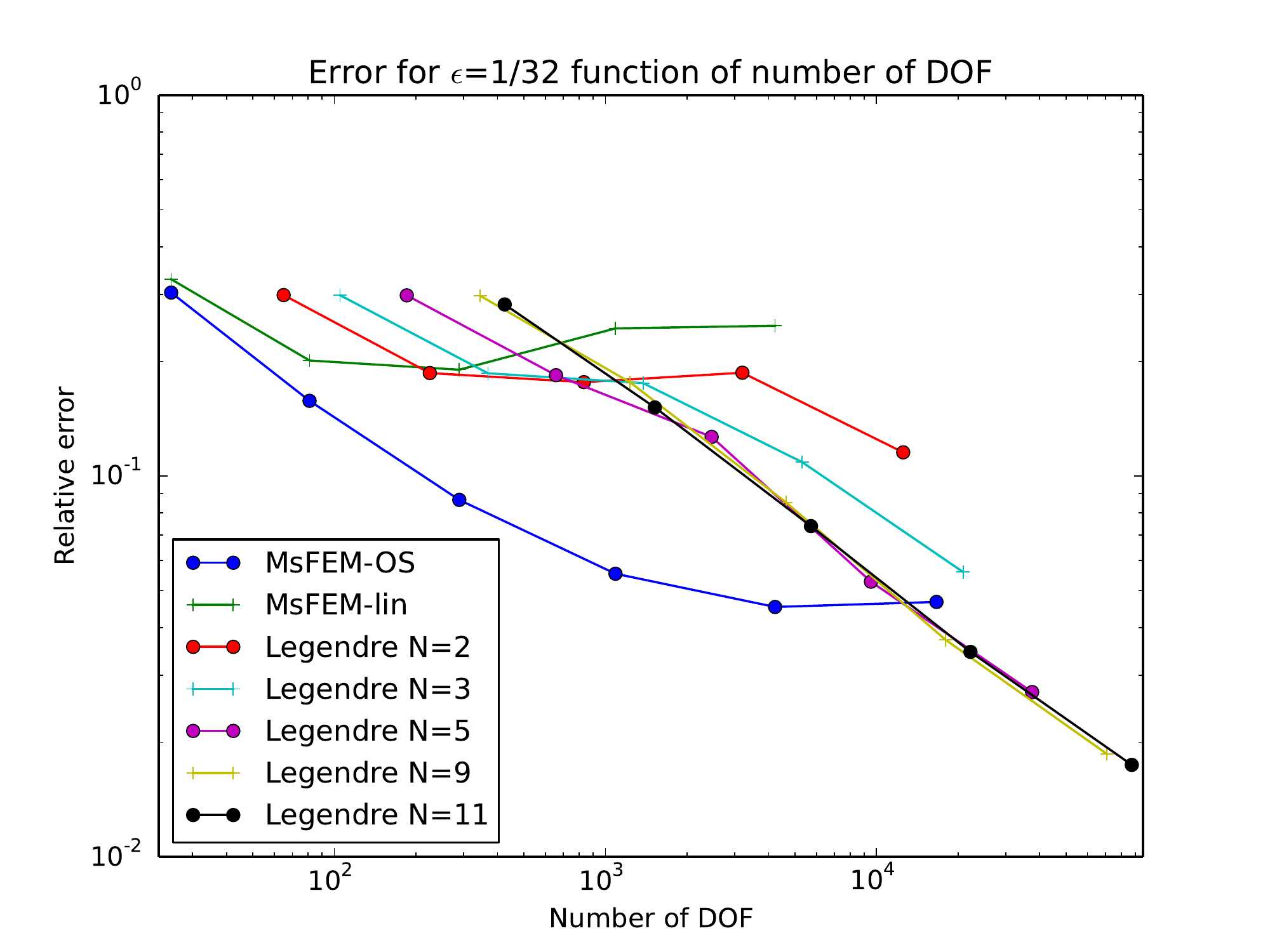}
  \caption{Comparison of our approach with classical MsFEM approaches. \label{res_leg_msfem__hou}}
\end{figure}

\begin{figure}[htbp]
  \begin{center}
    \includegraphics[width=0.6\textwidth]{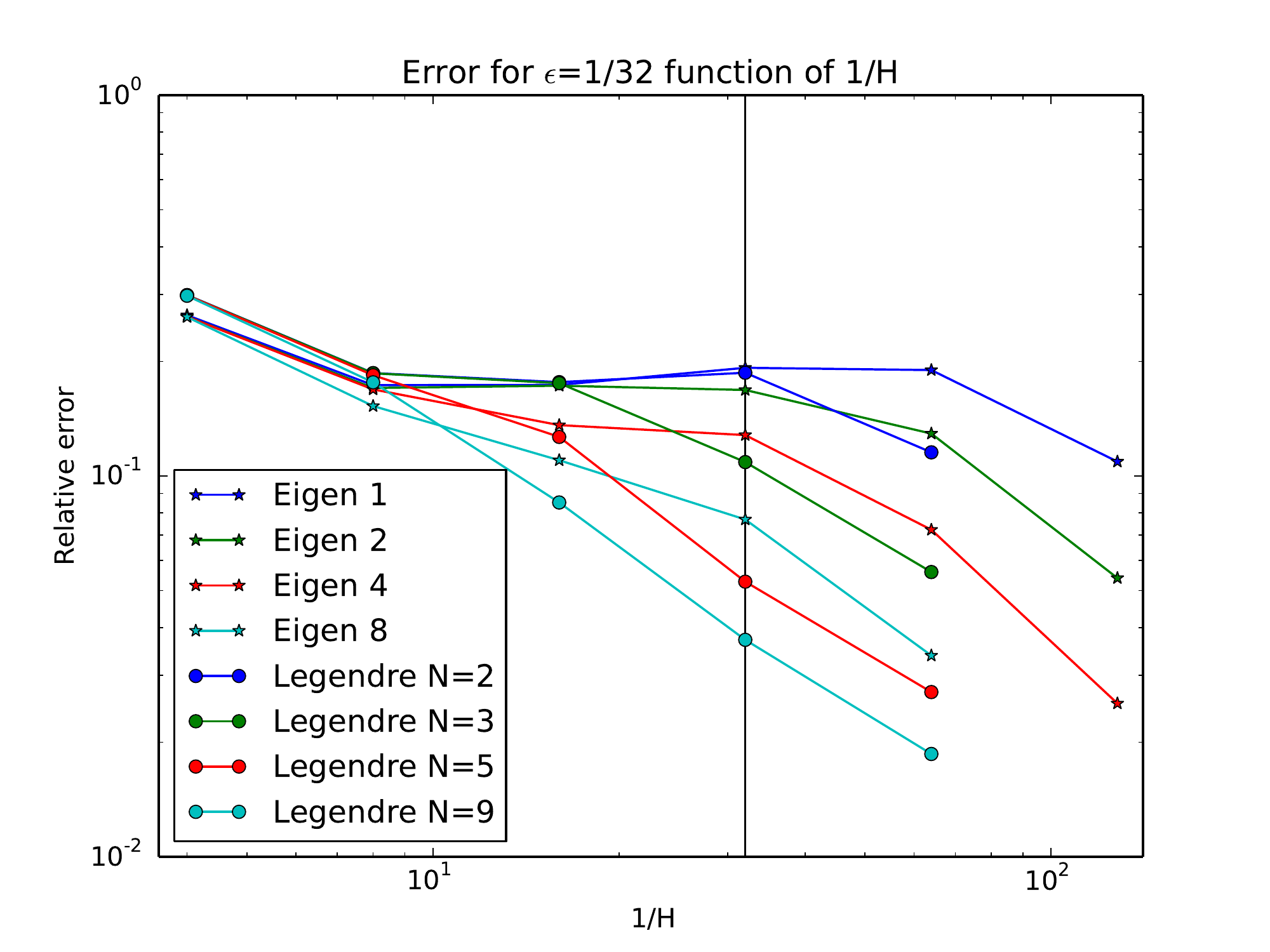}
  \end{center}
  \caption{Comparison of our approach with the ACMS method, at equal number of enrichments per edge (for instance, ``Eigen 1'' and ``Legendre N=2'' both correspond to adding one enrichment per edge wrt the MsFEM-lin approach). \label{res_leg_eigen_eps_32}}
\end{figure}
  
\medskip

Our next test compares the performance of our approach for triangular meshes and for quadrangular meshes. We set $\varepsilon=1/32$ and present the relative energy error as a function of $1/H$ (left of Figure~\ref{res_leg_tri_H}) and of the number of degrees of freedom (right of Figure~\ref{res_leg_tri_H}). Our conclusion is that, essentially, the approach performs equally well in both cases, thereby making possible the application to a large class of computational domains, with intricate geometries for which quadrangular meshes cannot be used.

\begin{figure}[htbp]
  \hskip -1truecm
  \includegraphics[width=0.6\textwidth]{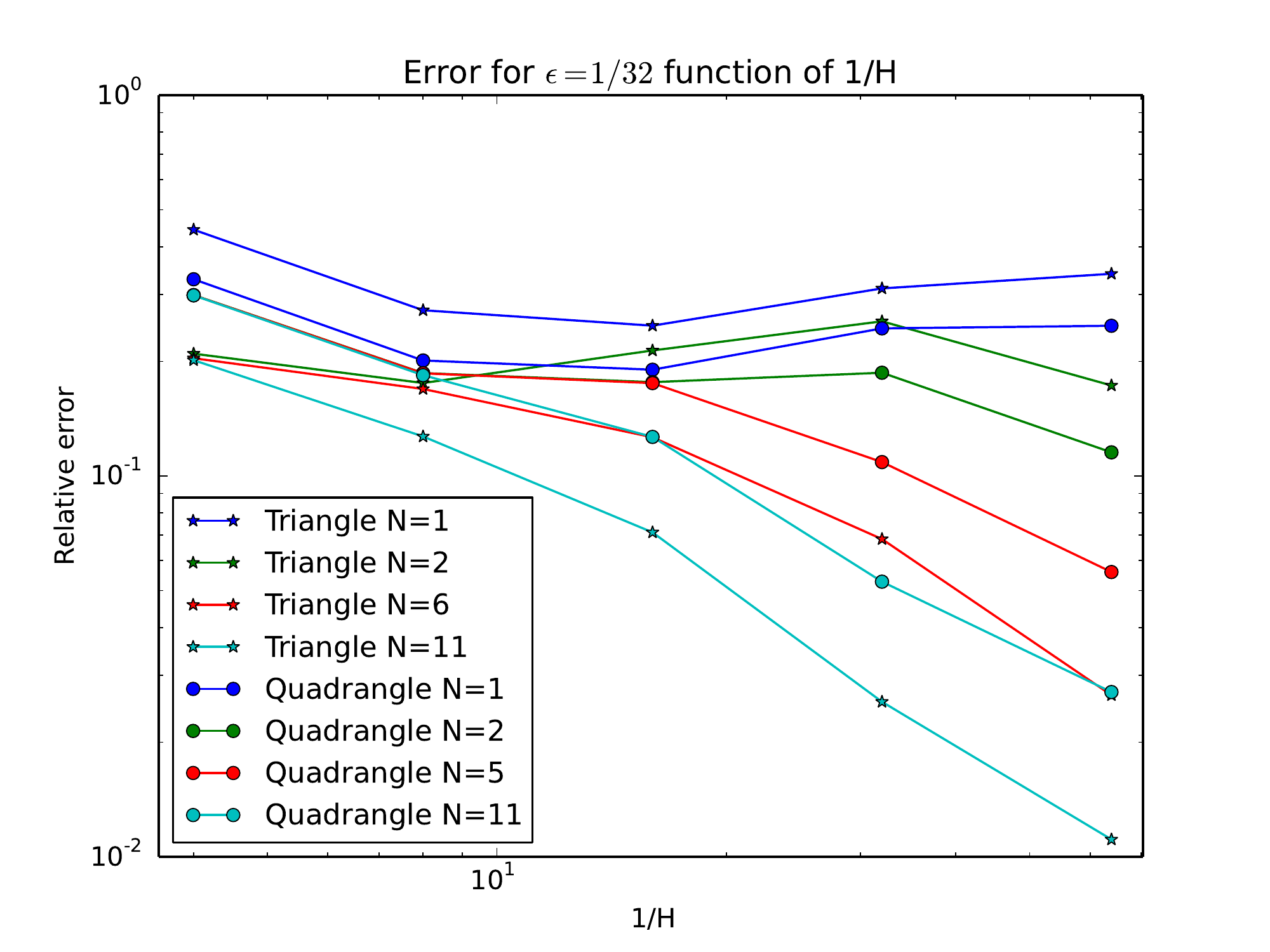}
  \includegraphics[width=0.6\textwidth]{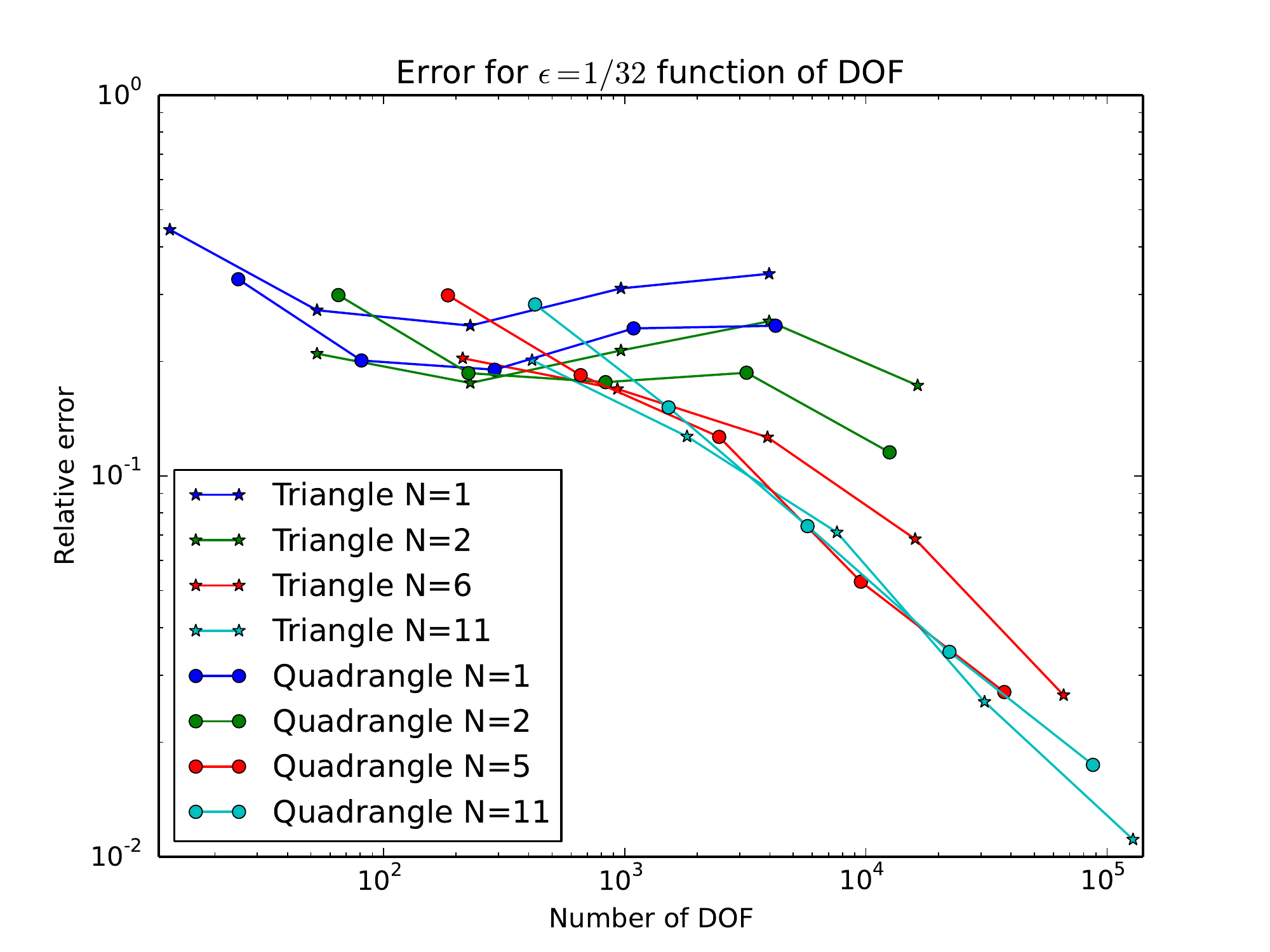}
  \caption{Our approach for triangles and quadrangles, in terms of $1/H$ (left) or of the number of degrees of freedom (right). The approaches ``Triangle N=1'' and ``Quadrangle N=1'' both correspond to the MsFEM-lin approach, on triangular (resp. quadrangular) meshes. The approaches ``Triangle N=2'' and ``Quadrangle N=2'' both correspond to adding one enrichment per edge wrt the MsFEM-lin approach. \label{res_leg_tri_H}}
\end{figure}

\subsection{\emph{A posteriori} estimator} \label{sec:num_a_posteriori}

We now investigate the performance of the \emph{a posteriori} estimate provided in Proposition~\ref{prop:error_indicator}. As in the previous section, we do not use bubble enrichments. This is why, instead of using as before the full energy error~\eqref{def:err_relative}, we now use the interface error defined by 
\begin{equation}
  \label{def:err_rel_gamma}
  \mathcal{E}_{\rm rel,\Gamma} = \sqrt{ \frac{a(u_\varepsilon^\Gamma-u_\varepsilon^{H,M,N},u_\varepsilon^\Gamma-u_\varepsilon^{H,M,N})}{a(u_\varepsilon^\Gamma,u_\varepsilon^\Gamma)} }.
\end{equation}
Note that, since $V_{B,H,M} = \{ 0 \}$, the function $u_\varepsilon^{H,M,N}$ belongs to $V_\Gamma$ and is meant to approximate $u_\varepsilon^\Gamma$.

We compare this actual relative error with the error indicator given in~\eqref{eq:a_posteriori_estimator_no_bubble}, and more precisely with the indicator of the interface error (where we have set $\eta=0$), that is 
\begin{equation} \label{def:err_rel_gamma_post}
  \mathcal{E}_{\rm post,\Gamma} 
  =
  \left\{ \sum_{K \in \mathcal{T}_H} \|f\|_{L^2(K)}^2 \left( \sum_{e \subset \partial K} \frac{H_e \, H_K}{N_e \, p_e} \right) + \sum_{e\subset \Gamma} \frac{H_e}{p_e} \left\| J_e \big( \nu^T A \nabla(u_{\Gamma,H,\{N_e\}}) \big) \right\|_{L^2(e)}^2 \right\}^{1/2}.
\end{equation}

\begin{remark} \label{rem:comput_rel_error}
As above, we compute the relative error~\eqref{def:err_rel_gamma} using the energy $\mathcal{E}$. The orthogonal decomposition~\eqref{eq:orthogonal_decomposition} ensures that $u_\eps^\Gamma$ is also the minimizer of the energy on $V_\Gamma$. Hence, the error $a(u_\varepsilon^\Gamma-u_\varepsilon^{H,M,N},u_\varepsilon^\Gamma-u_\varepsilon^{H,M,N})$ can also be expressed as the difference between the energy of our approximation (which belongs to $V_\Gamma$ since $V_{B,H,M}=\{0\}$), which is computed as explained in Remark~\ref{rem:energy_comput}, and the energy of $u_\eps^\Gamma$. We compute the energy of $u_\eps^\Gamma$ by computing explicitly a reference solution for $u_\eps^B$ (this simply requires to solve in parallel homogeneous Dirichlet problems in all the elements $K$). The energy of $u_\eps^\Gamma$ is equal to $\mathcal{E}(u_\eps) - \mathcal{E}(u_\eps^B)$. This procedure is simpler than computing $u_\eps^\Gamma$, a task for which we would need to store the value of $u_\eps$ on $\Gamma$.
\end{remark}

\begin{figure}[htbp]
  \includegraphics[width=0.5\textwidth]{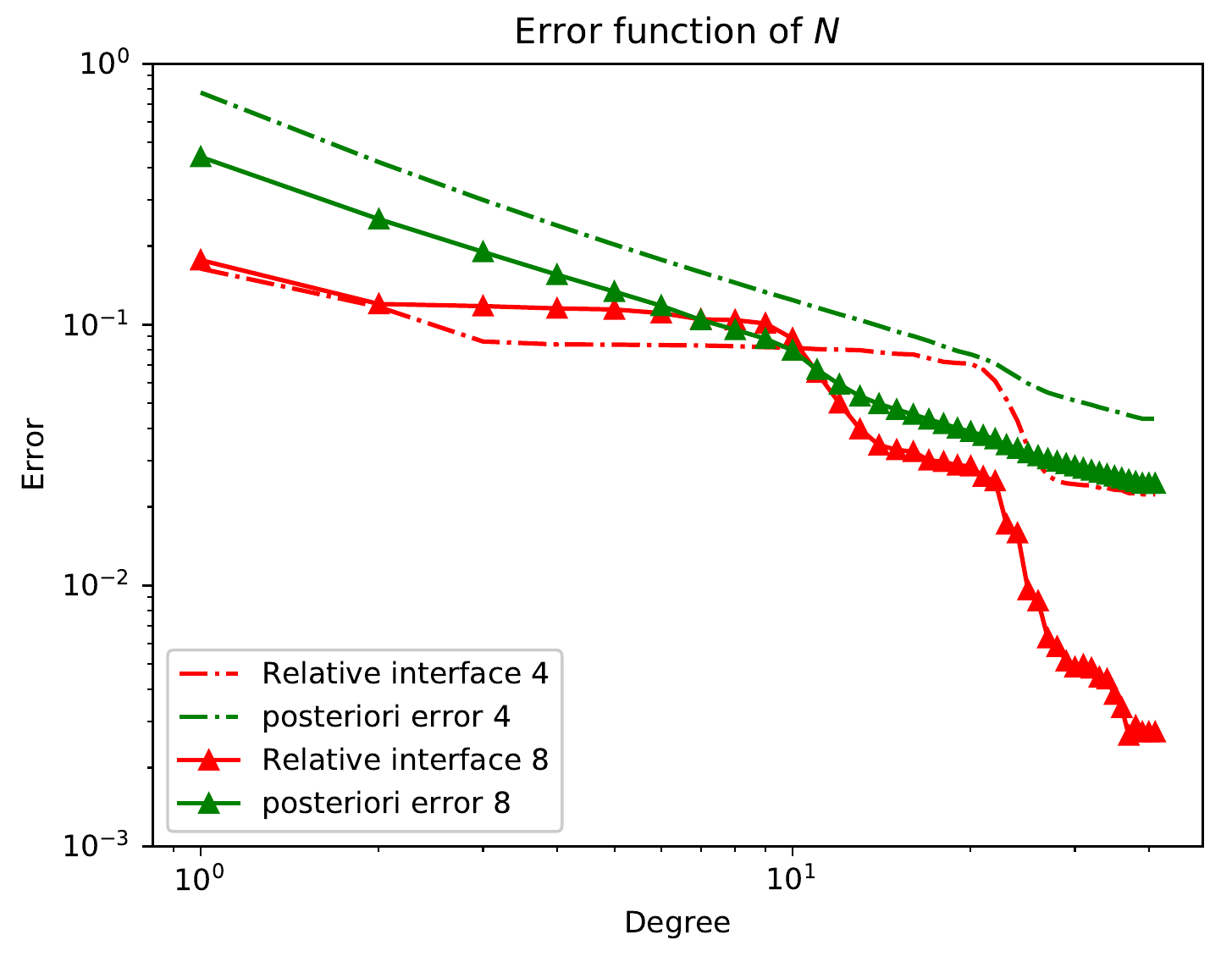}
  \includegraphics[width=0.5\textwidth]{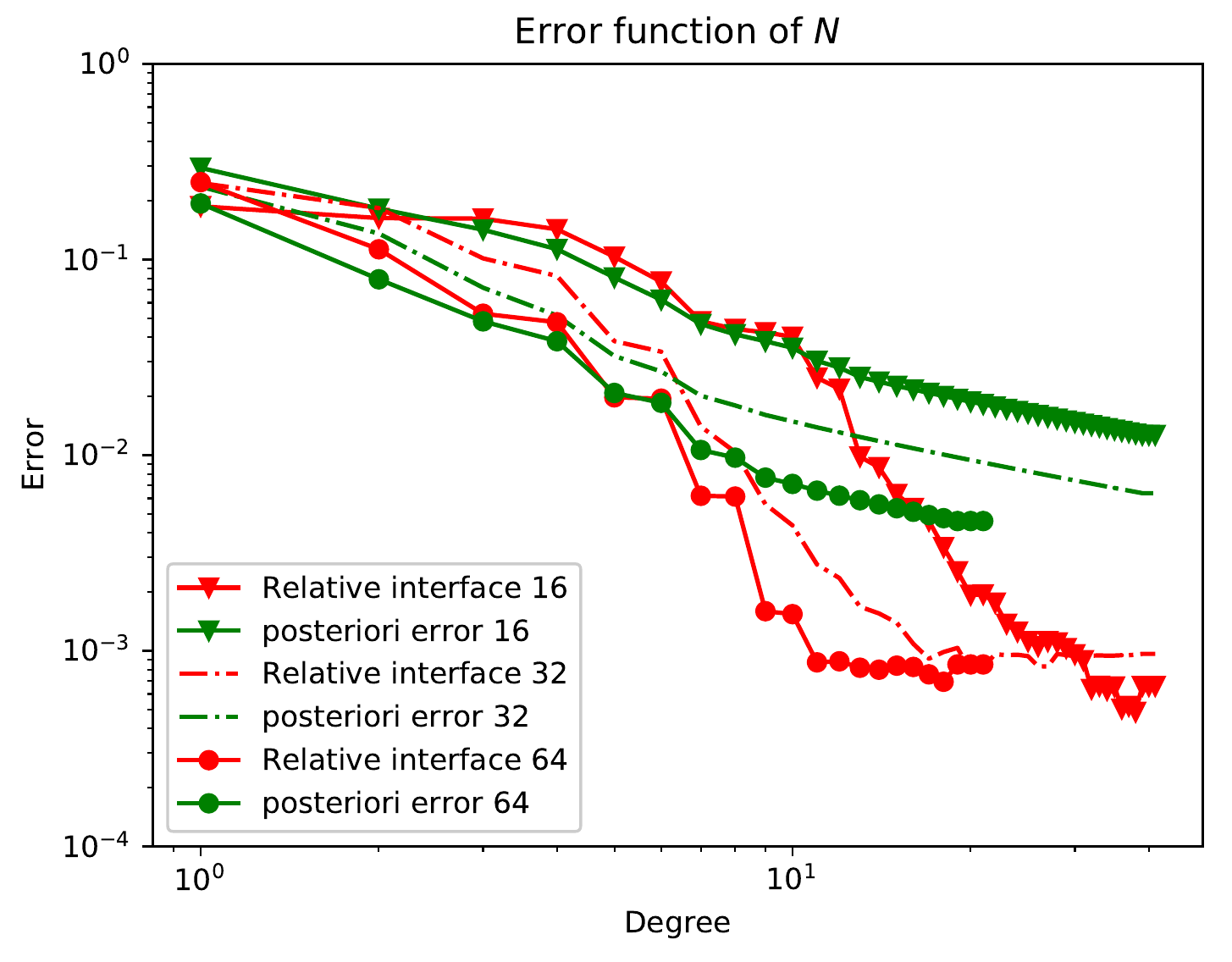}
  \caption{Left: \emph{a posteriori} error~\eqref{def:err_rel_gamma_post} and relative interface error~\eqref{def:err_rel_gamma} as a function of $N$ for $H=1/4$ and 1/8. Right: \emph{a posteriori} error~\eqref{def:err_rel_gamma_post} and relative interface error~\eqref{def:err_rel_gamma} as a function of $N$ for $H=1/16$, 1/32 and 1/64. \label{fig:res_N_posteriori}}
\end{figure}

\begin{figure}[htbp]
  \includegraphics[width=0.8\textwidth]{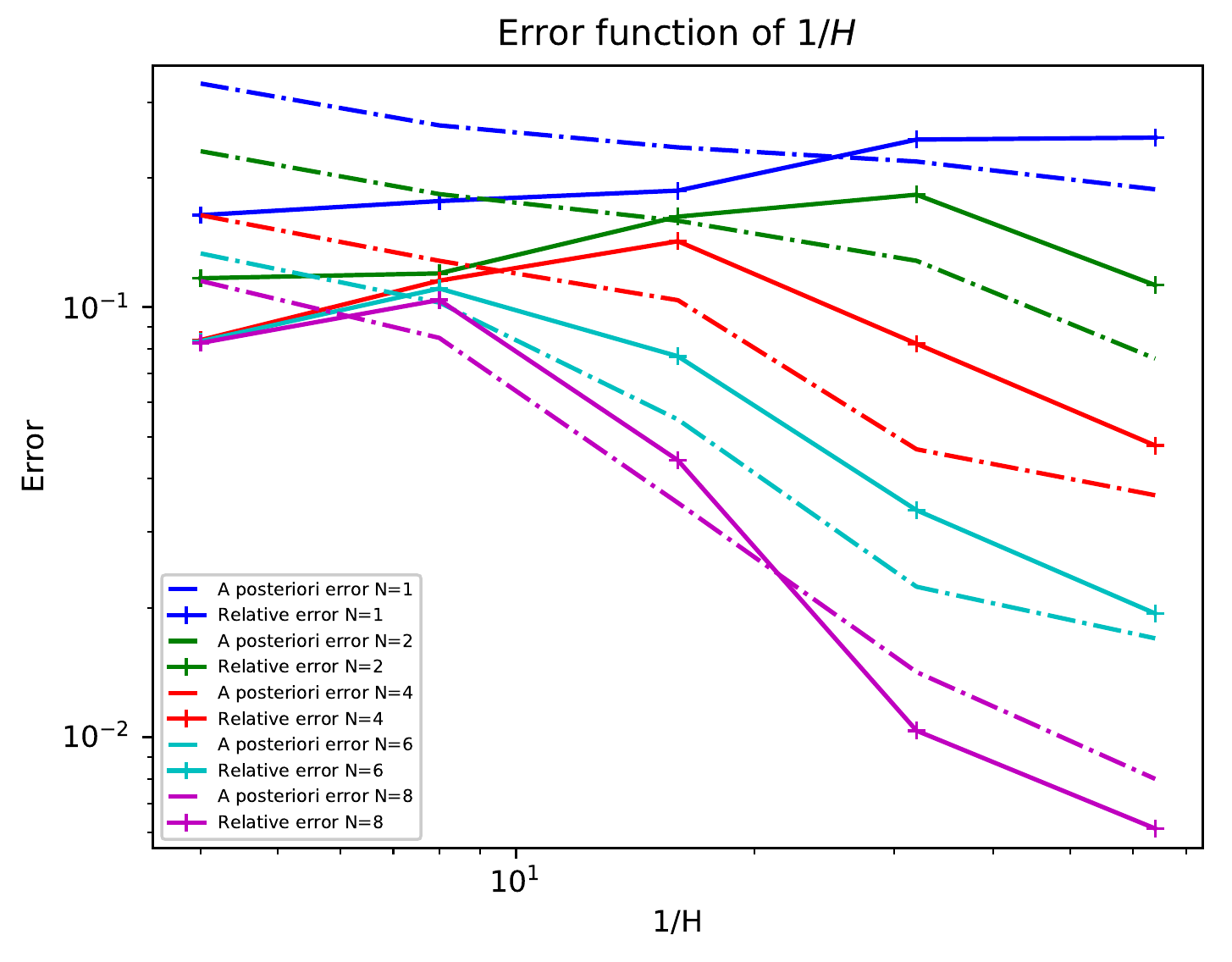}
  \caption{\emph{A posteriori} error~\eqref{def:err_rel_gamma_post} and relative interface error~\eqref{def:err_rel_gamma} as a function of $1/H$ for different polynomial degrees $N=\{1,2,4,6,8\}$. \label{fig:res_H_posteriori}}
\end{figure}

We consider here $\dps f(x,y) = -10 \exp \big[ -80\left((x-0.5)^2+(y-0.5)^2\right) \big]$ on $\Omega=(0,1)^2$, keep the definition~\eqref{eq:def_Aeps} for $A_\varepsilon$, and set $\varepsilon=1/32$. On the following figures, we compare the relative interface error~\eqref{def:err_rel_gamma} with the \emph{a posteriori} estimator~\eqref{def:err_rel_gamma_post} for several values of $N$ and $H$. Figures~\ref{fig:res_N_posteriori} and~\ref{fig:res_H_posteriori} show the behavior of both errors when $N$ increases and when $1/H$ increases, respectively.

In Figure~\ref{fig:res_N_posteriori}, we see that, for $H=1/4$ and $H=1/8$, the \emph{a posteriori} error estimator is an upper bound of the relative interface error and is a reliable indicator for any $N<10$. When $H=1/16$, $1/32$ and $1/64$, the \emph{a posteriori} error indicator seems to represent well the relative interface error for any $N \leq 10$. For higher polynomial degrees, the relative interface error decreases sharply and the \emph{a posteriori} indicator does not present such a behavior. The interpretation of the results for large values of $N$ is delicate, because several technicalities might affect the quality of these results: we do not know the energy of $u_\eps^\Gamma$, but only approximate it by the energy of some $u_\eps^{\Gamma,h}$ computed on a fine mesh; likewise, we only manipulate numerical approximations of the basis functions; third, when the difference of the energies is much smaller than the energies themselves, computing a relative error may become challenging.

We now turn to Figure~\ref{fig:res_H_posteriori}, which shows the behavior of the \emph{a posteriori} estimator when $H$ decreases for a fixed value of $N$. For any $H \leq 1/16$, the \emph{a posteriori} error seems to behave like the relative interface error for $N=\{1,2,4,6,8\}$. We can see that the \emph{a posteriori} estimator does not suffer from any resonance effect: it is decreasing with respect to $H$ for all values of $N$ tested. 

\medskip

\begin{figure}[htbp]
  \includegraphics[width=0.326\textwidth]{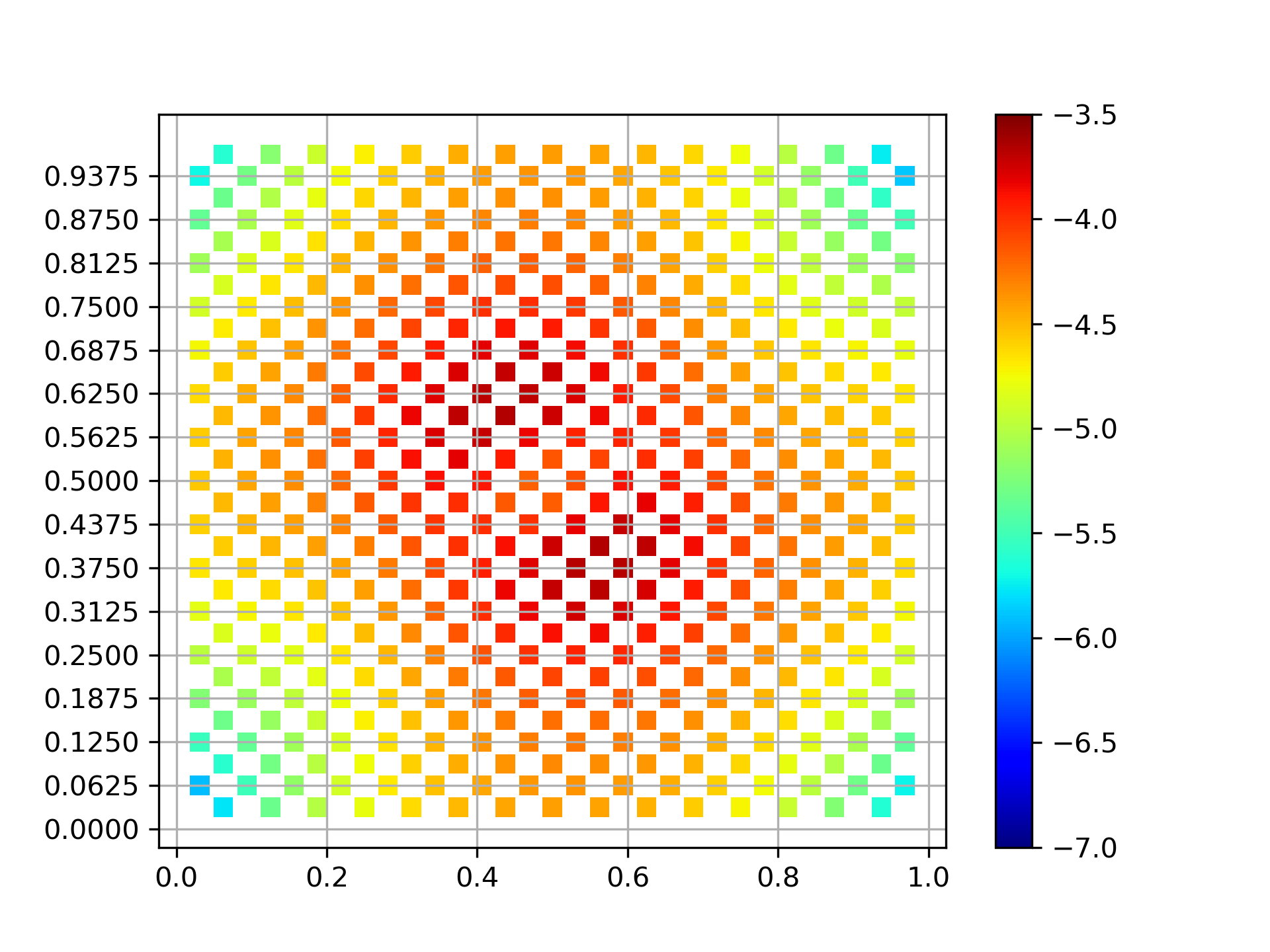}
  \includegraphics[width=0.326\textwidth]{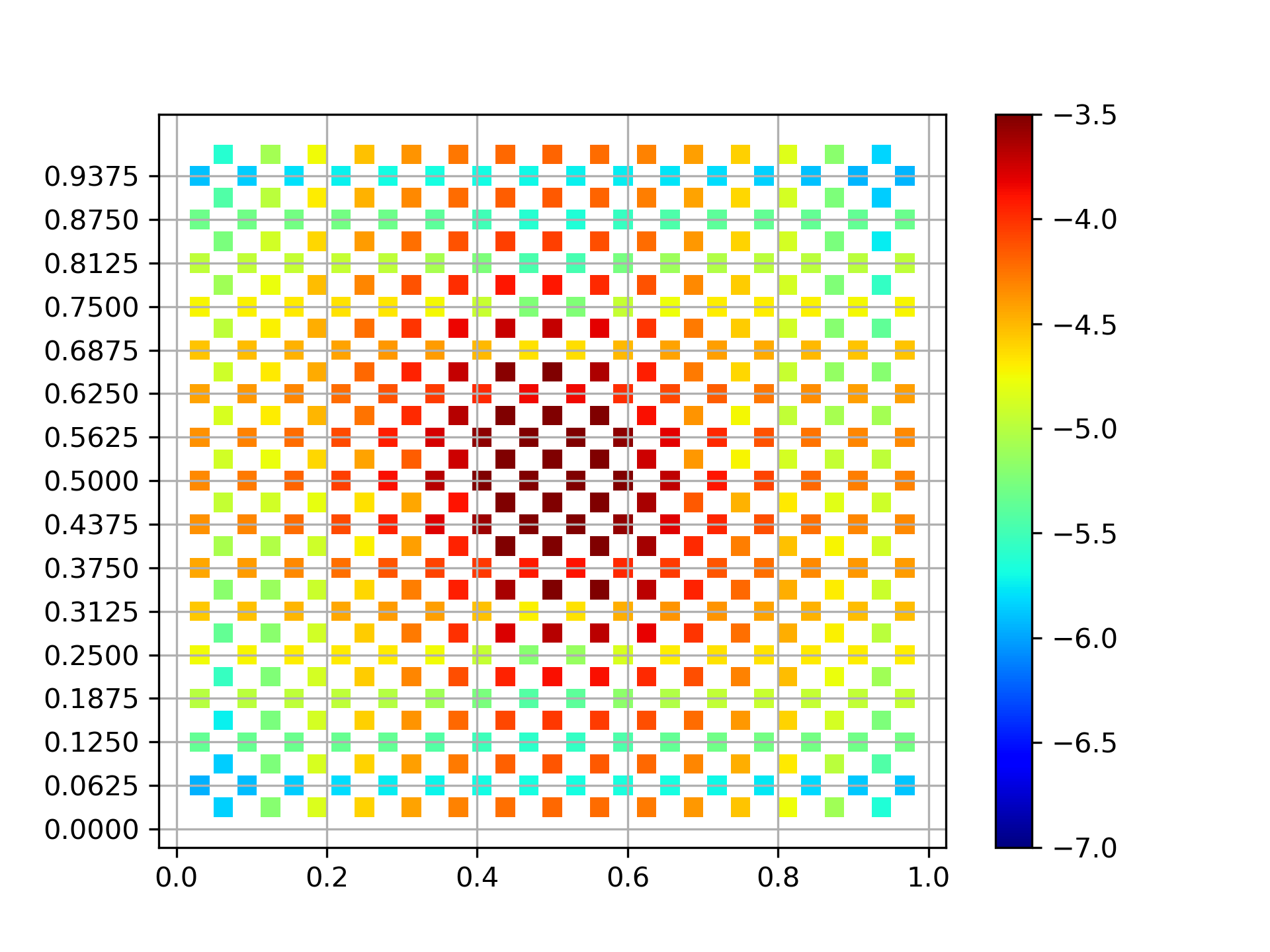}
  \includegraphics[width=0.326\textwidth]{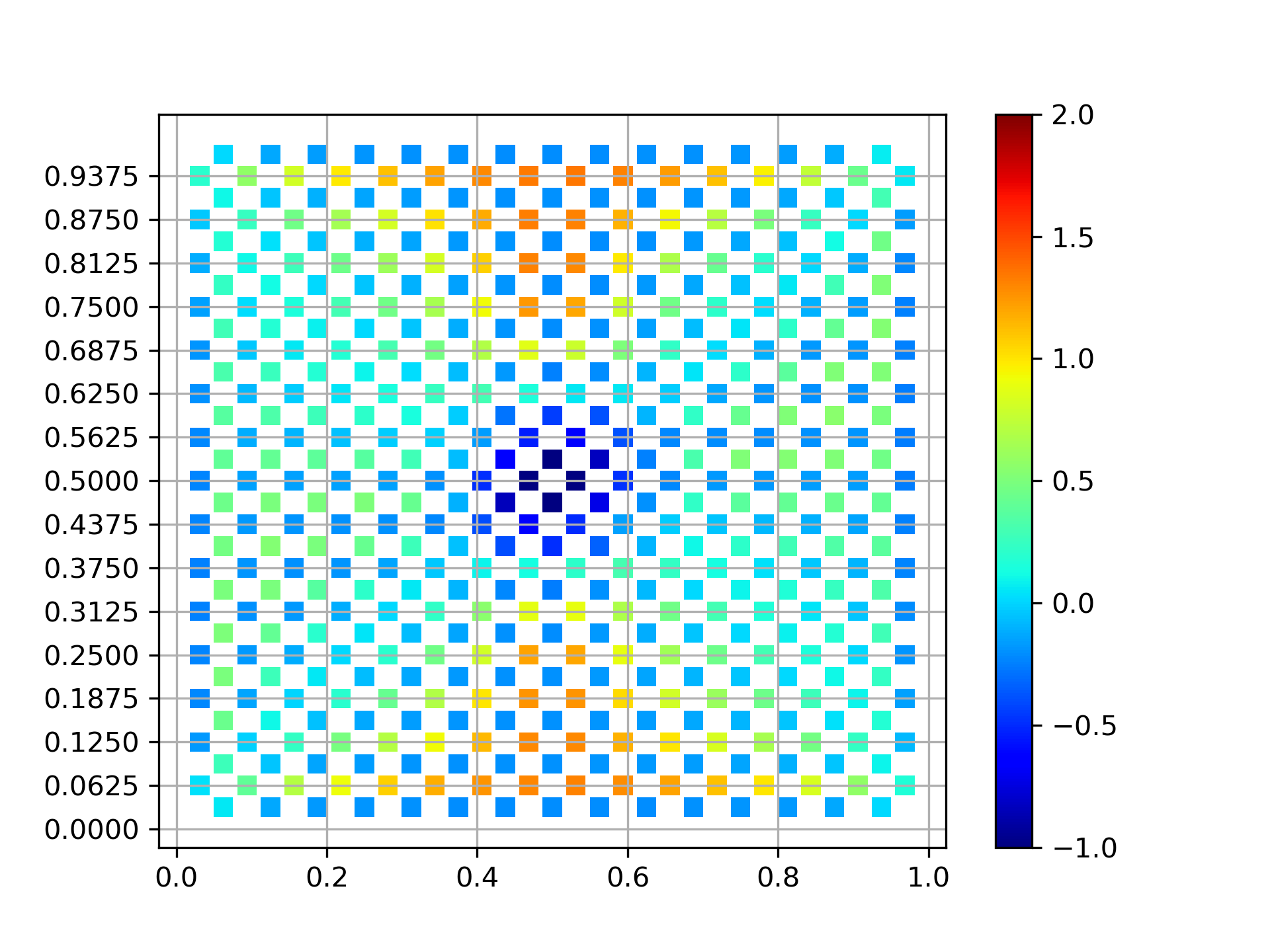} 
  \caption{Error maps edge by edge for $N=1$. Left: relative interface error; Center: \emph{a posteriori} estimator; Right: ratio of relative interface error and \emph{a posteriori} estimator (the plots are shown in a base-$10$ log scale). \label{fig:error_maps_N1}}
\end{figure}

\begin{figure}[htbp]
  \includegraphics[width=0.326\textwidth]{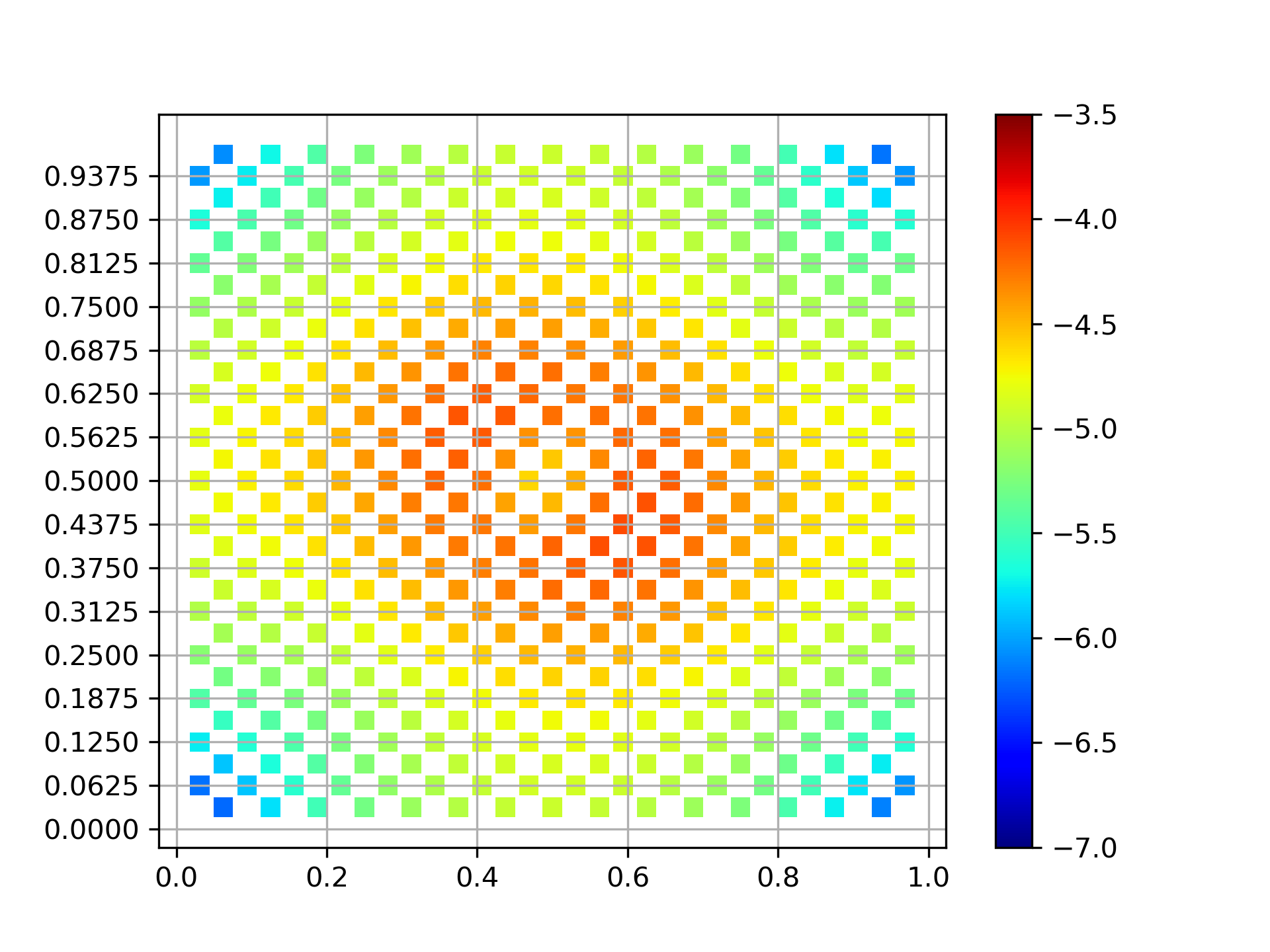}
  \includegraphics[width=0.326\textwidth]{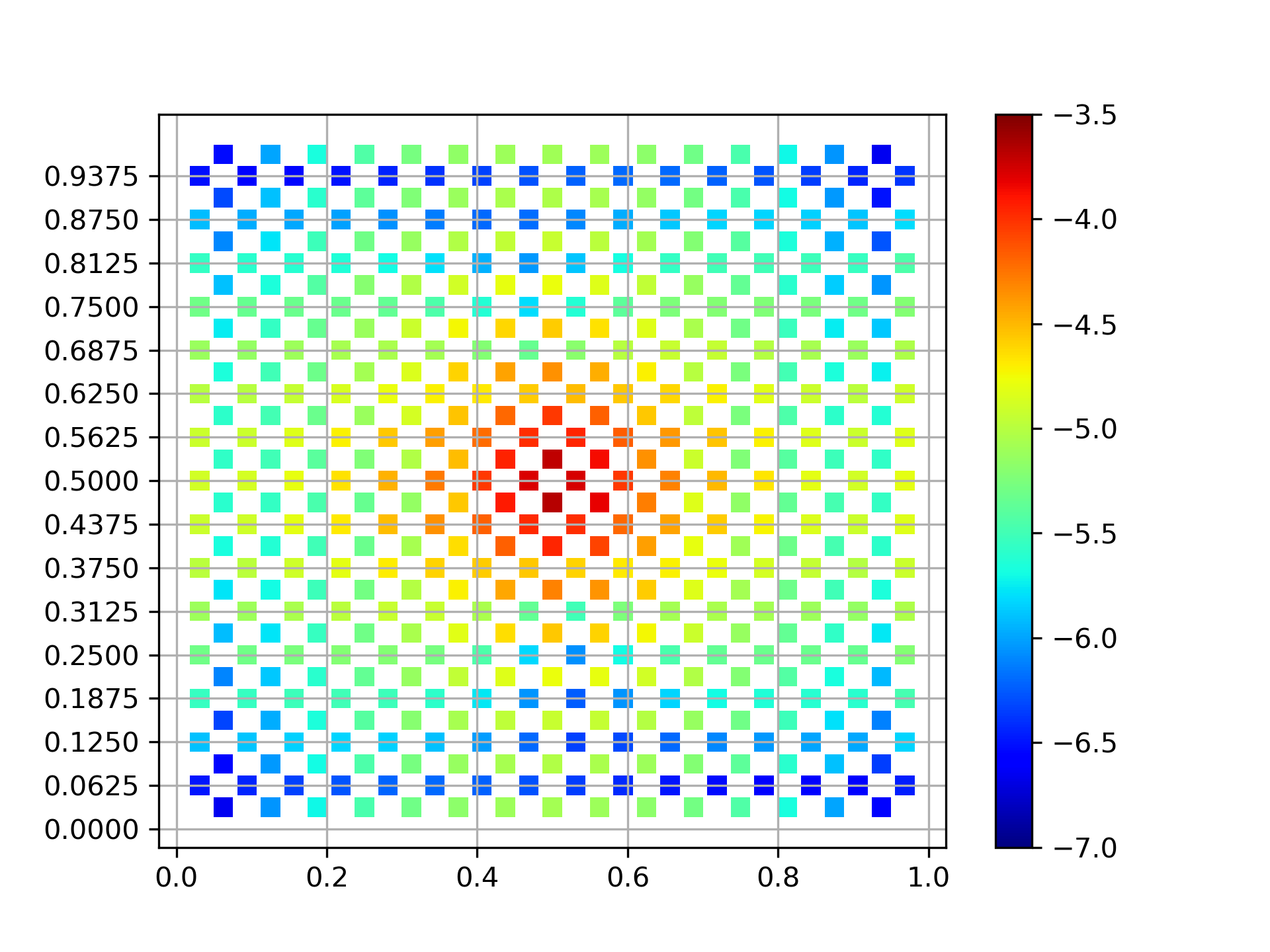}
  \includegraphics[width=0.326\textwidth]{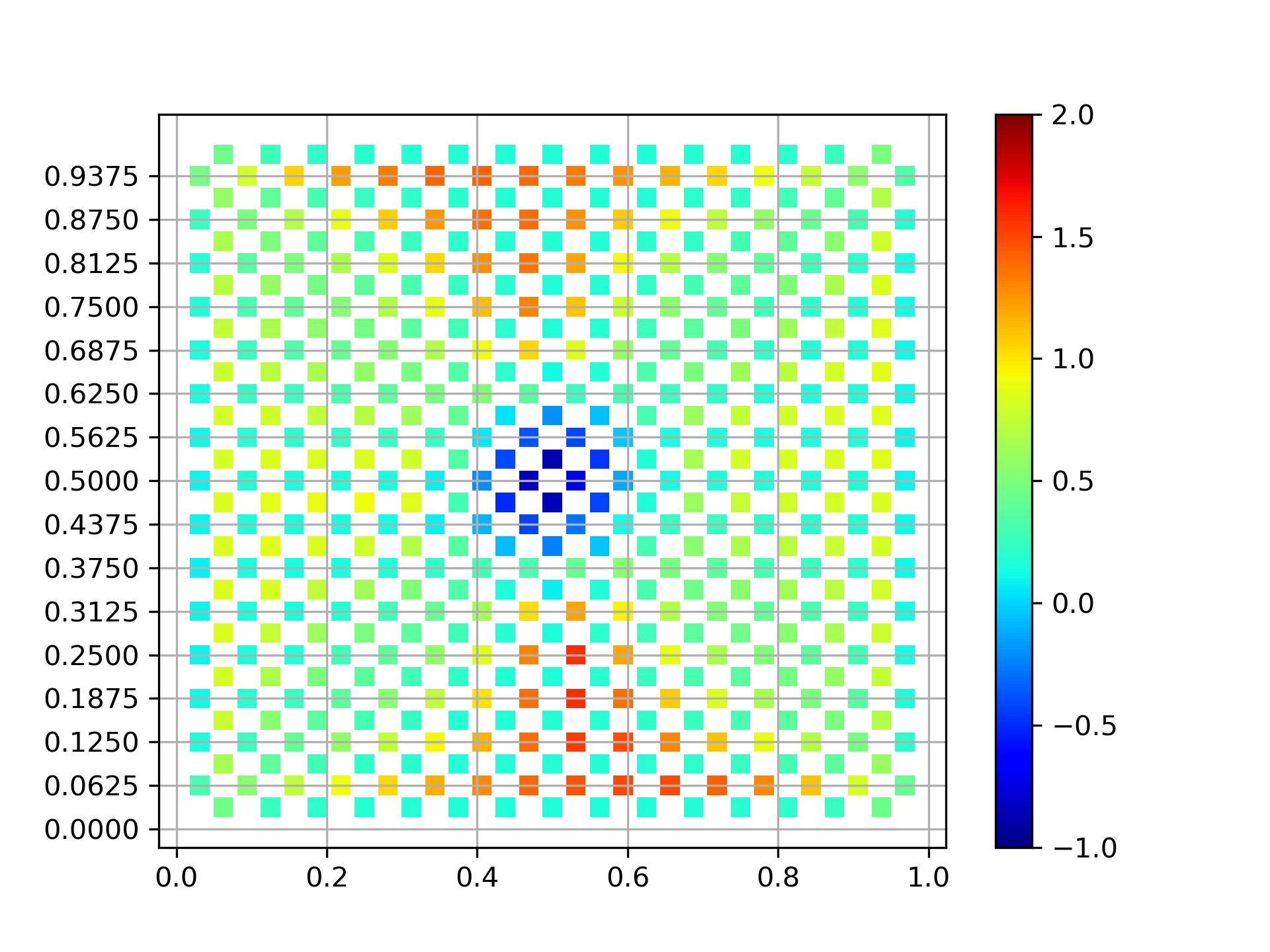} 
  \caption{Error maps edge by edge for $N=5$. Left: relative interface error; Center: \emph{a posteriori} estimator; Right: ratio of relative interface error and \emph{a posteriori} estimator (the plots are shown in a base-$10$ log scale). \label{fig:error_maps_N5}}
\end{figure}

\begin{figure}[htbp]			
  \includegraphics[width=0.326\textwidth]{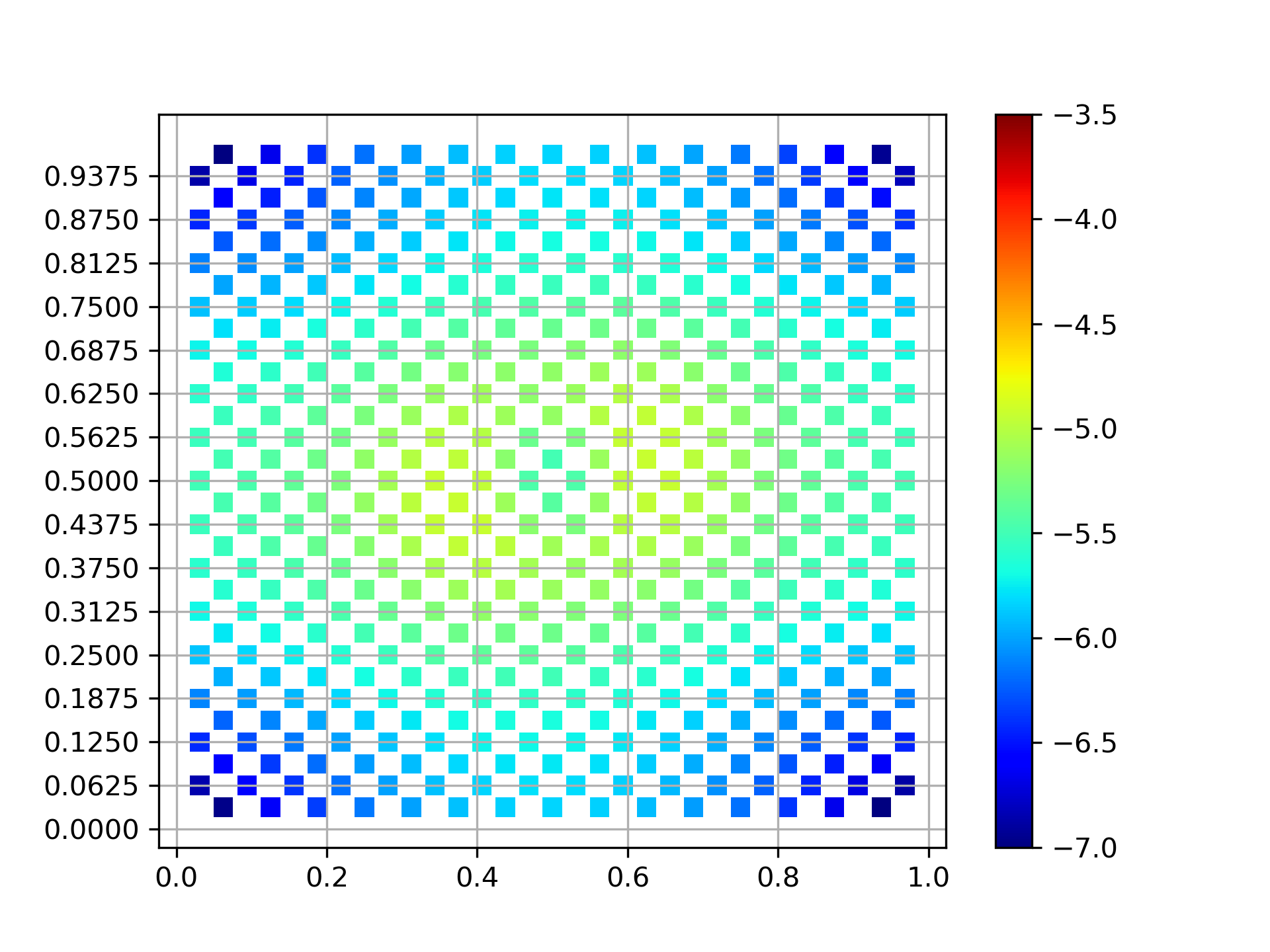}
  \includegraphics[width=0.326\textwidth]{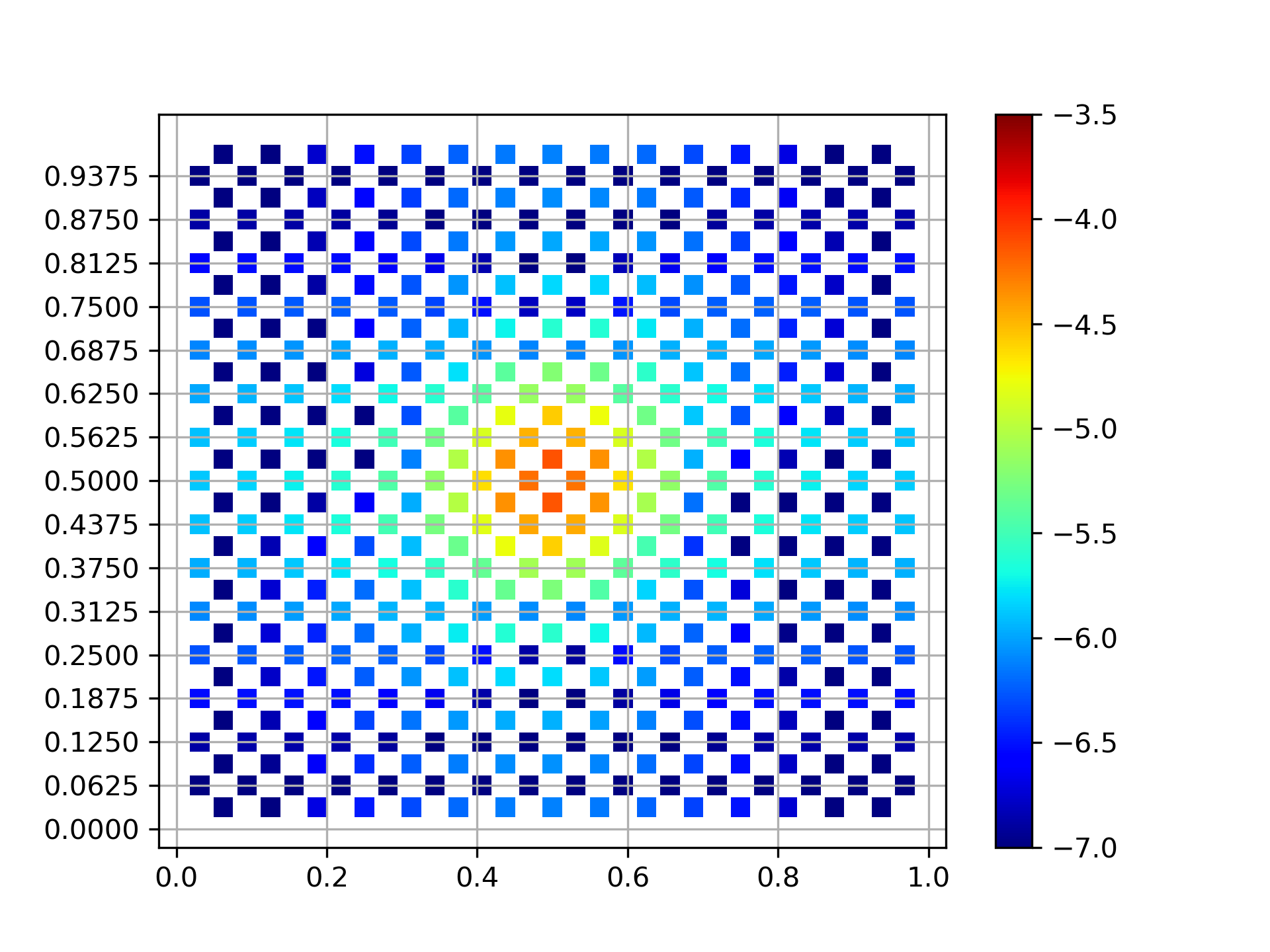}
  \includegraphics[width=0.326\textwidth]{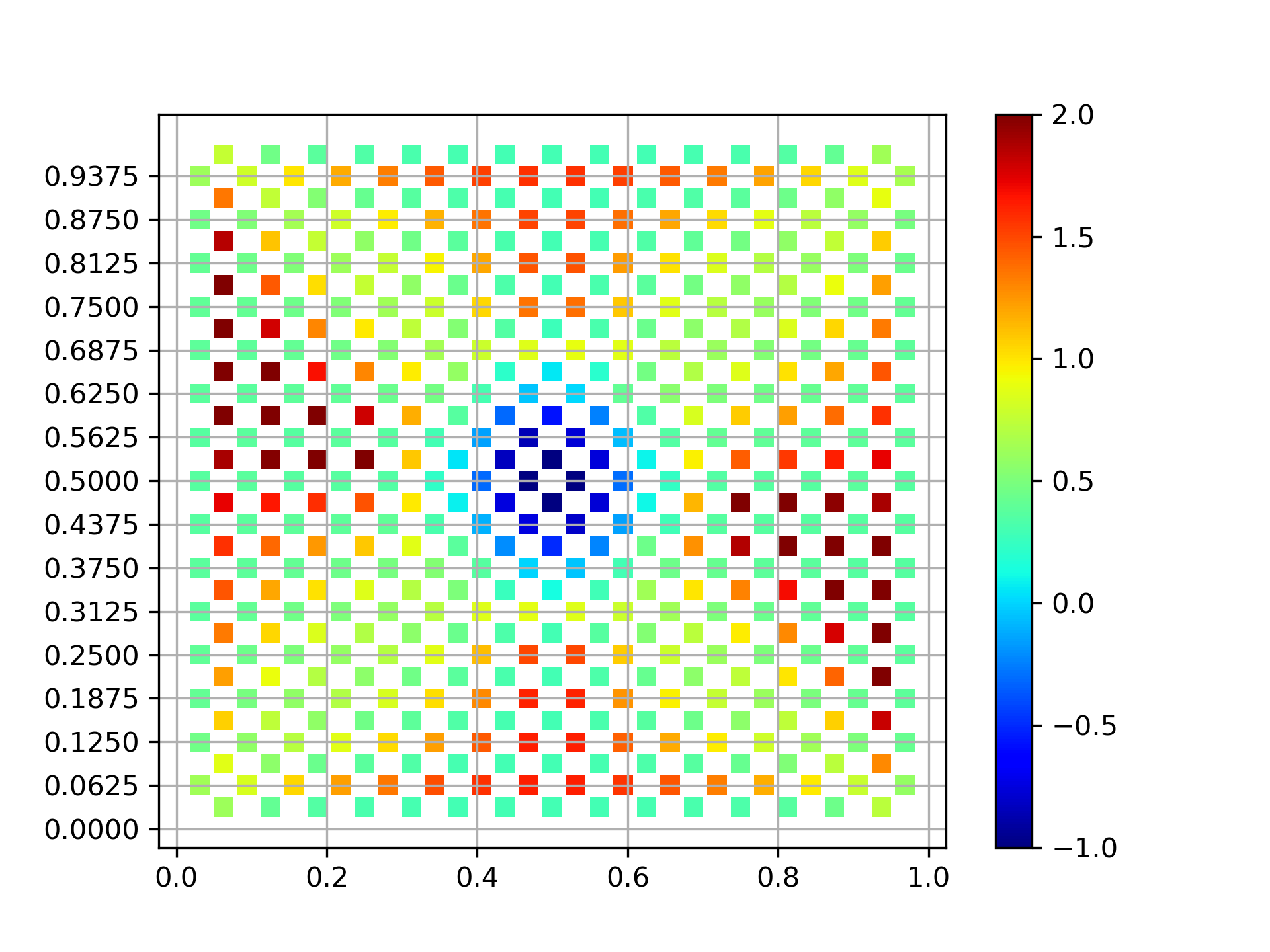} 
  \caption{Error maps edge by edge for $N=10$. Left: relative interface error; Center: \emph{a posteriori} estimator; Right: ratio of relative interface error and \emph{a posteriori} estimator (the plots are shown in a base-$10$ log scale). \label{fig:error_maps_N10}}
\end{figure}

One of the main interest of an \emph{a posteriori} estimator for which the error admits a local decomposition is precisely to allow for local refinement with respect to the parameters of the method: in our case, the polynomial degree $N_e$ of enrichments on each edge $e$ and the size $H_K$ of any element $K$. To that end, it is important to know whether the local behavior of the \emph{a posteriori} estimator represents well the local behavior of the actual error. This question is investigated on Figures~\ref{fig:error_maps_N1}, \ref{fig:error_maps_N5} and~\ref{fig:error_maps_N10}, where we show the error maps for $N=1$, 5 and 10 respectively (with $H=1/16$ fixed). We distribute the \emph{a posteriori} estimator~\eqref{def:err_rel_gamma_post} onto the edges. The first term of~\eqref{def:err_rel_gamma_post} is element based. For each edge, we therefore add the contributions of this first term associated to the two elements sharing the edge. The second term of~\eqref{def:err_rel_gamma_post} is simpler to handle since it is already edge based. Using such a localization procedure, we obtain an \emph{a posteriori} estimator which reads as a sum of contributions over the edges. Put differently, we write~\eqref{def:err_rel_gamma_post} as
$$
\mathcal{E}_{\rm post,\Gamma} = \sqrt{ \sum_{e\subset \Gamma} \big( \mathcal{E}_{\rm post,\Gamma}(e) \big)^2 }
$$
with
\begin{multline*}
\big( \mathcal{E}_{\rm post,\Gamma}(e) \big)^2
=
\frac{H_e}{p_e} \left\| J_e \big( \nu^T A \nabla(u_{\Gamma,H,\{N_e\}}) \big) \right\|_{L^2(e)}^2
\\
+
\sum_{K \in \mathcal{T}_H, \ e \subset \partial K} \frac{1}{\beta_K} \, \|f\|_{L^2(K)}^2 \left( \sum_{\widetilde{e} \subset \partial K} \frac{H_{\widetilde{e}} \, H_K}{N_{\widetilde{e}} \, p_{\widetilde{e}}} \right),
\end{multline*}
where $\beta_K$ is the number of edges of $K$ which belong to $\Gamma$ (for a quadrangular mesh for instance, $\beta_K=4$ for all elements $K$ except those at the boundary of $\Omega$). Note of course that the way we redistribute the error among the edges is arbitrary. Other choices could be made. This is a standard difficulty for residue-type estimators consisting of both element-based and edge-based contributions (see e.g.~\cite{CAR00}).

We plot on Figures~\ref{fig:error_maps_N1}, \ref{fig:error_maps_N5} and~\ref{fig:error_maps_N10} the resulting values $\mathcal{E}_{\rm post,\Gamma}(e)$ on a $16 \times 16$ coarse mesh. Regarding the actual error, we compute the relative energy error~\eqref{def:err_rel_gamma} elements by elements by comparing, on each element, the reference solution with the numerical approximation. Similarly to the first term of~\eqref{def:err_rel_gamma_post}, we may write the numerator of~\eqref{def:err_rel_gamma}, which is a quantity by definition distributed on all elements, as a sum of contributions among the edges (the denominator of~\eqref{def:err_rel_gamma} is kept unchanged and is never localized). On the right side of Figures~\ref{fig:error_maps_N1}, \ref{fig:error_maps_N5} and~\ref{fig:error_maps_N10}, we plot the error map showing the ratio between the actual local error on an edge $e$ and the local \emph{a posteriori} estimator $\mathcal{E}_{\rm post,\Gamma}(e)$.

When $N=1$ and $N=5$ (see the right plot on Figures~\ref{fig:error_maps_N1} and~\ref{fig:error_maps_N5}), we can see that the ratio between the local actual error and the local \emph{a posteriori} estimator does not significantly change over the domain $\Omega$. To be more precise, there is actually an exception near the center of $\Omega$, that may be due to the fact that $|\nabla u|$ is very small there, and thus challenging to approximate with a good accuracy, even when using a relative error. Since the ratio of actual error {\em vs} estimated error is close to a constant, it is thus possible to use the local \emph{a posteriori} estimator to drive an adaptive discretization procedure: the edges of the quadrangular mesh where $\mathcal{E}_{\rm post,\Gamma}(e)$ is large are indeed the edges where the actual error is large.

In contrast, when $N=10$ (see the right plot on Figure~\ref{fig:error_maps_N10}), the ratio between the actual and the predicted error widely varies over the domain $\Omega$. This is consistent with the above Figure~\ref{fig:res_N_posteriori} showing a significant difference between the global actual error and the global estimated error for large values of $N$. For this large value of $N$, the quantity $\mathcal{E}_{\rm post,\Gamma}(e)$ cannot reasonably be used to drive a reliable adaptation procedure.

\medskip

The above numerical tests hence show that the \emph{a posteriori} estimator defined in Proposition~\ref{prop:error_indicator} correctly represents the behavior of the actual error for $N<10$ and in the regime where $H$ is close to or slightly smaller than $\varepsilon$ (a regime still relevant for numerical multiscale approaches, since a classical P1 approach would need to take $H$ much smaller than $\eps$ to be accurate). In such a regime, the estimator can thus trustfully be used to locally refine the polynomial degree $N_e$ associated with the edge $e$ and the size $H_K$ of the element $K$. 

\appendix

\section{Proofs} 

This section is devoted to the proofs of Lemmas~\ref{lem:Error_bubble} and~\ref{lem:Error_Interface} and of Proposition~\ref{prop:error_indicator}. These proofs heavily rely on technical results about polynomial approximation theory, fractional Sobolev spaces and elliptic regularity, that we recall below as we proceed. Throughout this appendix, the constant denoted $C$ may change from one line to the next. When valid, the independence of that constant with respect to the mesh size and other quantities will always be underlined.

\subsection{Proof of Lemma~\ref{lem:Error_bubble}} \label{sec:proofs_un}
 
The proof of Lemma~\ref{lem:Error_bubble} makes use of the following approximation result, which is shown in~\cite[Equation~(5.8.27) p.~318]{canuto2010spectral} for the case of quadrangles and in~\cite[Section~5.9]{canuto2010spectral} for the case of triangles.

\begin{lemma} \label{lem:approx_poly}
Assume that $\left(\mathcal{T}_H\right)_H$ is a family of conformal meshes of $\Omega$ composed of a finite number of convex quadrangles (resp. triangles) with straight edges. Assume also that the meshes are regular in the sense of~\eqref{eq:def_mesh_regular}. For any quadrangle $K$ (resp. triangle $K$), let $\Pi_M^K$ be the $L^2(K)$-orthogonal projection on the vector space of polynomials of degree in each variable (resp. total degree) at most $M$. Then, for any non-negative integer $\ell$, there exists $C_\ell$ independent of $H$, $M$ and of the elements $K$ of the family of meshes such that, for any $v \in H^\ell(K)$,
\begin{equation} \label{eq:approx_poly}
\left\| v - \Pi_M^K(v) \right\|_{L^2(K)} \leq C_\ell \, \frac{H^{\min(\ell,M+1)}}{M^\ell} \, \| v \|_{H^\ell(K)}.
\end{equation}
\end{lemma}

\medskip

\noindent
{\bf Proof of Lemma~\ref{lem:Error_bubble}.} 
We first show that $\|u_B\|_E \leq  C H \|f\|_{L^2(\Omega)}/ \sqrt{\alpha_{\rm min}}$. We have that
$$
  \|u_B\|_E^2
  =
  \int_\Omega (\nabla u_B)^T A \nabla u_B
  =
  \int_\Omega f u_B
  =
  \sum_{K \in \mathcal{T}_H} \int_K f u_B.
$$
Using the Cauchy-Schwarz inequality and the Poincar\'e inequality (recall indeed that $u_B \in H^1_0(K)$ for any $K$), it holds that 
$$
\int_K f u_B \leq \| f \|_{L^2(K)} \, \|u_B\|_{L^2(K)} \leq C \, H \, \|f\|_{L^2(K)} \, |u_B|_{H^1(K)},
$$
for some universal constant $C$, where we recall that $|u_B|_{H^1(K)} = \| \nabla u_B \|_{L^2(K)}$ is the $H^1$ semi-norm on $K$. We hence have 
\begin{align*}
  \|u_B\|_E^2
  &\leq
  C H \sum_{K \in \mathcal{T}_H} \|f\|_{L^2(K)} |u_B|_{H^1(K)}
  \\
  & \leq
  CH \|f\|_{L^2(\Omega)} |u_B|_{H^1(\Omega)} 
  \\
  & \leq
  \frac{CH}{\sqrt{\alpha_{\rm min}}} \|f\|_{L^2(\Omega)} \|u_B\|_E, 
\end{align*}
from which we deduce that $\dps \|u_B\|_E \leq \frac{CH}{\sqrt{\alpha_{\rm min}}} \|f\|_{L^2(\Omega)}$ for some universal constant $C$. This proves~\eqref{eq:error_u_bubble0}, where we recall that $u_{B,H,M} = 0$.

\medskip

We now add bubble enrichments for each element $K \in \mathcal{T}_H$ with a uniform degree $M \geq 1$. For any $v_{B,H,M} \in V_{B,H,M}$, we have
\begin{align*}
  a(u_B - v_{B,H,M},u_B - v_{B,H,M})
  &=
  a(u_B - u_{B,H,M},u_B - u_{B,H,M})
  \\
  & + a(u_{B,H,M} - v_{B,H,M},u_{B,H,M} - v_{B,H,M})
  \\
  & + 2a(u_B - u_{B,H,M},u_{B,H,M} - v_{B,H,M}).
\end{align*}
The third term above vanishes in view of~\eqref{eq:approx_bubble}. The second term is non-negative since $a$ is coercive. We thus obtain that, for any $v_{B,H,M} \in V_{B,H,M}$,
\begin{equation} \label{z}
  a(u_B - u_{B,H,M}, u_B - u_{B,H,M}) \leq a(u_B - v_{B,H,M}, u_B - v_{B,H,M}).
\end{equation}
Recall that $V_{B,H,M}$ is the span of the functions $\{ \phi_{K,i}^B \}_{i=1,\dots,{\cal N}_M}$ (see~\eqref{eq:def_VBHM}) that solve the problem $-\div (A \nabla \phi_{K,i}^B) = P_i$ in the element $K$ with the boundary condition $\phi_{K,i}^B = 0$ on $\partial K$, where $\{ P_i \}_{i=1,\dots,{\cal N}_M}$ is a basis of polynomial functions with degree at most $M$.

Let $\dis \Pi_M(f) = \sum_{K \in \mathcal{T}_H} 1_K \, \Pi_M^K(f)$, where $\Pi_M^K$ is defined in Lemma~\ref{lem:approx_poly}. The unique $z \in V_B$ such that $\dis a(z,v) = \int_\Omega \Pi_M(f) \, v$ for any $v \in V_B$ can be written as a linear combination of the functions $\{ \phi_{K,i}^B \}_{i=1,\dots,{\cal N}_M, \ K \in \mathcal{T}_H}$. It thus belongs to $V_{B,H,M}$, and we denote it $v_{B,H,M} \in V_{B,H,M}$, which thus satisfies
$$
\forall v \in V_B, \quad a(v_{B,H,M},v) = \int_\Omega \Pi_M(f) \, v.
$$
Using the definition of $u_B$, we next obtain that, for any $v\in V_B$,
\begin{align*}
a(u_B - v_{B,H,M},v)
&=
\int_\Omega f v - \int_\Omega \Pi_M(f) \, v
\nonumber
\\
&=
\sum_{K \in \mathcal{T}_H} \int_K \left(f-\Pi_M^K(f)\right) v
\nonumber
\\
&=
\sum_{K \in \mathcal{T}_H} \int_K \left(f-\Pi_M^K(f)\right) \left(v-\Pi_M^K(v)\right).
\end{align*}
Choosing now $v=u_B-v_{B,H,M}$ in the above equality yields
\begin{align*}
& a(u_B - v_{B,H,M},u_B - v_{B,H,M})
\\
& \leq
\sum_{K \in \mathcal{T}_H} \left\| f-\Pi_M^K(f) \right\|_{L^2(K)} \ \left\| u_B - v_{B,H,M} -\Pi_M^K\left(u_B - v_{B,H,M} \right) \right\|_{L^2(K)}
\\
& \leq
C_\ell \, \frac{H^{\min(\ell,M+1)+1}}{M^{\ell+1}} \sum_{K \in \mathcal{T}_H} \| f \|_{H^\ell(K)} \ \| u_B - v_{B,H,M} \|_{H^1(K)},
\end{align*}
where we have used in the last line, for each of the two factors, the polynomial projection properties stated in Lemma~\ref{lem:approx_poly} (using that $f \in H^\ell(K)$ for the first factor and that $u_B - v_{B,H,M} \in H^1(K)$ for the second factor). Using a Poincar\'e inequality in $H^1_0(K)$ and that $H \leq 1$, we obtain $\| u_B - v_{B,H,M} \|_{H^1(K)} \leq \overline{C} \, (1+H) \, \| \nabla (u_B - v_{B,H,M}) \|_{L^2(K)} \leq C \, \| \nabla (u_B - v_{B,H,M}) \|_{L^2(K)}$. We thus deduce that
\begin{align*}
& a(u_B - v_{B,H,M},u_B - v_{B,H,M} )
\\
& \leq
C_\ell \, \frac{H^{\min(\ell,M+1)+1}}{M^{\ell+1}} \sum_{K \in \mathcal{T}_H} \| f \|_{H^\ell(K)} \ \| \nabla (u_B - v_{B,H,M}) \|_{L^2(K)}
\\
& \leq
C_\ell \, \frac{H^{\min(\ell,M+1)+1}}{M^{\ell+1}} \, \| f \|_{H^\ell(\Omega)} \, \| \nabla (u_B - v_{B,H,M}) \|_{L^2(\Omega)}
\\
& \leq
\frac{C_\ell}{\sqrt{\alpha_{\rm min}}} \, \frac{H^{\min(\ell,M+1)+1}}{M^{\ell+1}} \, \| f \|_{H^\ell(\Omega)} \, \sqrt{a(u_B - v_{B,H,M}, u_B - v_{B,H,M})},
\end{align*}
where we have used a discrete Cauchy-Schwarz inequality in the third line and the lower bound on $A$ in the last line. We hence obtain
\begin{equation} \label{zz}
\sqrt{a(u_B - v_{B,H,M}, u_B - v_{B,H,M})} \leq \frac{C_\ell}{\sqrt{\alpha_{\rm min}}} \, \frac{H^{\min(\ell,M+1)+1}}{M^{\ell+1}} \, \| f \|_{H^\ell(\Omega)}.
\end{equation}
Inserting~\eqref{zz} into~\eqref{z}, we obtain~\eqref{eq:error_u_bubble}. This concludes the proof of Lemma~\ref{lem:Error_bubble}.
\qed

\subsection{Proof of Lemma~\ref{lem:Error_Interface}} \label{sec:proofs_deux}

For the proof of Lemma~\ref{lem:Error_Interface}, we separately consider the case of quadrangles and the case of triangles. For the former case, we need the following approximation result (see in~\cite{canuto2010spectral} estimates~(5.8.26) and~(5.8.27) and the discussion following the latter).

\begin{lemma} \label{lem:inter_poly}
Assume that $\left(\mathcal{T}_H\right)_H$ is a family of conformal meshes of $\Omega$ composed of a finite number of convex quadrangles with straight edges, and that the meshes are regular in the sense of~\eqref{eq:def_mesh_regular}. For any quadrangle $K$, let $i_N^K$ be the Legendre interpolant at the $(1+N)^2$ Gauss-Lobatto points in $K$ ($i_N^K$ is thus a polynomial function in $\mathbb{Q}_N$). Let $s \geq 3/2$ and $N \geq 1$. Then there exists $C_s$ independent of $H$, $N$ and of the elements $K$ of the family of meshes such that, for any $v \in H^s(K)$,
$$
\left| v - i_N^K(v) \right|_{H^1(K)} \leq C_s \ \frac{H^{\min(s,N+1)-1}}{N^{s-1}} \ \| v \|_{H^s(K)},
$$
where we recall that $| \cdot |_{H^1(K)}$ is the $H^1$ semi-norm on $K$.
\end{lemma}

Note that a function $v \in H^s(K)$ with $s>1$ is continuous (recall that we consider a two-dimensional setting), thus $i_N^K(v)$ is well-defined.

\medskip

\noindent
{\bf Proof of Lemma~\ref{lem:Error_Interface} for quadrangles.} 
Using arguments similar to those used to prove~\eqref{z}, we have, for any $v_{\Gamma,H,N} \in V_{\Gamma,H,N}$,
\begin{equation} \label{eq:jeudi1}
a(u_\Gamma - u_{\Gamma,H,N}, u_\Gamma - u_{\Gamma,H,N})
\leq
a(u_\Gamma - v_{\Gamma,H,N}, u_\Gamma - v_{\Gamma,H,N}).
\end{equation}
Since $u \in H^s(\Omega)$ with $s > 1$, $u$ is continuous on $\Omega$, hence on $\Gamma$. We denote by $i^\Gamma_N(u)$ the interpolant (at the Gauss-Lobatto points on $\Gamma$) of $u$ on the set of continuous functions on $\Gamma$ which are piecewise equal to polynomial functions of degree lower than or equal to $N$.

Let $w = E_\Omega(i^\Gamma_N(u)) \in V_{\Gamma,H,N}$ denote the harmonic lifting of $i^\Gamma_N(u)$, that is the solution to $-\div (A \nabla w)=0$ on each coarse element $K$ with the Dirichlet boundary conditions $w = i^\Gamma_N(u)$ on $\Gamma$. 
We then have
\begin{align}
& a\Big(u_\Gamma - E_\Omega(i^\Gamma_N(u)), u_\Gamma - E_\Omega(i^\Gamma_N(u))\Big)
\nonumber
\\
& =
\sum_{K \in \mathcal{T}_H} \int_K \big( \nabla u_\Gamma - \nabla E_\Omega(i^\Gamma_N(u)) \big)^T A \big( \nabla u_\Gamma - \nabla E_\Omega(i^\Gamma_N(u)) \big)
\nonumber
\\
& \leq
\sum_{K \in \mathcal{T}_H} \int_K \big( \nabla u - \nabla I_{\Gamma,H,N}(u) \big)^T A \big( \nabla u - \nabla I_{\Gamma,H,N}(u) \big)
\nonumber
\\
& \leq 
\alpha_{\rm max} \sum_{K \in \mathcal{T}_H} \left| u - I_{\Gamma,H,N}(u) \right|_{H^1(K)}^2,
\label{eq:jeudi2}
\end{align}
where $I_{\Gamma,H,N}(u) \in H^1_0(\Omega)$ is defined piecewise on each $K$ as the Legendre interpolant of $u |_K$ at the Gauss-Lobatto points in $K$ (it is thus a polynomial function in $\mathbb{Q}_N$). The first inequality of~\eqref{eq:jeudi2} holds for the following three reasons:
\begin{itemize}
\item First, $E_\Omega(i^\Gamma_N(u))$ and $I_{\Gamma,H,N}(u)$ agree on $\Gamma$ for quadrangular mesh elements. Recall indeed that the Gauss-Lobatto points of each edge of $\partial K$ are a subset of the Gauss-Lobatto points of $K$. On each edge, $E_\Omega(i^\Gamma_N(u))$ and $I_{\Gamma,H,N}(u)$ are thus two polynomial functions of degree lower than or equal to $N$ which are equal on the $(1+N)$ Gauss-Lobatto points of the edge, and are thus equal.
\item Second, $u_\Gamma$ and $u$ agree on $\Gamma$, by definition of $u_\Gamma$.
\item Third, we observe, using the short-hand notation $v_1 = u - I_{\Gamma,H,N}(u)$ and $v_2 = u_\Gamma - E_\Omega(i^\Gamma_N(u))$, that
\begin{align*}
  & \int_K (\nabla v_1)^T A \nabla v_1
  \\
  &=
  \int_K (\nabla v_2)^T A \nabla v_2 + \int_K (\nabla (v_1-v_2))^T A \nabla (v_1-v_2) + 2 \int_K (\nabla (v_1-v_2))^T A \nabla v_2
  \\
  & \geq
  \int_K (\nabla v_2)^T A \nabla v_2,
\end{align*}
where we have used that, in the second line, the second term is non-negative and the third term vanishes (since $v_1-v_2=0$ on $\partial K$ and $-\div (A \nabla v_2)=0$ in $K$). This shows the first inequality of~\eqref{eq:jeudi2}.
\end{itemize}  
Using Lemma~\ref{lem:inter_poly} (where the Legendre interpolant $i_N^K(u)$ there is denoted $I_{\Gamma,H,N}(u)$ here), we see that
$$
\left( \sum_{K \in \mathcal{T}_H} \left| u - I_{\Gamma,H,N}(u) \right|_{H^1(K)}^2 \right)^{1/2}
\leq 
C_s \ \frac{H^{\min(s,N+1)-1}}{N^{s-1}} \ \| u \|_{H^s(\Omega)}.
$$
Collecting this bound with~\eqref{eq:jeudi1} (that we use for $v_{\Gamma,H,N} = E_\Omega(i^\Gamma_N(u))$, which indeed belongs to $V_{\Gamma,H,N}$) and~\eqref{eq:jeudi2}, we deduce~\eqref{eq:error_u_gamma_quad}, which concludes the proof of Lemma~\ref{lem:Error_Interface} for quadrangles.
\qed

\medskip

We now turn to the proof of Lemma~\ref{lem:Error_Interface} for the case of triangles. The proof given below actually also holds in the case of quadrangles. We have however kept the above proof specific to the case of quadrangles because the choice of $v_{\Gamma,H,N}$ is therein constructive, in contrast to the proof below. The following result plays in the general case the role of Lemma~\ref{lem:inter_poly} in the case of quadrangles.

\begin{lemma}[see proof of Theorem~4.6 of~\cite{babuvska1987hp}]
\label{lem:Babuska}
Consider a conformal mesh $\mathcal{T}_H$, regular in the sense of~\eqref{eq:def_mesh_regular}, and composed of triangular (resp. quadrangular) elements with meshsize $H$. Let $u \in H^s(\Omega) \cap H^1_0(\Omega)$ with $s>3/2$. Let $V^N_{H,0} = \{ v \in C^0(\overline{\Omega}) \cap H^1_0(\Omega); \ v|_K \in P_N^K\}$ where $P_N^K$ is the set of polynomial functions on $K$ that are of partial (resp. total) degree lower than or equal to $N$. We then have
$$
\min_{v \in V^N_{H,0}} \|u-v\|_{H^1(\Omega)}\leq C_s \frac{H^{\min(s,N+1)-1}}{N^{s-1}} \|u\|_{H^s(\Omega)},
$$
where $C_s$ is independent of $H$, $N$ and $u$.
\end{lemma}

\medskip

\begin{remark}
Theorem~4.6 in~\cite{babuvska1987hp} considers the problem of approximating the solution $u \in H^1_0(\Omega)$ to $-\Delta u + u = f$ in $\Omega$. Its proof relies on some approximation results, which have their own interest and can be stated as in Lemma~\ref{lem:Babuska} for any $u \in H^s(\Omega) \cap H^1_0(\Omega)$. 
\end{remark}

\noindent
{\bf Proof of Lemma~\ref{lem:Error_Interface} for triangles (and alternative proof for quadrangles).}
Using arguments similar to those used to prove~\eqref{z}, we have, for any $v_{\Gamma,H,N} \in V_{\Gamma,H,N}$,
\begin{equation}
  \label{eq:galerkine_orth}
  a(u_\Gamma - u_{\Gamma,H,N}, u_\Gamma - u_{\Gamma,H,N})
  \leq
  a(u_\Gamma - v_{\Gamma,H,N}, u_\Gamma - v_{\Gamma,H,N}).
\end{equation}
Using Lemma~\ref{lem:Babuska}, there exists a function $P(u) \in V^N_{H,0}$ such that, for any $s>3/2$,
\begin{equation}
  \label{eq:bab_approx}
  \|u - P(u) \|_{H^1(\Omega)} \leq C_s \ \frac{H^{\min(s,N+1)-1}}{N^{s-1}} \ \|u\|_{H^s(\Omega)}.
\end{equation}
We consider the harmonic lifting $w = E_\Omega(P(u))$ of the restriction of $P(u)$ on $\Gamma$. The function $w$ is defined on $\Omega$, and is the solution to $-\div (A \nabla w)=0$ on each coarse element $K$ with the Dirichlet boundary conditions $w = P(u)$ on $\Gamma$. Note that $P(u)$ is smooth on each edge, and globally continuous on $\Gamma$. It thus belongs to $H^{1/2}(\Gamma)$, which implies that $w$ is well-defined and belongs to $H^1(\Omega)$. Moreover, on each edge, $P(u)$ is a polynomial function of degree lower than or equal to $N$. We therefore have that $w \in V_{\Gamma,H,N}$. 

We now write
\begin{align}
  & a\Big(u_\Gamma - E_\Omega(P(u)), u_\Gamma - E_\Omega(P(u))\Big)
  \nonumber
  \\
  &=
  \sum_{K \in \mathcal{T}_H} \int_K \left( \nabla u_\Gamma - \nabla E_\Omega(P(u)) \right)^T A \left( \nabla u_\Gamma - \nabla E_\Omega(P(u)) \right)
  \nonumber
  \\
  &\leq 
  \sum_{K \in \mathcal{T}_H}
  \int_K \left( \nabla u - \nabla P(u)) \right)^T A \left( \nabla u - \nabla P(u) \right)
  \nonumber
  \\
  &\leq 
  \alpha_{\max} \sum_{K \in \mathcal{T}_H} | u - P(u) |_{H^1(K)}^2
  \nonumber
  \\
  &=
  \alpha_{\max} | u - P(u) |_{H^1(\Omega)}^2,
  \label{eq:youpi}
\end{align}
where the first inequality above again comes from the fact that $u_\Gamma - E_\Omega(P(u))$ is harmonic in each element $K$ and agrees with $u - P(u)$ on $\partial K$ (similar arguments were used to show the first inequality of~\eqref{eq:jeudi2}). Collecting~\eqref{eq:galerkine_orth} (that we use for $v_{\Gamma,H,N} = E_\Omega(P(u))$, which indeed belongs to $V_{\Gamma,H,N}$, as recalled above), \eqref{eq:youpi} and~\eqref{eq:bab_approx}, we conclude the general proof of Lemma~\ref{lem:Error_Interface}.
\qed

\subsection{Proof of Proposition~\ref{prop:error_indicator}} \label{sec:proofs_trois}

The proof of Proposition~\ref{prop:error_indicator} requires the following three results, namely Lemmas~\ref{lem:regul_elliptic_a_mat}, \ref{lem:regul_elliptic_a_mat_bis} and~\ref{lem:clem_interp}, which are stated for triangular meshes or quadrangular meshes.

\medskip

We first need the following elliptic regularity result (see~\cite[Theorems~2.2.2.3 and~3.2.1.2]{grisvard2011elliptic} and also~\cite[p.~176]{savare98}).

\begin{lemma}[from~\cite{grisvard2011elliptic}] \label{lem:regul_elliptic_a_mat}
Consider the reference element $K_{\rm ref}$ of the mesh $\mathcal{T}_H$, and assume that $K_{\rm ref}$ is convex. Let $g_{\rm ref} \in L^2(K_{\rm ref})$ and let $A_{\rm ref}$ be a symmetric matrix-valued diffusion coefficient satisfying the ellipticity condition~\eqref{eq:A_ellip} in $K_{\rm ref}$ and such that $\dps A_{\rm ref} \in \left( C^1(\overline{K_{\rm ref}}) \right)^{d \times d}$. Consider $z_{\rm ref} \in H^1_0(K_{\rm ref})$ solution to
\begin{equation} \label{eq:inria1_mat}
    - \div ( A_{\rm ref} \, \nabla z_{\rm ref} ) = g_{\rm ref} \quad \text{in $K_{\rm ref}$}.
\end{equation}
Then the function $z_{\rm ref}$ belongs to $H^2(K_{\rm ref})$ and there exists $C_{A_{\rm ref}}$, which only depends on $\alpha_{\rm min}$ and $\| A_{\rm ref} \|_{C^1(\overline{K_{\rm ref}})}$ such that, for any edge $e_{\rm ref} \subset \partial K_{\rm ref}$, we have
  \begin{equation} \label{eq:inria2_mat}
    \| \nabla z_{\rm ref} \|_{H^{1/2}(e_{\rm ref})} \leq C_{A_{\rm ref}} \, \| g_{\rm ref} \|_{L^2(K_{\rm ref})}.
  \end{equation}
\end{lemma}
The bound~\eqref{eq:inria2_mat} is not shown in~\cite{grisvard2011elliptic}, but it is a direct consequence of the following facts. Consider the operator $T$ from $E = H^1_0(K_{\rm ref}) \cap H^2(K_{\rm ref})$ to $F = L^2(K_{\rm ref})$ which, to any element $v \in E$, associates $T v = -\div ( A_{\rm ref} \, \nabla v )$. The operator $T$ is of course linear and continuous from $E$ to $F$ and injective (i.e. one-to-one). It is also surjective (i.e. onto) in view of~\cite[Theorems~2.2.2.3 and~3.2.1.2]{grisvard2011elliptic}. It is thus bijective. As a consequence of the open mapping theorem (see e.g.~\cite[Corollary~2.7]{brezis}), $T^{-1}$ is continuous from $F$ into $E$, which yields, using a trace estimate, the bound~\eqref{eq:inria2_mat}.

\medskip

The above result implies the following one.

\begin{lemma} \label{lem:regul_elliptic_a_mat_bis}
Consider an element $K$ of diameter $H$ in the mesh $\mathcal{T}_H$. We assume that $K$ is convex and that $H \leq 1$. Let $\widehat{K}$ be the image of $K$ by the map $x \mapsto x/H$. Let $\widehat{g} \in L^2(\widehat{K})$ and let $\widehat{A}$ be a symmetric matrix-valued diffusion coefficient satisfying the ellipticity condition~\eqref{eq:A_ellip} in $\widehat{K}$ and such that $\dps \widehat{A} \in \left( C^1\left( \overline{\widehat{K}} \right) \right)^{d \times d}$. Consider $\widehat{z} \in H^1_0(\widehat{K})$ solution to
\begin{equation} \label{eq:inria1_mat_bis}
    - \div ( \widehat{A} \, \nabla \widehat{z} ) = \widehat{g} \quad \text{in $\widehat{K}$}.
\end{equation}
Then the function $\widehat{z}$ belongs to $H^2(\widehat{K})$ and there exists $C_{\widehat{A}}$, which only depends on the regularity of the mesh (in the sense of~\eqref{eq:def_mesh_regular}), $\alpha_{\rm min}$ and $\| \widehat{A} \|_{C^1(\overline{\widehat{K}})}$ (and is thus independent of $H$), such that, for any edge $\widehat{e} \subset \partial \widehat{K}$, we have
  \begin{equation} \label{eq:inria2_mat_bis}
    \| \nabla \widehat{z} \|_{H^{1/2}(\widehat{e})} \leq C_{\widehat{A}} \, \| \widehat{g} \|_{L^2(\widehat{K})}.
  \end{equation}
\end{lemma}

\medskip

\noindent
{\bf Proof of Lemma~\ref{lem:regul_elliptic_a_mat_bis}.}
The proof is performed using Lemma~\ref{lem:regul_elliptic_a_mat}. We recall that, as assumed in~\eqref{eq:def_mesh_regular}, there exists an affine transformation $\widehat{F} : K_{\rm ref} \mapsto \widehat{K}$ such that $\| \nabla \widehat{F} \|_{L^\infty} \leq \gamma$ and $\| \nabla \widehat{F}^{-1} \|_{L^\infty} \leq \gamma$. We define $z_{\rm ref}(x) = \widehat{z}(\widehat{F}(x))$ and $g_{\rm ref}(x) = \widehat{g}(\widehat{F}(x))$ for any $x$ in the reference element $K_{\rm ref}$ of unit diameter. Using that $\nabla \widehat{F}$ is a constant, we compute that
$$
- \div (A_{\rm ref} \, \nabla z_{\rm ref}) = g_{\rm ref} \quad \text{in $K_{\rm ref}$},
$$
where $A_{\rm ref} = (\nabla \widehat{F}^{-1})^T \, \widehat{A} \, \nabla \widehat{F}^{-1}$. Using the above bounds on $\nabla \widehat{F}$ and $\nabla \widehat{F}^{-1}$, we observe that the symmetric matrix $A_{\rm ref}$ is bounded from below by a constant only depending on $\alpha_{\rm min}$ and $\gamma$, and that $\| A_{\rm ref} \|_{C^1(\overline{K_{\rm ref}})} \leq \gamma^2 \| \widehat{A} \|_{C^1(\overline{\widehat{K}})}$. Using Lemma~\ref{lem:regul_elliptic_a_mat}, we obtain that $\dps \| \nabla z_{\rm ref} \|_{H^{1/2}(e_{\rm ref})} \leq C \, \| g_{\rm ref} \|_{L^2(K_{\rm ref})}$, where $C$ only depends on $\gamma$, $\alpha_{\rm min}$ and $\| \widehat{A} \|_{C^1(\overline{\widehat{K}})}$. By a change of variable, we obtain~\eqref{eq:inria2_mat_bis}.
\qed

\bigskip

Our third and last technical lemma is the following approximation result. Consider a mesh $\mathcal{T}_H$ and choose a maximal polynomial degree $p_K \in \N^\star$ for any element $K\in \mathcal{T}_H$. We assume that these degrees are uniformly comparable on neighboring elements, in the sense that
\begin{equation} \label{eq:degre_element_proche}
\forall K, K' \in {\mathcal T}_H \ \text{s.t.} \ \overline{K} \cap \overline{K'} \neq \emptyset, \quad \frac{p_K}{\gamma} \leq p_{K'} \leq \gamma \, p_K,
\end{equation}
where $\gamma$ is the mesh regularity constant of~\eqref{eq:def_mesh_regular}. We then have the following result.

\begin{lemma}[Scott-Zhang type interpolation result, see~Theorem~3.3 of~\cite{melenk2003hp_article}]
\label{lem:clem_interp}
Assume that $\mathcal{T}_H$ is a conformal mesh which is regular in the sense of~\eqref{eq:def_mesh_regular}. For any element $K\in \mathcal{T}_H$, we choose a maximal degree $p_K \in \N^\star$ and we assume that these degrees $\{ p_K \}$ satisfy~\eqref{eq:degre_element_proche}. Then there exists a continuous interpolation operator $\mathcal{SZ}$ from $H^1_0(\Omega)$ to $H^1_0(\Omega) \cap \mathcal{S}(\{ p_K \})$, where
$$
\mathcal{S}(\{ p_K \}) = \{u \in C^0(\overline{\Omega}); \ \text{$u|_K$ is a polynomial function of degree at most $p_K$} \}.
$$
Furthermore, there exists a constant $C$ which only depends on the mesh regularity constant $\gamma$ of~\eqref{eq:def_mesh_regular} such that, for any $u \in H^1_0(\Omega)$ and any edge $e \subset \Gamma$, it holds that
\begin{equation} \label{eq:clem_interp}
  \| u-\mathcal{SZ}(u) \|_{L^2(e)} \leq C \left(\frac{H_e}{p_e}\right)^{1/2} |u|_{H^1(\omega_e)},
\end{equation}
where $\omega_e$ is the union of all the elements that share a vertex with the edge $e$, $H_e$ is the length of the edge $e$ and $p_e = \min \{p_K \ | \ e \subset \partial K \}$.
\end{lemma}

We are now in position to prove Proposition~\ref{prop:error_indicator}.

\medskip

\noindent
{\bf Proof of Proposition~\ref{prop:error_indicator}.}
The proof falls in two steps: we first estimate $u_\Gamma - u_{\Gamma,H,\{N_e\}}$ and next $u_B - u_{B,H,\{M_K\}}$.

\medskip

\noindent \textbf{Step 1: interface approximation.}
For the numerical solution $u_{\Gamma,H,\{N_e\}} \in V_{\Gamma,H,\{N_e\}}$, we write, using an integration by parts over every element $K$ and the definition~\eqref{eq:def_phi_gamma} of the basis functions of $V_{\Gamma,H,\{N_e\}}$, that, for any $w_\Gamma \in V_\Gamma$,
\begin{align}
a(u_{\Gamma,H,\{N_e\}},w_\Gamma)
&=
\int_\Omega (\nabla w_\Gamma)^T A \nabla u_{\Gamma,H,\{N_e\}}
\nonumber
\\
&=
\sum_{K\in\mathcal{T}_H} \sum_{e \subset \partial K} \int_e \left( \nu^T A \nabla u_{\Gamma,H,\{N_e\}} \right) w_\Gamma
\nonumber
\\
&=
\sum_{e \subset \Gamma} \int_e w_\Gamma \, J_e\left(\nu^T A \nabla u_{\Gamma,H,\{N_e\}} \right),
\label{eq:a_posteriori_relation_1_FL}
\end{align}
where, we recall, $J_e(\psi)$ denotes the jump of a given function $\psi$ across the edge $e$ and $\nu$ is a normal vector to the edge.

Using~\eqref{eq:approx_interface}, we write that, for any $v_\Gamma \in V_\Gamma$ and any $v_{\Gamma,H,\{N_e\}} \in V_{\Gamma,H,\{N_e\}}$,
\begin{align}
& a(u_\Gamma-u_{\Gamma,H,\{N_e\}},v_\Gamma)
\nonumber
\\
&=
a(u_\Gamma-u_{\Gamma,H,\{N_e\}},v_\Gamma-v_{\Gamma,H,\{N_e\}})
\nonumber
\\
&=
a(u_\Gamma,v_\Gamma-v_{\Gamma,H,\{N_e\}}) - a(u_{\Gamma,H,\{N_e\}},v_\Gamma-v_{\Gamma,H,\{N_e\}})
\nonumber
\\
&=
\int_\Omega f \, (v_\Gamma - v_{\Gamma,H,\{N_e\}}) - \sum_{e\subset\Gamma} \int_e (v_\Gamma - v_{\Gamma,H,\{N_e\}}) \, J_e \left(\nu^T A \nabla u_{\Gamma,H,\{N_e\}} \right),
\label{eq:jeudi3}
\end{align}
where, in the last line, we have used the definition of the exact solution $u_\Gamma$ and~\eqref{eq:a_posteriori_relation_1_FL} for $w_\Gamma = v_\Gamma-v_{\Gamma,H,\{N_e\}}$.

\medskip

We now make the following specific choices. Since we aim at estimating $u_\Gamma-u_{\Gamma,H,\{N_e\}}$, the natural choice for $v_\Gamma$ (see the left-hand side of~\eqref{eq:jeudi3}) is $v_\Gamma = u_\Gamma-u_{\Gamma,H,\{N_e\}}$, a choice we will make at the very end of the present Step~1 (see just above~\eqref{eq:retour8}). Next, since the difference $v_\Gamma - v_{\Gamma,H,\{N_e\}}$ appears in the right hand-side of~\eqref{eq:jeudi3} and we intend to have this right-hand side as small as possible for our estimator (see the bounds~\eqref{eq:retour6} and~\eqref{eq:retour7} below), we wish to choose $v_{\Gamma,H,\{N_e\}}$ as close as possible to $v_\Gamma$, under the constraint that $v_{\Gamma,H,\{N_e\}}$ should belong to $V_{\Gamma,H,\{N_e\}}$, and thus should be the harmonic extension of some function that is piecewise polynomial on $\Gamma$. To this end, we are going to define $v_{\Gamma,H,\{N_e\}}$ as the extension (see~\eqref{eq:choice_v_gamma} below) of some function $\varphi$ defined on $\Gamma$ and that approximates $v_\Gamma$ in some sense. 

We now proceed in details and define the function $\varphi$ on $\Gamma$ by
\begin{equation} \label{eq:def_ww}
\varphi = \mathcal{SZ}(v_\Gamma)|_\Gamma + \sum_{e \in \Gamma} \Pi_{N_e}^{e,0}(v_\Gamma-\mathcal{SZ}(v_\Gamma)),
\end{equation}
where, for any edge $e \subset \Gamma$, $\Pi_{N_e}^{e,0}$ is the $L^2(e)$ projection on the polynomial functions that vanish at both ends of the edge and of degree lower than or equal to $N_e$ on $e$ (by construction, for any function $\psi$, $\Pi_{N_e}^{e,0}(\psi)$ is supported on the edge $e$). Since $v_\Gamma$ and $\mathcal{SZ}(v_\Gamma)$ (see below) belong to $H^1_0(\Omega)$, they belong to $L^2(e)$ and $\Pi_{N_e}^{e,0}(v_\Gamma-\mathcal{SZ}(v_\Gamma))$ is well-defined. In~\eqref{eq:def_ww}, $\mathcal{SZ}$ is the Scott-Zhang type interpolant defined in Lemma~\ref{lem:clem_interp}, where we choose, for each element $K$, the polynomial degree
\begin{equation} \label{eq:choix_pk_SZ}
  p_K = \min \{N_e \ | \ e \subset \partial K\}.
\end{equation}

We next observe that, on each edge $e$, $\varphi$ is a polynomial function of degree lower than or equal to $N_e$ (a property that will be useful below to ensure that the harmonic extension of $\varphi$ belongs to the right space). This is obviously the case for the second term in~\eqref{eq:def_ww}. This is also the case for the first term, which is indeed a polynomial function of degree $p_{K^1_e}$ (resp. $p_{K^2_e}$) on $K^1_e$ (resp. $K^2_e$), where $K^1_e$ and $K^2_e$ are the two elements sharing the edge $e$. By construction (see~\eqref{eq:choix_pk_SZ}), we have $p_{K^1_e} \leq N_e$ and likewise for $p_{K^2_e}$. 

Since $\varphi$ is globally continuous on $\Gamma$ (because $\mathcal{SZ}(v_\Gamma)$ is continuous on $\overline{\Omega}$ and $\Pi_{N_e}^{e,0}$ is a polynomial that vanishes at the edge boundaries) and smooth on each edge, it belongs to $H^{1/2}(\Gamma)$ and we can consider its harmonic lifting
\begin{equation}
  \label{eq:choice_v_gamma}
  v_{\Gamma,H,\{N_e\}} = E_\Omega(\varphi) = E_\Omega \Big( \mathcal{SZ}(v_\Gamma)|_\Gamma + \sum_{e \in \Gamma} \Pi_{N_e}^{e,0}(v_\Gamma-\mathcal{SZ}(v_\Gamma)) \Big),
\end{equation}
which belongs to $H^1(\Omega)$. Since $\varphi$ is a polynomial function of degree lower than or equal to $N_e$ on any edge $e$, we observe that $v_{\Gamma,H,\{N_e\}}$ belongs to the approximation space $V_{\Gamma,H,\{N_e\}}$. 

\medskip

In passing, we observe that, by construction, the degrees $\{ p_K \}$ satisfy~\eqref{eq:degre_element_proche}. Consider indeed two neighboring elements $K$ and $K'$. Then, denoting $\widetilde{e}$ the edge shared by $K$ and $K'$, we have
$$
\frac{p_K}{p_{K'}}
=
\frac{\min \{N_e \ | \ e \subset \partial K\}}{\min \{N_{e'} \ | \ e' \subset \partial K'\}}
=
\frac{\min \{N_e \ | \ e \subset \partial K\}}{N_{\widetilde{e}}} \ \frac{N_{\widetilde{e}}}{\min \{N_{e'} \ | \ e' \subset \partial K'\}}
\leq
\gamma,
$$
where we have used the property~\eqref{eq:degre_edge_proche}. We likewise have that $p_K/p_{K'} \geq 1/\gamma$. Since the degrees $\{ p_K \}$ satisfy~\eqref{eq:degre_element_proche}, we will be in position to use the approximation result~\eqref{eq:clem_interp} in the sequel.

\medskip

For any $v_\Gamma \in V_\Gamma$, we thus define $v_{\Gamma,H,\{N_e\}} \in V_{\Gamma,H,\{N_e\}}$ by~\eqref{eq:choice_v_gamma}. In the sequel of the proof, we bound $v_\Gamma - v_{\Gamma,H,\{N_e\}}$ in $L^2$ norm in the bulk of each element and on its boundaries. This is the purpose of Step~1a (see~\eqref{eq:parc2} and~\eqref{eq:retour1} below). We next use these bounds to successively majorize in Step~1b the two terms in the right-hand side of~\eqref{eq:jeudi3}, which yields the estimate~\eqref{eq:retour8} below of the interface approximation and concludes Step~1.

\medskip

\noindent \textbf{Step 1a.}
To bound the first term of~\eqref{eq:jeudi3}, we need to estimate $\| v_\Gamma - v_{\Gamma,H,\{N_e\}} \|_{L^2(K)}$ for any element $K \in \mathcal{T}_H$. To this end, we introduce the unique solution $z$ in $H^1_0(K)$ to
\begin{equation} \label{eqn:define_z}
- \div ( A \nabla z ) = v_\Gamma - v_{\Gamma,H,\{N_e\}} \quad \text{in $K$}.
\end{equation}
Since $K$ is convex and $\dps A \in \left( C^1(\overline{\Omega}) \right)^{d \times d}$, we know that $z \in H^2(K)$, by elliptic regularity (see a similar result in Lemma~\ref{lem:regul_elliptic_a_mat_bis}).

Using the definition of $z$, we have
\begin{align*}
  & \| v_\Gamma-v_{\Gamma,H,\{N_e\}} \|_{L^2(K)}^2
  \\
  &=
  \int_K (\nabla (v_\Gamma-v_{\Gamma,H,\{N_e\}}))^T A \nabla z  - \sum_{e \subset \partial K} \int_e \left( v_\Gamma-v_{\Gamma,H,\{N_e\}} \right) \, \nu^T A \nabla z
  \\
  &=
  - \sum_{e \subset \partial K} \int_e \left( v_\Gamma-v_{\Gamma,H,\{N_e\}} \right) \, \nu^T A \nabla z,
\end{align*}
where the first term of the second line vanishes since $A$ is symmetric, both $v_\Gamma$ and $v_{\Gamma,H,\{N_e\}}$ are harmonic and $z$ vanishes on $\partial K$. Since $\dps A \in \left( C^1(\overline{\Omega}) \right)^{d \times d}$ and $z \in H^2(K)$, we have that $\nu^T A \nabla z \in H^{1/2}(e)$.
Setting
\begin{equation} \label{eq:def_w}
w = \nu^T A \nabla z
\end{equation}
and using our specific choice~\eqref{eq:choice_v_gamma}, we get
\begin{align}
  & \| v_\Gamma-v_{\Gamma,H,\{N_e\}} \|_{L^2(K)}^2
  \nonumber
  \\
  &=
  - \sum_{e \subset \partial K} \int_e \Big( v_\Gamma - \mathcal{SZ}(v_\Gamma) - \Pi^{e,0}_{N_e}(v_\Gamma-\mathcal{SZ}(v_\Gamma)) \Big) \, w
  \nonumber
  \\
  &=
  - \sum_{e \subset \partial K} \int_e \Big( v_\Gamma - \mathcal{SZ}(v_\Gamma) - \Pi^{e,0}_{N_e}(v_\Gamma-\mathcal{SZ}(v_\Gamma)) \Big) \, \left( w - \Pi^{e,0}_{N_e}(w) \right)
  \nonumber
  \\
  &\leq
  \sum_{e \subset \partial K} \left\| v_\Gamma - \mathcal{SZ}(v_\Gamma) - \Pi^{e,0}_{N_e}(v_\Gamma-\mathcal{SZ}(v_\Gamma)) \right\|_{L^2(e)} \, \left\| w - \Pi^{e,0}_{N_e}(w) \right\|_{L^2(e)},
  \label{eq:inria6}
\end{align}
where the second equality stems from the fact that $\Pi^{e,0}_{N_e}$ is a $L^2(e)$ orthogonal projection. We successively bound the two factors of~\eqref{eq:inria6}.

\medskip

Successively using the definition of the projection and Lemma~\ref{lem:clem_interp}, we obtain, for the first factor of~\eqref{eq:inria6},
\begin{align} 
  \left\| v_\Gamma - \mathcal{SZ}(v_\Gamma) - \Pi^{e,0}_{N_e}(v_\Gamma-\mathcal{SZ}(v_\Gamma)) \right\|_{L^2(e)}
  &\leq
  \left\| v_\Gamma - \mathcal{SZ}(v_\Gamma) \right\|_{L^2(e)}
  \nonumber
  \\
  &\leq
  C \sqrt{\frac{H_e}{p_e}} \, |v_\Gamma|_{H^1(\omega_e)},
  \label{eq:retour1}
\end{align}
with
\begin{equation} \label{eq:def_pe}
p_e
=
\min\{ p_{K^1_e}, p_{K^2_e} \}
=
\min\{ N_{\widetilde{e}} \ | \ \widetilde{e} \subset \partial K^1_e \cup \partial K^2_e \},
\end{equation}
where $K^1_e$ and $K^2_e$ are the two elements sharing the edge $e$ and $p_{K^1_e}$ and $p_{K^2_e}$ are the degrees chosen in~\eqref{eq:choix_pk_SZ} for the construction of the Scott-Zhang type interpolation operator $\mathcal{SZ}$. The explicit expression~\eqref{eq:def_pe} yields the value of $p_e$ cited in the statement of Proposition~\ref{prop:error_indicator}.

\medskip

We now turn to the second factor of~\eqref{eq:inria6}. Introduce the image $\widehat{K}$ of $K$ by the map $x \mapsto x/H_K$, and define $\widehat{w}(x) = w(H_K \, x)$ on $\widehat{e}$. We then have
\begin{equation}
\left\| w - \Pi^{e,0}_{N_e}(w) \right\|_{L^2(e)} = \sqrt{H_K} \, \left\| \widehat{w} - \Pi^{\widehat{e},0}_{N_e}(\widehat{w}) \right\|_{L^2(\widehat{e})}.
\label{eq:jeudi5}
\end{equation}
We are going to bound the right-hand side of~\eqref{eq:jeudi5} by interpolation. By definition of the projection, we have
\begin{equation} \label{eq:retour1_}
\forall \widehat{w} \in L^2(\widehat{e}), \quad \left\| \widehat{w} - \Pi^{\widehat{e},0}_{N_e}(\widehat{w}) \right\|_{L^2(\widehat{e})} \leq \left\| \widehat{w} \right\|_{L^2(\widehat{e})}.
\end{equation}
Second, we have
\begin{equation} \label{eq:retour2}
\forall \widehat{w} \in H^1_0(\widehat{e}), \quad \left\| \widehat{w} - \Pi^{\widehat{e},0}_{N_e}(\widehat{w}) \right\|_{L^2(\widehat{e})} \leq \left\| \widehat{w} - I_{N_e}(\widehat{w}) \right\|_{L^2(\widehat{e})} \leq \frac{C}{N_e} \, |\widehat{w}|_{H^1(\widehat{e})},
\end{equation}
where $I_{N_e}$ is the interpolant of degree $N_e$ at the Gauss Lobatto points of the edge $\widehat{e}$. Note that the first inequality in~\eqref{eq:retour2} critically relies on the fact that $\widehat{w}$ vanishes at the two vertices of the edge (hence $I_{N_e}(\widehat{w})$ also vanishes at the two vertices, and thus can be compared with $\Pi^{\widehat{e},0}_{N_e}(\widehat{w})$). The second inequality in~\eqref{eq:retour2} is for instance given in~\cite[Eq.~(5.4.33)]{canuto2010spectral} and in~\cite[Corollaire~IV.1.13]{bernardi-maday-rapetti}. 

By Sobolev interpolation between $L^2(\widehat{e})$ and $H^1_0(\widehat{e})$ (see Appendix~\ref{sec:H_un_demi}), we deduce from~\eqref{eq:retour1_} and~\eqref{eq:retour2} that, for any $\eta > 0$, there exists $C_\eta$ such that 
\begin{equation} \label{eq:retour3}
\forall \widehat{w} \in H^{1/2-\eta}(\widehat{e}), \quad \left\| \widehat{w} - \Pi^{\widehat{e},0}_{N_e}(\widehat{w}) \right\|_{L^2(\widehat{e})} \leq \frac{C_\eta}{N_e^{1/2-\eta}} \, \| \widehat{w} \|_{H^{1/2-\eta}(\widehat{e})}.
\end{equation}
We thus deduce from~\eqref{eq:jeudi5} and~\eqref{eq:retour3} that, for our function $w$ of interest given by~\eqref{eq:def_w},
\begin{equation} \label{eq:H-1on2}
  \left\| w - \Pi^{e,0}_{N_e}(w) \right\|_{L^2(e)} \leq C_\eta \, \frac{\sqrt{H_K}}{N_e^{1/2-\eta}} \, \| \widehat{w} \|_{H^{1/2}(\widehat{e})},
\end{equation}
where we recall that $\widehat{w}(x) = w(H_K \, x)$ on $\widehat{e}$.

Recalling that $w = \nu^T A \nabla z$ where $z \in H^1_0(K)$ satisfies~\eqref{eqn:define_z}, which we recast as
$$
- \div ( A \nabla z ) = g \quad \text{in $K$}, \qquad g = v_\Gamma - v_{\Gamma,H,\{N_e\}},
$$
we introduce $\widehat{z}$, $\widehat{A}$ and $\widehat{g}$ defined on $\widehat{K}$ by $\widehat{z}(x) = z(H_K \, x)$, $\widehat{A}(x) = A(H_K \, x)$ and $\widehat{g}(x) = g(H_K \, x)$, and compute that
\begin{equation} \label{eqn:define_zz}
- \div ( \widehat{A} \nabla \widehat{z} ) = H_K^2 \, \widehat{g} \quad \text{in $\widehat{K}$}.
\end{equation}
Furthermore, $\widehat{w}(x) = w(H_K \, x) = \nu^T A(H_K \, x) (\nabla z)(H_K \, x) = H_K^{-1} \, \nu^T \widehat{A}(x) (\nabla \widehat{z})(x)$ on $\widehat{e}$. We hence write that
\begin{multline*}
\| \widehat{w} \|_{H^{1/2}(\widehat{e})} \leq H_K^{-1} \| \widehat{A} \nabla \widehat{z} \|_{H^{1/2}(\widehat{e})} \leq H_K^{-1} \| \widehat{A} \|_{C^1(\widehat{e})} \, \| \nabla \widehat{z} \|_{H^{1/2}(\widehat{e})} \\ \leq H_K^{-1} \| A \|_{C^1(\overline{\Omega})} \, \| \nabla \widehat{z} \|_{H^{1/2}(\widehat{e})}.
\end{multline*}
Using next Lemma~\ref{lem:regul_elliptic_a_mat_bis} on~\eqref{eqn:define_zz}, we deduce that
$$
\| \widehat{w} \|_{H^{1/2}(\widehat{e})} \leq H_K^{-1} \| A \|_{C^1(\overline{\Omega})} \, C_{\widehat{A}} \, H_K^2 \, \| \widehat{g} \|_{L^2(\widehat{K})},
$$
where $C_{\widehat{A}}$ only depends on $\gamma$, $\alpha_{\rm min}$ and $\| \widehat{A} \|_{C^1(\overline{\widehat{K}})}$. We thus obtain that
$$
\| \widehat{w} \|_{H^{1/2}(\widehat{e})} \leq C(\gamma,\alpha_{\rm min},\| A \|_{C^1(\overline{\Omega})}) \, H_K \, \| \widehat{g} \|_{L^2(\widehat{K})} = C(\gamma,\alpha_{\rm min},\| A \|_{C^1(\overline{\Omega})}) \, \| g \|_{L^2(K)}. 
$$
Collecting this bound with~\eqref{eq:H-1on2} and recalling that $g = v_\Gamma - v_{\Gamma,H,\{N_e\}}$, we obtain
\begin{equation} \label{eq:parc}
\left\| w - \Pi^{e,0}_{N_e}(w) \right\|_{L^2(e)} \leq C \, \frac{\sqrt{H_K}}{N_e^{1/2-\eta}} \, \| v_\Gamma - v_{\Gamma,H,\{N_e\}} \|_{L^2(K)},
\end{equation}
where $C$ only depends on $\gamma$, $\eta$, $\alpha_{\rm min}$ and $\| A \|_{C^1(\overline{\Omega})}$. Collecting~\eqref{eq:inria6}, \eqref{eq:retour1} and~\eqref{eq:parc}, we deduce that
\begin{equation} \label{eq:parc2}
\| v_\Gamma-v_{\Gamma,H,\{N_e\}} \|_{L^2(K)}
\leq
C \, \sum_{e \subset \partial K} \sqrt{\frac{H_e \, H_K}{p_e}} \, \frac{1}{N_e^{1/2-\eta}} \, |v_\Gamma|_{H^1(\omega_e)}.
\end{equation}

\medskip

\noindent \textbf{Step 1b.}
Using the above bound~\eqref{eq:parc2}, we are now in position to bound the first term of~\eqref{eq:jeudi3} by
\begin{align}
  & \left| \int_\Omega f (v_\Gamma-v_{\Gamma,H,\{N_e\}}) \right|
  \nonumber
  \\
  &\leq
  \sum_{K\subset \mathcal{T}_H} \|f\|_{L^2(K)} \|v_\Gamma-v_{\Gamma,H,\{N_e\}}\|_{L^2(K)}
  \nonumber
  \\
  &\leq
  C \sum_{K\subset \mathcal{T}_H} \|f\|_{L^2(K)} \left(\sum_{e\subset \partial K} \sqrt{\frac{H_e \, H_K}{p_e}} \, \frac{1}{N_e^{1/2-\eta}} \, |v_\Gamma|_{H^1(\omega_e)} \right)
  \nonumber
  \\
  &\leq
  C \sqrt{\sum_{K\subset \mathcal{T}_H} \|f\|_{L^2(K)}^2 \left( \sum_{e\subset \partial K} \sqrt{\frac{H_e \, H_K}{p_e}} \, \frac{1}{N_e^{1/2-\eta}} \right)^2} \sqrt{\sum_{K\subset \mathcal{T}_H} |v_\Gamma|_{H^1(\omega_K)}^2}
  \nonumber
  \\
  &\leq
  C \, |v_\Gamma |_{H^1(\Omega)} \, \sqrt{\sum_{K\subset \mathcal{T}_H} \|f\|_{L^2(K)}^2 \left( \sum_{e\subset \partial K} \frac{H_e \, H_K}{N_e^{1-2\eta} \, p_e} \right)},
  \label{eq:retour6}
\end{align} 
where $\dps \omega_K = \cup_{e\subset \partial K} \omega_e$, and where $C$ only depends on $\gamma$, $\eta$, $\alpha_{\rm min}$ and $\| A \|_{C^1(\overline{\Omega})}$.

\medskip

We now consider the second term of~\eqref{eq:jeudi3}. We have shown in~\eqref{eq:retour1} that
$$
\| v_\Gamma - v_{\Gamma,H,\{N_e\}} \|_{L^2(e)} \leq C \sqrt{\frac{H_e}{p_e}} \, |v_\Gamma|_{H^1(\omega_e)},
$$
with $p_e$ given by~\eqref{eq:def_pe}. We therefore have
\begin{align*}
  & \left| \sum_{e\subset \Gamma} \int_e (v_\Gamma-v_{\Gamma,H,\{N_e\}}) J_e\left( \nu^T A \nabla u_{\Gamma,H,\{N_e\}} \right) \right|
  \\
  & \leq \sum_{e\subset \Gamma} \left\| J_e\left(\nu^T A \nabla u_{\Gamma,H,N} \right) \right\|_{L^2(e)} \, \|v_\Gamma-v_{\Gamma,H,\{N_e\}}\|_{L^2(e)}
  \\
  & \leq C \sum_{e\subset \Gamma} \sqrt{\frac{H_e}{p_e}} \ \left\| J_e\left(\nu^T A \nabla u_{\Gamma,H,\{N_e\}} \right) \right\|_{L^2(e)} |v_\Gamma|_{H^1(\omega_e)}.
\end{align*}
Using the discrete Cauchy-Schwarz inequality, we deduce that
\begin{align}
  & \left| \sum_{e\subset \Gamma} \int_e (v_\Gamma-v_{\Gamma,H,\{N_e\}}) J_e\left( \nu^T A \nabla u_{\Gamma,H,\{N_e\}} \right) \right|
  \nonumber
  \\
  & \leq C \sqrt{\sum_{e\subset \Gamma} \frac{H_e}{p_e} \ \left\| J_e\left(\nu^T A \nabla u_{\Gamma,H,\{N_e\}} \right) \right\|_{L^2(e)}^2} \, \sqrt{ \sum_{e\subset \Gamma} |v_\Gamma|_{H^1(\omega_e)}^2}
  \nonumber
  \\
  & \leq C |v_\Gamma|_{H^1(\Omega)} \sqrt{\sum_{e\subset \Gamma} \frac{H_e}{p_e} \ \left\|J_e\left(\nu^T A \nabla u_{\Gamma,H,\{N_e\}} \right) \right\|_{L^2(e)}^2}.
  \label{eq:retour7}
\end{align}
Collecting~\eqref{eq:jeudi3}, \eqref{eq:retour6} and~\eqref{eq:retour7}, we obtain, for any $v_\Gamma \in V_\Gamma$, that
\begin{multline*}
  \frac{a(u_\Gamma-u_{\Gamma,H,\{N_e\}},v_\Gamma)}{|v_\Gamma|_{H^1(\Omega)}}
  \leq C \left\{ \sqrt{\sum_{K\subset \mathcal{T}_H} \|f\|_{L^2(K)}^2 \left( \sum_{e\subset \partial K} \frac{H_e \, H_K}{N_e^{1-2\eta} \, p_e} \right)} \right.
  \\
  \left. + \sqrt{\sum_{e\subset \Gamma} \frac{H_e}{p_e} \ \left\| J_e\left(\nu^T A \nabla u_{\Gamma,H,\{N_e\}} \right) \right\|_{L^2(e)}^2} \right\}.
\end{multline*}
We use the above estimate for the choice $v_\Gamma = u_\Gamma-u_{\Gamma,H,\{N_e\}}$, which obviously belongs to $V_\Gamma$.
We thus deduce that 
\begin{multline} \label{eq:retour8}
\| u_\Gamma - u_{\Gamma,H,\{N_e\}} \|_E \leq C \left\{ \sum_{K\subset \mathcal{T}_H} \|f\|_{L^2(K)}^2 \left( \sum_{e\subset \partial K} \frac{H_e \, H_K}{N_e^{1-2\eta} \, p_e} \right) \right. \\
\left. + \sum_{e\subset \Gamma} \frac{H_e}{p_e} \, \left\| J_e\left(\nu^T A \nabla u_{\Gamma,H,\{N_e\}} \right) \right\|_{L^2(e)}^2 \right\}^{1/2}. 
\end{multline}

\medskip

\noindent \textbf{Step 2: bubble approximation.}
In what follows, we establish an {\em a posteriori} estimate on $u_B - u_{B,H,\{M_K\}}$ in the case when bubble enrichments are considered, i.e. in the case when $M_K \geq 1$ for any element $K$. If no enrichements are used (that is in the case when we approximate $u_B \in V_B$ by $u_{B,H,M} = 0$), then we simply use the right-hand side of~\eqref{eq:error_u_bubble0}, or more precisely $\dps \left( \sum_{K\in\mathcal{T}_H} H_K^2 \, \| f \|_{L^2(K)}^2 \right)^{1/2}$, as \emph{a posteriori} estimator.

For the numerical solution $u_{B,H,\{M_K\}} \in V_{B,H,\{M_K\}}$, we write, using an integration by parts over each element $K$, that, for any $w_B \in V_B$,
\begin{align}
a(u_{B,H,\{M_K\}},w_B)
&=
\int_\Omega \left(\nabla w_B\right)^T A \nabla u_{B,H,\{M_K\}}
\nonumber
\\
&=
- \sum_{K\in\mathcal{T}_H} \int_K w_B \, \div \left( A \nabla u_{B,H,\{M_K\}} \right). 
\label{eq:samedi1}
\end{align}
Successively using~\eqref{eq:approx_bubble} and~\eqref{eq:samedi1} for $w_B = v_B - v_{B,H,\{M_K\}}$, we deduce that, for any $v_B \in V_B$ and any $v_{B,H,\{M_K\}} \in V_{B,H,\{M_K\}}$,
\begin{align}
& a(u_B-u_{B,H,\{M_K\}},v_B)
\nonumber
\\
&=
a(u_B-u_{B,H,\{M_K\}},v_B-v_{B,H,\{M_K\}})
\nonumber
\\
&=
a(u_B,v_B-v_{B,H,\{M_K\}}) - a(u_{B,H,\{M_K\}},v_B-v_{B,H,\{M_K\}})
\nonumber
\\
&=
\int_\Omega f \, (v_B - v_{B,H,\{M_K\}}) + \sum_{K\in\mathcal{T}_H} \int_K (v_B - v_{B,H,\{M_K\}}) \, \div \left( A \nabla u_{B,H,\{M_K\}} \right)
\nonumber
\\
&=
\sum_{K\in\mathcal{T}_H} \int_K (v_B - v_{B,H,\{M_K\}}) \, \big( f + \div \left( A \nabla u_{B,H,\{M_K\}} \right) \big).
\label{eq:samedi2}
\end{align}
We now make the following specific choices in~\eqref{eq:samedi2}. Since we aim at estimating $u_B-u_{B,H,\{M_K\}}$, the natural choice for $v_B$ (see the left-hand side of~\eqref{eq:samedi2}) is
\begin{equation}
\label{eq:choice2}
v_B = u_B-u_{B,H,\{M_K\}}.
\end{equation}
Since the difference $v_B - v_{B,H,\{M_K\}}$ appears in the right hand-side of~\eqref{eq:samedi2} and we aim at having this right-hand side as small as possible for our estimator, this choice of $v_B$ in turn suggests to define $v_{B,H,\{M_K\}}$ as the solution on each element $K$ to
\begin{equation}
\label{eq:choice3}
- \div (A \nabla v_{B,H,\{M_K\}}) = \Pi^K_{M_K}(z) \ \ \text{in $K$}, \qquad v_{B,H,\{M_K\}} = 0 \ \ \text{on $\partial K$},
\end{equation}
where
\begin{equation}
\label{eq:choice3_aa}
z = - \div (A \nabla v_B) \ \ \text{in $K$},
\end{equation}
and where $\Pi^K_{M_K}$ is the $L^2(K)$-projection on the polynomials of degree at most $M_K$ on the element $K$.

In view of~\eqref{eq:choice2} and of the definition of $u_B$, we see that, in $K$,
\begin{equation} \label{eq:def_zz}
z = - \div (A \nabla u_B) + \div (A \nabla u_{B,H,\{M_K\}}) = f + \div (A \nabla u_{B,H,\{M_K\}}).
\end{equation}
In view of~\eqref{eq:def_phi_B}, we thus see that the right-hand side $\Pi^K_{M_K}(z)$ in~\eqref{eq:choice3} satisfies
\begin{equation} \label{eq:def_Piz}
\Pi^K_{M_K}(z) = \Pi^K_{M_K}(f) + \div (A \nabla u_{B,H,\{M_K\}}).
\end{equation}
In the sequel of the present Step~2, we successively bound $v_B - v_{B,H,\{M_K\}}$ in Step~2a, in order to next estimate the right-hand side of~\eqref{eq:samedi2} in Step~2b.

\medskip

\noindent \textbf{Step 2a.}
Subtracting~\eqref{eq:choice3} to~\eqref{eq:choice3_aa}, we have
$$
- \div ( A \nabla (v_B - v_{B,H,\{M_K\}}) ) = z - \Pi^K_{M_K}(z) \ \ \text{in $K$}, \quad v_B - v_{B,H,\{M_K\}} = 0 \ \ \text{on $\partial K$}.
$$
Successively using $v_B - v_{B,H,\{M_K\}}$ as test function in the variational formulation of that equation, the Poincar\'e inequality on the mesh element $K$ and the coercivity of the problem, we obtain
\begin{equation}
\| v_B - v_{B,H,\{M_K\}} \|_{L^2(K)}
\leq
C \, \frac{H_K^2}{\alpha_{\rm min}} \, \left\| z - \Pi^K_{M_K}(z) \right\|_{L^2(K)}
=
C \, \frac{H_K^2}{\alpha_{\rm min}} \, \left\| f - \Pi^K_{M_K}(f) \right\|_{L^2(K)},
\label{eq:samedi3}
\end{equation}
where the last equality is obtained by subtracting~\eqref{eq:def_Piz} to~\eqref{eq:def_zz}.

\medskip

\noindent \textbf{Step 2b.}
Inserting~\eqref{eq:samedi3} in~\eqref{eq:samedi2}, and using there that $v_B$ is given by~\eqref{eq:choice2}, we get that
\begin{align}
& a(u_B-u_{B,H,\{M_K\}},u_B-u_{B,H,\{M_K\}})
\nonumber
\\
&\leq
\sum_{K\in\mathcal{T}_H} \| v_B - v_{B,H,\{M_K\}} \|_{L^2(K)} \, \| f + \div \left( A \nabla u_{B,H,\{M_K\}} \right) \|_{L^2(K)}
\nonumber
\\
& \leq
\frac{C}{\alpha_{\rm min}} \sum_{K\in\mathcal{T}_H} H_K^2 \, \left\| f - \Pi^K_{M_K}(f) \right\|_{L^2(K)} \, \| f + \div \left( A \nabla u_{B,H,M} \right) \|_{L^2(K)}.
\label{eq:samedi5}
\end{align}
Since $f \in H^{\ell_K}(K)$ for some integer $\ell_K \geq 0$ satisfying $\ell_K \leq \overline{\ell}$, we know from Lemma~\ref{lem:approx_poly} that
\begin{equation} \label{eq:estim_f}
\left\| f - \Pi_{M_K}^K(f) \right\|_{L^2(K)} \leq C(\overline{\ell}) \, \frac{H_K^{\min(\ell_K,M_K+1)}}{M_K^{\ell_K}} \, \| f \|_{H^{\ell_K}(K)}
\end{equation}
with $\dps C(\overline{\ell}) = \max_{0 \leq \ell \leq \overline{\ell}} C_\ell$, where $C_\ell$ is the constant in~\eqref{eq:approx_poly}. Inserting~\eqref{eq:estim_f} in~\eqref{eq:samedi5} yields
\begin{multline} \label{eq:retour9}
\| u_B - u_{B,H,\{M_K\}} \|_E
\leq \\
\frac{C}{\sqrt{\alpha_{\rm min}}} \left\{ \sum_{K\in\mathcal{T}_H} H_K^2 \, \frac{H_K^{\min(\ell_K,M_K+1)}}{M_K^{\ell_K}} \, \| f + \div \left( A \nabla u_{B,H,\{M_K\}} \right) \|_{L^2(K)} \, \| f \|_{H^{\ell_K}(K)} \right\}^{1/2},
\end{multline}
where the right-hand side is explicitly computable, apart from the unknown constant $C$ (which only depends on $\overline{\ell}$).

\medskip

\noindent \textbf{Step 3.} Collecting~\eqref{eq:retour8} and~\eqref{eq:retour9}, and using the decomposition~\eqref{eq:decoupling} of the error, we obtain~\eqref{eq:a_posteriori_estimator}. This concludes the proof of Proposition~\ref{prop:error_indicator}.
\qed

\section{The Sobolev spaces $H^{1/2}$ and related interpolation results} \label{sec:H_un_demi}

For the convenience of the reader, we collect in this short appendix some standard results on the Sobolev spaces $H^{1/2}$ and on related interpolation results that we need for the proof of the {\em a posteriori} estimate~\eqref{eq:a_posteriori_estimator}. Such results are classical and we refer e.g. to~\cite[Section~1]{grisvard2011elliptic}, \cite[Chapter~1]{magenes1968problemes}, \cite[Section~3]{mclean2000strongly}, \cite[Section~2.3 through 2.8]{SauterSchwab2004} and~\cite[Chapter~33]{book_tartar}. 

First we recall the definition, see e.g.~\cite[Definition B.30]{Ern_FE_book}, of the $H^s$ space for $0 < s < 1$. 

\begin{definition}
\label{def:h_demi}
Let $0 < s < 1$. For any open domain $\omega \subset \RR^n$ and any $u \in L^2(\omega)$, we define the norm
$$
\| u \|^2_{H^s(\omega)} = \| u \|^2_{L^2(\omega)} + | u |^2_{H^s(\omega)},
$$
where
$$
| u |^2_{H^s(\omega)} = \int_\omega \int_\omega \frac{|u(x)-u(y)|^2}{|x-y|^{2s+n}} \, dx dy,
$$
and define the space
$$
H^s(\omega) = \left\{ u \in L^2(\omega), \quad \| u \|_{H^s(\omega)} < \infty \right\}.
$$
\end{definition}
In dimension $n=1$, and if $s>1/2$, functions in $H^s(\omega)$ are continuous on $\overline{\omega}$, and the injection $H^s(\omega) \subset C^0(\overline{\omega})$ is continuous (see~\cite[Chapter~1, Theorem~9.8]{magenes1968problemes}). Moreover, still in that case ($n=1$ and $s>1/2$), and assuming that $\omega = \omega_1 \cup \omega_2 \cup \{z\}$, where $z$ is the intersection point of two segments $\omega_1$ and $\omega_2$, we have
$$
u \in H^s(\omega) \quad \text{if and only if} \quad \left\{ \begin{array}{c}
    u|_{\omega_1} \in H^s(\omega_1), \ \ u|_{\omega_2} \in H^s(\omega_2) \\
    \text{and $u$ is continuous at $z$.}
    \end{array} \right.
$$
Furthermore, the norm $\|\cdot\|_{H^s(\omega)}$ is equivalent to $\|\cdot\|_{H^s(\omega_1)}+\|\cdot\|_{H^s(\omega_2)}$.

\medskip

In dimension $n=1$ and for $0 \leq s <1/2$, we have
$$
u \in H^s(\omega) \quad \text{if and only if} \quad u|_{\omega_1} \in H^s(\omega_1) \text{ and } u|_{\omega_2} \in H^s(\omega_2).
$$
Furthermore, the norm $\|\cdot\|_{H^s(\omega)}$ is equivalent to $\|\cdot\|_{H^s(\omega_1)}+\|\cdot\|_{H^s(\omega_2)}$.

\medskip

The critical case $s=1/2$ deserves more attention. Let $\eps > 0$, and denote $\sigma_1$ (resp. $\sigma_2$) the unit vector parallel to $\omega_1$ (resp. $\omega_2$) respectively pointing toward $z$. It holds that
\begin{equation} \label{eq:accoler}
u \in H^s(\omega) \quad \text{if and only if} \quad \left\{ \begin{array}{c}
    u|_{\omega_1} \in H^s(\omega_1), \quad u|_{\omega_2} \in H^s(\omega_2), \quad \text{and} \\
\dps \int_0^\eps \frac{|u(z-t\sigma_1)-u(z+t\sigma_2)|^2}{t} \ dt < \infty.
\end{array} \right.
\end{equation}
This definition is independent of $\eps$ (of course provided $\eps$ is small enough so that $z-t\sigma_1 \in \omega_1$ and $z+t\sigma_2 \in \omega_2$ for any $t \in (0,\eps)$). We emphasize that the norms $\| \cdot \|_{H^{1/2}(\omega)}$ and $\| \cdot \|_{H^{1/2}(\omega_1)}+\| \cdot \|_{H^{1/2}(\omega_2)}$ are {\em not equivalent}.

\medskip

We next recall (see~\cite[Chapter~1, Theorem~11.1]{magenes1968problemes}) that the set $C^\infty_c(\omega)$ of smooth functions with compact support in $\omega$ is dense in $H^s(\omega)$ if and only if $s \leq 1/2$. 

\medskip

We now turn to the space $H^{1/2}_{00}(\omega)$ (sometimes called the Lions-Magenes space, see~\cite[Chapter~33]{book_tartar}), which is formally the space of functions in $H^{1/2}(\omega)$ which can be extended by zero and remain of regularity $H^{1/2}$. Restricting our presentation to the case $\omega = (0,1)$, and following~\cite[Chapter~1, Theorem~11.7]{magenes1968problemes}, we introduce
$$
H^{1/2}_{00}(\omega) = \left\{ u \in H^{1/2}(\omega), \quad \frac{u}{\sqrt{x(1-x)}} \in L^2(\omega) \right\},
$$
with the norm
$$
\| u \|_{H^{1/2}_{00}(\omega)} = \left( \| u \|^2_{H^{1/2}(\omega)} + \left\| \frac{u}{\sqrt{x(1-x)}} \right\|^2_{L^2(\omega)} \right)^{1/2}.
$$
Note that the function $x \in (0,1) \to x(1-x)$ is positive on $\omega$ and vanishes at the end-points of $\omega$ with a non-trivial derivative. The space $H^{1/2}_{00}(\omega)$ is strictly contained in $H^{1/2}(\omega)$. This definition is consistent with~\eqref{eq:accoler} in the sense that, if $u \in H^{1/2}_{00}(\omega)$, then the extension $\overline{u}$ of $u$, defined say on $(-1,2)$ by $\overline{u} = u$ on $\omega=(0,1)$ and $\overline{u} = 0$ elsewhere, indeed belongs to $H^{1/2}(-1,2)$. Conversely, if $\overline{u} \in H^{1/2}(-1,2)$, then $u$ belongs to $H^{1/2}_{00}(\omega)$.

\medskip

We now turn to interpolation properties.

\begin{lemma}[see~Theorem~5.1 of Chapter~1 of~\cite{magenes1968problemes}]
\label{lem:interp_sobolev}
Let $(\mathcal{X},\mathcal{Y})$ be a couple of separable Hilbert spaces with $\mathcal{X} \subset \mathcal{Y}$, such that $\mathcal{X}$ is dense in $\mathcal{Y}$ and such that the injection from $\mathcal{X}$ to $\mathcal{Y}$ is continuous. Let $(X,Y)$ be another couple of Hilbert spaces with analogous properties. Denote by $\mathcal{L}(X,\mathcal{X})$ the set of linear continuous operators from $X$ to $\mathcal{X}$, and likewise for $\mathcal{L}(Y,\mathcal{Y})$. Let $\pi$ be an operator satisfying $\pi \in \mathcal{L}(X,\mathcal{X}) \cap \mathcal{L}(Y,\mathcal{Y})$. Then, for all $0<\theta<1$, we have
$$
\pi \in \mathcal{L}([X,Y]_\theta,[\mathcal{X},\mathcal{Y}]_\theta),
$$
where the interpolated space $[X,Y]_\theta$ is defined in~\cite[Chapter~1, Definition~2.1]{magenes1968problemes}.
\end{lemma}

It is then well-known that, for any open domain $\omega \subset \RR^n$ and any $0<s<1$, one can define $H^s(\omega)$ by interpolation as $[H^1(\omega),L^2(\omega)]_s$.

Furthermore, for any $0<s<1$ with $s \neq 1/2$, we have $[H^1_0(\omega),L^2(\omega)]_s = H^s_0(\omega)$ (see~\cite[Chapter~1, Theorem~11.6]{magenes1968problemes}), where $H^s_0(\omega)$ is the closure of $C^\infty_c(\omega)$ for the $H^s$-norm (we recall, as pointed out above, that $H^s_0(\omega) = H^s(\omega)$ for any $0 < s \leq 1/2$ and that $H^s_0(\omega)$ is a strict subset of $H^s(\omega)$ for any $1/2 < s < 1$).

The case $s=1/2$ is again critical. For this value, we have $[H^1_0(\omega),L^2(\omega)]_{1/2} = H^{1/2}_{00}(\omega)$ (see~\cite[Chapter~1, Theorem~11.7]{magenes1968problemes}), which is a strict subset of $H^{1/2}_0(\omega) = H^{1/2}(\omega)$.

\medskip

\noindent {\bf Acknowledgments.} The work of CLB, FL and PLR is partly supported by ONR and EOARD. CLB and FL acknowledge the continuous support from these two agencies, in particular under the current Grants ONR N00014-20-1-2691 and EOARD FA8655-20-1-7043. Part of this work has been completed while PLR was visiting the University of Washington in Seattle. The hospitality of that institution and the support of a ``Bourse de Mobilité'' of the Ecole Doctorale SIE at Université Paris-Est are gratefully acknowledged. Some preliminary material, eventually included herein, was originally presented in the plenary address of CLB at DD25, Saint John's, Newfoundland, July 2018. CLB wishes to thank the scientific program committee for their invitation. The authors thank A.~Lozinski for stimulating and enlightening discussions on the work reported here, and for carefully reading a preliminary version of this manuscript. The authors finally thank L.~Chamoin and M.~Gander for enlightening discussions on {\em a posteriori} error estimators and domain decomposition methods, respectively.

\bibliographystyle{plain}
\bibliography{biblio_msfem_legendre}

\begin{thebibliography}{10}

\bibitem{babuvska1987hp}
I.~Babuska and M.~Suri.
\newblock The $hp$ version of the finite element method with quasiuniform
  meshes.
\newblock {\em ESAIM: Mathematical Modelling and Numerical Analysis},
  21(2):199--238, 1987.

\bibitem{Bennighof:2004:An-aa}
J.K. Bennighof and R.B. Lehoucq.
\newblock An automated multilevel substructuring method for eigenspace
  computation in linear elastodynamics.
\newblock {\em SIAM J. Sci. Comput.}, 25(6):2084--2106, 2004.

\bibitem{bernardi-maday-rapetti}
C.~Bernardi, Y.~Maday, and F.~Rapetti.
\newblock {\em Discr\'etisations variationnelles de probl\`emes aux limites
  elliptiques}, volume~45 of {\em Math\'ematiques et Applications}.
\newblock Springer, 2004.

\bibitem{bour:92}
F.~Bourquin.
\newblock Component mode synthesis and eigenvalues of second order operators:
  Discretization and algorithm.
\newblock {\em ESAIM: Mathematical Modelling and Numerical Analysis},
  26:385--423, 1992.

\bibitem{BS}
S.C. Brenner and L.R. Scott.
\newblock {\em The mathematical theory of {F}inite {E}lement methods},
  volume~15.
\newblock Springer, 2008.

\bibitem{brezis}
H.~Brezis.
\newblock {\em Functional Analysis, Sobolev Spaces and Partial Differential
  Equations}.
\newblock Springer, 2011.

\bibitem{canuto2010spectral}
C.~Canuto, M.Y. Hussaini, A.~Quarteroni, and T.A. Zang.
\newblock {\em Spectral methods}.
\newblock Springer, 2006.

\bibitem{CAR00}
C.~Carstensen and S.A. Funken.
\newblock Constants in {C}l\'ement-interpolation error and residual based a
  posteriori error estimates in finite element methods.
\newblock {\em East-West Journal of Numerical Mathematics}, 8(3):153--175,
  2000.

\bibitem{Hou_book_msfem_2009}
Y.~Efendiev and T.~Hou.
\newblock {\em {M}ultiscale {F}inite {E}lement {M}ethods: {T}heory and
  {A}pplications}, volume~4 of {\em Surveys and Tutorials in the Applied
  Mathematical Sciences}.
\newblock Springer New York, first edition, 2009.

\bibitem{Ern_FE_book}
A.~Ern and J.-L. Guermond.
\newblock {\em Theory and practice of finite elements}, volume 159 of {\em
  Applied Mathematical Sciences}.
\newblock Springer-Verlag, New York, 2004.

\bibitem{omnes}
Q.~Feng.
\newblock {\em Development of a multiscale finite element method for
  incompressible flows in heterogeneous media}.
\newblock PhD thesis, Universit\'e Paris Saclay, 2019.
\newblock Available at {\tt https://tel.archives-ouvertes.fr/tel-02325512}.

\bibitem{gander_loneland}
M.J. Gander and A.~Loneland.
\newblock {SHEM}: An optimal coarse space for {RAS} and its multiscale
  approximation.
\newblock In C.-O. Lee, X.-C. Cai, D.E. Keyes, H.H. Kim, A.~Klawonn, E.-J.
  Park, and O.B. Widlund, editors, {\em Domain Decomposition Methods in Science
  and Engineering}, volume 116 of {\em Lecture Notes in Computational Science
  and Engineering}, pages 281--288. Springer, 2016.

\bibitem{gander2015analysis}
M.J. Gander, A.~Loneland, and T.~Rahman.
\newblock Analysis of a new harmonically enriched multiscale coarse space for
  domain decomposition methods.
\newblock {\em arXiv preprint arXiv:1512.05285}, 2015.

\bibitem{gao2018high}
K.~Gao, S.~Fu, and E.T. Chung.
\newblock A high-order multiscale finite-element method for time-domain
  acoustic-wave modeling.
\newblock {\em J. Comput. Phys.}, 360:120--136, 2018.

\bibitem{gervasio_spectral_1997}
P.~Gervasio, E.~Ovtchinnikov, and A.~Quarteroni.
\newblock The spectral projection decomposition method for elliptic equations
  in two dimensions.
\newblock {\em SIAM J. Numer. Anal.}, 34(4):1616--1639, 1997.

\bibitem{grisvard2011elliptic}
P.~Grisvard.
\newblock {\em Elliptic problems in nonsmooth domains}.
\newblock Pitman, 1985.

\bibitem{freefem++}
F.~Hecht.
\newblock New development in {F}ree{F}em++.
\newblock {\em J. Numer. Math.}, 20(3-4):251--265, 2012.

\bibitem{hetmaniuk2014error}
U.~Hetmaniuk and A.~Klawonn.
\newblock Error estimates for a two-dimensional special finite element method
  based on component mode synthesis.
\newblock {\em Electron. Trans. Numer. Anal}, 41:109--132, 2014.

\bibitem{hetmaniuk2010special}
U.~Hetmaniuk and R.B. Lehoucq.
\newblock A special finite element method based on component mode synthesis.
\newblock {\em ESAIM: Mathematical Modelling and Numerical Analysis},
  44(3):401--420, 2010.

\bibitem{Hou:1997:A-maa}
T.~Hou and X.~Wu.
\newblock A multiscale finite element method for elliptic problems in composite
  materials and porous media.
\newblock {\em J. Comput. Phys.}, 134:169--189, 1997.

\bibitem{magenes1968problemes}
J.-L. Lions and E.~Magenes.
\newblock {\em Probl{\`e}mes aux limites non homog{\`e}nes et applications},
  volume~1.
\newblock Dunod, 1968.
\newblock (English version: Non-homogeneous boundary value problems and
  applications. Volume~I. Translated from the French by P.~Kenneth. Die
  Grundlehren der mathematischen Wissenschaften, Band 181. Springer-Verlag, New
  York-Heidelberg, 1972).

\bibitem{maalqvist2014localization}
A.~Malqvist and D.~Peterseim.
\newblock Localization of elliptic multiscale problems.
\newblock {\em Math. Comp.}, 83(290):2583--2603, 2014.

\bibitem{mclean2000strongly}
W.~McLean.
\newblock {\em Strongly elliptic systems and boundary integral equations}.
\newblock Cambridge University Press, 2000.

\bibitem{melenk2003hp_article}
J.M. Melenk.
\newblock {$hp$}-interpolation of nonsmooth functions and an application to
  $hp$-a posteriori error estimation.
\newblock {\em SIAM J. Numer. Anal.}, 43(1):127--155, 2005.

\bibitem{quva:99}
A.~Quarteroni and A.~Valli.
\newblock {\em Domain Decomposition Methods for Partial Differential
  Equations}.
\newblock Numerical Mathematics and Scientific Computation. Oxford University
  Press, Oxford, UK, first edition, 1999.

\bibitem{SauterSchwab2004}
S.~Sauter and C.~Schwab.
\newblock {\em Boundary element methods}, volume~39 of {\em Springer {S}eries
  in {C}omputational {M}athematics}.
\newblock Springer, 2010.

\bibitem{savare98}
G.~Savar\'e.
\newblock Regularity results for elliptic equations in {L}ipschitz domains.
\newblock {\em Journal of Functional Analysis}, 152:176--201, 1998.

\bibitem{book_tartar}
L.~Tartar.
\newblock {\em An {I}ntroduction to {S}obolev {S}paces and {I}nterpolation
  {S}paces}, volume~3 of {\em Lecture Notes of the Unione Matematica Italiana}.
\newblock Berlin, Springer-Verlag, 2007.

\bibitem{Toselli2005}
A.~Toselli and O.~Widlund.
\newblock {\em Domain decomposition methods -- algorithms and theory},
  volume~34 of {\em Springer {S}eries in {C}omputational {M}athematics}.
\newblock Springer, 2005.

\end{thebibliography}

\end{document}